\documentclass[12pt]{amsart}
\usepackage{amsmath}
\usepackage{amssymb}
\usepackage{amscd}
\usepackage{amsthm}
\usepackage{epsfig}
\usepackage{color}
\numberwithin{equation}{section}
\theoremstyle{plain}
\newtheorem{theorem}[equation]{Theorem}
\newtheorem{claim}[equation]{Claim}
\newtheorem{thm}[equation]{Theorem}

\newtheorem{prop}[equation]{Proposition}
\newtheorem{problem}[equation]{Problem}
\newtheorem{assumption}[equation]{Assumption}

\newtheorem{corollary}[equation]{Corollary}
\newtheorem{conjecture}[equation]{Conjecture}

\theoremstyle{remark}
\newtheorem{remark}[equation]{Remark}
\newtheorem{rem}[equation]{Remark}
\theoremstyle{definition}
\newtheorem{definition}[equation]{Definition}
\newtheorem{defn}[equation]{Definition}

\newtheorem{question}[equation]{Question}
\newtheorem{observation}[equation]{Observation}

\newtheorem{example}[equation]{Example}

\def\C{\mathbb C}
\def\R{\mathbb R}
\def\H{\mathbb H}
\def\Z{\mathbb Z}

\def\S{\mathbb S}
\def\B{\mathbb B}

\def\N{\mathbb N}

\newcommand{\acts}{\curvearrowright}

\newcommand{\Isom}{\operatorname{Isom}}

\newcommand{\mob}{\operatorname{{\mathbf{Mob}}(\S^n)}}
\newcommand{\MOB}{\operatorname{{\mathbf{Mob}}(\S^{n+1})}}
\newcommand{\Mob}{\operatorname{{\mathbf{Mob}}}}
\newcommand{\mbd}{\operatorname{{\mathbf{Mob}}(\S^2)}}
\newcommand{\mbt}{\operatorname{{\mathbf{Mob}}(\S^3)}}
\newcommand{\mbc}{\operatorname{{\mathbf{Mob}}(\S^4)}}
\newcommand{\mbp}{\operatorname{{\mathbf{Mob}}(\S^5)}}

\renewcommand{\int}{\operatorname{int}}
\newcommand{\ext}{\operatorname{ext}}
\newcommand{\diam}{\operatorname{diam}}
\newcommand{\sproof}{{\em Sketch of the proof:}}
\newcommand{\length}{\operatorname{length}}

\def\D{\partial}
\newcommand{\al}{\alpha}

\def\Del{\Delta}
\def\eps{\epsilon}
\def\ga{\gamma}
\def\Ga{\Gamma}

\def\la{\lambda}
\def\La{\Lambda}

\def\<{\langle}
\def\>{\rangle}

\def\Om{\Omega}

\def\ra{\rightarrow}

\def\si{\sigma}
\def\Si{\Sigma}

\def\be{\beta}
\def\del{\delta}

\def\t{\tilde}

\def\hook{\hookrightarrow}

\newcommand{\qf}{quasifuchsian~}

\begin{document}

\title{Kleinian groups in higher dimensions}

\author{Michael Kapovich}

\date{\today}

\begin{abstract}
This is a survey of higher-dimensional Kleinian groups, i.e., discrete isometry groups of the
hyperbolic $n$-space $\H^n$ for $n\ge 4$. Our main emphasis is on the topological and geometric aspects of
higher-dimensional Kleinian groups and their contrast with the discrete groups of isometry of $\H^3$.
\end{abstract}

\address{Department of Mathematics, University of California,
Davis, CA 95616} \email{kapovich@math.ucdavis.edu}

\maketitle

\centerline{{\em To the memory of Sasha Reznikov}}

\section{Introduction}

The goal of this survey is to give an overview (mainly from the topological perspective) of the theory of Kleinian groups in higher
dimensions. The survey grew out of a series of lectures I gave in the University of Maryland in the Fall of 1991.
An early (much shorter) version of this paper appeared as the preprint \cite{Kapovich(1992)}. In this survey
I collect well-known facts as well as less-known and new results. Hopefully, this will
make the survey interesting to both non-experts and experts. We also refer the reader to Tukia's short survey \cite{Tukia(1995)}
 of higher-dimensional Kleinian groups.

There is a vast variety of Kleinian groups in higher dimensions: It appears that there is no hope for
a comprehensive structure theory similar to the theory of discrete groups of isometries of $\H^3$.
I do not know a good guiding principle for the taxonomy of higher-dimensional Kleinian groups.
In this paper the higher-dimensional Kleinian groups are organized according to the
topological complexity of their limit sets.
In this setting one of the key questions that I will address is the interaction between
the geometry and topology of the limit set and the algebraic and topological properties
of the Kleinian group.

This paper is organized as follows. In Section \ref{basics} we consider the most basic concepts of the theory of Kleinian groups, e.g.
domain of discontinuity, limit set, geometric finiteness, etc.
In Section \ref{waysandmeans} we discuss various ways to construct Kleinian groups
and list the tools of the theory of Kleinian groups in higher dimensions.
In Section \ref{homology} we review the homological algebra
used in the paper. %; most readers should skip this section during the first reading of the survey.
In Section \ref{toprigidity} we state topological rigidity results of Farrell and Jones and the coarse compact core theorem
for higher-dimensional Kleinian groups. In Section \ref{equivalences} we discuss
various notions of equivalence between Kleinian groups: From the weakest (isomorphism) to the strongest (conjugacy).
In Section \ref{0} we consider groups with zero-dimensional limit sets; such groups are relatively well-understood.
Convex-cocompact groups with 1-dimensional limit sets are discussed in Section \ref{kk}. Although the topology
of the limit sets of such groups is well-understood, their group-theoretic structure is a mystery.
We know very little about Kleinian groups with higher-dimensional limit sets, thus we restrict
the discussion to Kleinian groups whose limit sets are topological spheres (Section \ref{spheres}). We then discuss
Ahlfors finiteness theorem and its failure in higher dimensions (Section \ref{grail}).
 We then consider the representation varieties of Kleinian groups (Section \ref{varieties}).
 Lastly  we discuss algebraic and  topological constraints on Kleinian groups in higher dimensions (Section \ref{constrains}).

\medskip
{\bf Acknowledgments.} During this work I was partially supported
by  various NSF grants, especially DMS-8902619 at the University
of Maryland and  DMS-04-05180 at UC Davis. Most of this work was
done when I was visiting the Max Plank Institute for Mathematics
in Bonn. I am grateful to C.~McMullen, T.~Delzant, A.~Nabutovsky
and J.~Souto for several suggestions, and to L.~Potyagailo for a
number of comments, suggestions and corrections. I am also
grateful to the referee for numerous corrections.

\tableofcontents

\section{Basic definitions}
 \label{basics}

{\bf M\"obius transformations.} For the lack of space, our
discussion of the basics of Kleinian groups below is somewhat
sketchy. For the detailed treatment we refer the reader to
\cite{Benedetti-Petronio, Bowditch(1993b), Kapovich00, KAG,
Ratcliffe(1994)}. We let $\B^{n+1}$ denote the closed ball
$\H^{n+1}\cup \S^n$; its boundary $\S^n$ is identified via the
stereographic projection with $\overline{\R^{n}} = {\R}^{n}\cup
\{\infty \}$. A {\em horoball} $B$ in $\H^{n+1}$ is a round ball
in $\H^{n+1}$ which is tangent to the boundary sphere $\S^n$. The
point of tangency is called the (hyperbolic) {\em center} of $B$.

Let $\mob$ denote the group of all
{\em M\"obius} transformations of the $n$-sphere
${\S}^{n}$, i.e., compositions of inversions in $\S^n$. The group $\mob$ admits an extension to the hyperbolic space  $\H^{n+1}$,
so that $\mob=\Isom(\H^{n+1})$, the isometry group of $\H^{n+1}$.

For elements $\ga\in \mob$ define the {\em displacement function}
$$
d_\ga (x):= d(x, \ga(x)), \quad x\in \H^{n+1}.
$$
The elements $\ga$ of $\mob$ are classified as:

1. {\em Hyperbolic:} The function $d_\ga$ is bounded away from zero. Its minimum is attained on a geodesic $A_\ga\subset \H^{n+1}$
invariant under $\ga$. The ideal end-points of $A_{\ga}$ are the fixed points of $\ga$ in $\S^n$.

2. {\em Parabolic}: The  function $d_\ga$ is positive but has zero infimum on $\H^{n+1}$; such elements have precisely one fixed
point in $\S^n$.

3. {\em Elliptic}: $\ga$ fixes a point in $\H^{n+1}$.

\medskip
The group $\mob$ is isomorphic to an index 2 subgroup in the Lorentz
group $O(n+1,1)$, see e.g. \cite{Ratcliffe(1994)}.
In particular, $\mob$ is a matrix group. Selberg's lem\-ma \cite{Selberg} implies that every
finitely generated group of matrices contains a finite index subgroup which is torsion-free.
A group $\Ga$ is said to {\em virtually} satisfy a property $X$ if it contains a finite index subgroup $\Ga'\subset \Ga$,
such that $\Ga'$ satisfies $X$. Therefore, every  finitely generated group of matrices is {\em virtually torsion-free}.
Moreover, every finitely-generated matrix group is {\em residually finite},
i.e., the intersection of all its finite-index subgroups is trivial, see  \cite{Malcev(1940), Selberg}.
This, of course, applies to the finitely generated subgroups of $\mob$ as well.

\begin{defn}
A discrete subgroup $\Ga\subset \mob$ is called a {\em Kleinian group}.
 \end{defn}

{\bf Dynamical notions.}
The {\em discontinuity set} $\Omega(\Gamma)$ of a group $\Gamma \subset \mob$,
is the largest open subset in $\S^n$ where $\Ga$ acts properly discontinuously.
Its complement $\S^n\setminus \Omega(\Gamma)$ is the {\em limit set} $\La(\Ga)$ of the group $\Ga$. Equivalently,
the limit set of a Kleinian group can be described as the accumulation set in the sphere $\S^n$ of an orbit $\Ga\cdot o$.
Here $o$ is an arbitrary point in $\H^{n+1}$. A Kleinian group is called {\em elementary} if its limit set is finite,
i.e., is either empty, or consists of one or of two points.

We will use the notation $M^n(\Ga)$ for the $n$-dimensional quotient $\Omega(\Ga)/\Ga$ and $\bar{M}^{n+1}(\Ga)$
for the $n+1$-dimensional quotient $(\H^{n+1}\cup \Om(\Gamma))/\Gamma$.

For a closed subset $\La\subset \S^n$, let $Hull(\La)$ denote its {\em convex hull} in $\H^n$, i.e.,
the smallest closed convex subset $H$ of $\H^n$ such that
$$
cl_{\B^{n+1}} (H)\cap \S^n= \La.
$$
Clearly, if $\La$ is a point, then $Hull(\La)$ does not exist.
Otherwise, $Hull(\La)$ exists and is unique.
%see for instance \cite{Bowditch(1993b), Ratcliffe(1994)}.
We declare $Hull(\La)$ to be empty in the case when $\La$ is a single point.

One way to visualize the convex hull $Hull(\La)$ is to consider the {\em projective model} of the
hyperbolic space, where the geodesic lines are straight line segments contained in the interior of $\B^{n+1}$.
 Therefore, the Euclidean notion of convexity coincides with
the hyperbolic notion. This implies that the convex hull in this model can be described as follows:
$Hull(\La)$ is the intersection of the  Euclidean convex hull of $\La$ with the interior of  $\B^{n+1}$.

Suppose that $\La=\La(\Ga)$ is the limit set of a Kleinian group $\Ga\subset \mob$.
The quotient $Hull(\La)/\Ga$ is called the {\em convex core} of the orbifold $N=\H^{n+1}/\Ga$.
It is characterized by the property that it is the smallest closed convex
subset in $N$, whose inclusion to $N$ is a homotopy-equivalence.
For $\epsilon > 0$ consider the open $\epsilon$-neighborhood
$Hull_{\epsilon}(\La)$ of $Hull(\La)$ in $\H^{n+1}$. Since $Hull_{\epsilon}(\La)$ is $\Ga$-invariant,
we can form the quotient  $Hull_{\epsilon}(\La)/\Ga$.  Then $Hull_{\epsilon}(\La)/\Ga$
is the $\eps$-neighborhood of the convex core.

{\bf Geometric finiteness.} We now arrive to one of the key notions in the theory of Kleinian groups:

\begin{defn}
A Kleinian group $\Ga\subset \mob$ is called {\em geometrically finite} if:

(1) $\Ga$ is finitely generated, and

\smallskip
(2) $vol(Hull_\eps(\La(\Ga))/\Ga) < \infty$.
\end{defn}

In a number of important special cases, e.g. when $\Ga$ is torsion-free, or $n=2$, or when $\La(\Ga)=\S^n$,
the assumption (1) follows from (2), see \cite{Bowditch(1993b)}. However,
E.\! Hamilton \cite{Ham} constructed an example of a Kleinian
group $\Ga\subset \Mob(\S^3)$ for which (2) holds but (1) fails. This group contains finite order elements
of arbitrarily high order; by Selberg's lemma such groups cannot be finitely-generated.

A Kleinian group $\Ga\subset \mob$ is called a {\em lattice}
if $\H^{n+1}/\Ga$ has finite volume. Equivalently, $\La(\Ga)=\S^n$ and $\Ga$ is geometrically finite.
A lattice is {\em cocompact } (or {\em uniform}) if $\H^{n+1}/\Ga$ is compact.
%It is immediate that for every lattice $\Ga\subset \mob$, $\Om(\Ga)=\emptyset$.

One can characterize geometrically finite groups in terms of their limit sets. Before stating this theorem we
need two more definitions.

\begin{defn}
A limit point $\xi\in \La(\Ga)$ is called a {\em conical limit point} if there exists a geodesic
$\al\subset \H^n$ asymptotic to $\xi$, a point $o\in \H^n$, a number $r<\infty$, and a sequence $\ga_i\in \Ga$
so that

1. $\lim_i \ga_i(o)=\xi$.

2. $d(\ga_i(o), \al)\le r$.
\end{defn}

The reason for this name comes from the shape of the $r$-neighborhood of the vertical geodesic  $\al$ in the
upper half-space model of $\H^{n+1}$: It is a Euclidean cone with the axis $\al$. Equivalently, one can describe
the conical limit points of nonelementary groups as follows (see \cite{Beardon-Maskit, Bowditch(1993b)}):

\smallskip
$\xi \in \Lambda (\Ga)$ is a conical limit point if and only if for every
$\eta \in \Lambda (\Ga) \setminus \{ \xi \}$ there exists a point $\psi$ and a sequence $\ga_i \in \Ga$ such
that:

1. $\lim_{i} \ga_i (\zeta) = \xi$ for every $\zeta \in \Lambda (\Ga)\setminus \{\psi\}$.

2. $\lim_i \ga_i^{-1}(\xi)\ne \lim_i \ga_i^{-1}(\eta)$.

\noindent The set of conical limit points of a Kleinian group $\Ga$ is denoted $\La_c(\Ga)$.

\begin{defn}
A point $\xi \in \Lambda (\Ga)$ is called a {\em bounded parabolic point}
if it is the fixed point of a parabolic subgroup $\Pi \subset \Ga$ and
$(\Lambda (\Ga) - \{ \xi\})/\Pi$ is compact.
\end{defn}

Below is a {\em dynamical} characterizations of geometrically finite groups:

\begin{theorem}
\label{BiMa}
(A.\! Beardon and B.\! Maskit \cite{Beardon-Maskit}, B.\! Bowditch \cite{Bowditch(1993b)})
A Kleinian group $\Ga$ is geometrically finite if and only if each limit point $\xi\in \La(\Ga)$ is
either a {\em conical limit point} or a {\em bounded parabolic point}.
\end{theorem}

C.\! Bishop \cite{Bishop(1996)} proved  that one can drop the word {\em bounded} in the above theorem.
We refer the reader to Bowditch's paper \cite{Bowditch(1993b)} for the proof of other
criteria of geometric finiteness collected in Theorems \ref{T1}, \ref{T2}, \ref{T3} below.
(The case $n=2$ is treated in \cite{Marden} and \cite{Thurston(1978-81)}.)

\begin{thm}\label{T1}
1. If a  Kleinian subgroup $\Ga\subset \mob$ admits a convex fundamental
polyhedron with finitely many faces then it is geometrically finite.

2. Let $\Ga\subset \mob$ be a geometrically finite Kleinian group so that either
(a) $n\le 2$, or (b) $\Ga$ contains no parabolic elements,
or (c) $\Ga$ is a lattice.

Then $\Ga$ admits a convex fundamental polyhedron with finitely many faces.
\end{thm}

On the other hand, there are geometrically finite subgroups of $\mbt$
which do not admit a convex fundamental polyhedron with finitely many faces, see \cite{Apanasov}.

\begin{thm}\label{T2}
Let $\Ga\subset \mob$ be a Kleinian subgroup containing no parabolic elements. Then the following are equivalent:

(a) $\Ga$ is geometrically finite.

(b) $Hull(\La(\Ga))/\Ga$ is compact.

(c) $\bar{M}^{n+1}(\Ga)$ is compact.
\end{thm}

If $\Ga$ is geometrically finite and contains no parabolic elements, it is called {\em convex--cocompact}.
We will frequently use the fact that every convex-cocompact Kleinian group is {\em Gromov-hyperbolic},
 see e.g. \cite{Bridson-Haefliger}.

The criterion given in Theorem \ref{T2} generalizes to the case of groups with parabolic elements, although the
statement becomes more complicated:

\begin{thm}\label{T3} The following are equivalent:

(a)  $\Ga$ is geometrically finite.

(b) There exists a pairwise disjoint $\Ga$-invariant collection of open horoballs $B_i\subset \H^{n+1}, i\in I$,
which are centered at fixed points of parabolic subgroups of $\Ga$, such that the quotient
$$
\left(Hull(\La(\Ga))\setminus \bigcup_{i\in I} B_i\right)/\Ga
$$
is compact.

(c) Let $\Pi_i, i\in I$ be the collection of maximal (virtually) parabolic subgroups of $\Ga$. For
each $i$ there exists a $\Pi_i$--invariant convex subset $C_i\subset \B^{n+1}$, so that
the quotient
$$
\left(\H^{n+1}\cup \Om(\Ga)\setminus \bigcup_{i\in I} C_i\right)/\Ga
$$
is compact. If $\Om(\Ga)=\emptyset$, then one can take $C_i=B_i$, a horoball in $\H^{n+1}$.
\end{thm}

If $n=1$, then every finitely generated   Kleinian group is geometrically finite. The proof is rather elementary,
see e.g. \cite{Casson-Bleiler(1988)}. For $n\ge 2$ this  implication is no longer true. The  first (implicit)
examples were given by L.~Bers, they are {\em singly-degenerate groups}:

\begin{defn}
A finitely generated nonelementary Kleinian subgroup of $\mbd$ is {\em singly degenerate} if its domain of discontinuity is
simply-connected, i.e., homeomorphic to the 2-disk.
\end{defn}

L.~Bers \cite{Bers(1970)} proved that singly degenerate Kleinian
groups exist and are never geometrically finite. The first {\em
explicit} examples of finitely generated geometrically infinite
Kleinian subgroups $\Ga$ of $\mbd$ were given by T.~J{\o}rgensen
\cite{Jorgensen(1977)}. In J{\o}rgensen's examples, $\Ga$ appears
as a normal subgroup of a lattice $\hat{\Ga}\subset \mbd$ with
$\hat{\Ga}/\Ga\cong \Z$. Remarkably, all known examples of
finitely-generated geometrically infinite Kleinian subgroups of
$\mob$ can be traced to the 2-dimensional examples. More
precisely, every {\em known} finitely-generated geometrically
infinite Kleinian subgroup $\Ga\subset \mob$ admits a
decomposition as the graph of groups
$$
({\mathcal G}, \Ga_v, \Ga_e),
$$
where at least one of the vertex groups $\Ga_v$ is either a geometrically infinite subgroup contained in $\mbd$,
or is a quasiconformal deformations of such.

\begin{problem}
Construct examples of finitely-generated geometrically infinite subgroups of $\mob$, $n\ge 3$, which
do not have the 2-dimensional origin as above.
\end{problem}

%The best candidates for such examples seem to be among subgroups of the convex-cocompact groups constructed in
%\cite{Kapovich05}.

\begin{assumption}
From now on we will assume that all Kleinian groups are
finitely generated and torsion free, unless stated otherwise.
\end{assumption}

Note that the second part of this assumption is not very restrictive because of Selberg's lemma.

\medskip
{\bf Cusps and tubes.} The $\Ga$-conjugacy classes $[\Pi]$ of maximal parabolic
subgroups $\Pi$ of a Kleinian group $\Ga$ are called {\em cusps} of $\Ga$. More geometrically, cusps of $\Ga$ can be described
using the {\em thick-thin decomposition} of the quotient manifold $M=\H^{n+1}/\Ga$. Given a positive
number $\eps>0$, let $M_{(0,\eps]}$ denote the collection of points $x$ in $M$ such that  there exists
a homotopically nontrivial loop $\al$ based at $x$, so that the length of $\al$ is at most $\eps$.
Then $M_{(\eps,\infty)}$ is the complement of $M_{(0,\eps]}$ in $M$. According to Kazhdan-Margulis lemma
\cite{Kazhdan-Margulis}, there exists a number $\mu=\mu_{n+1}>0$ such that for every Kleinian group $\Ga$
and every $0<\eps\le \mu$, every component
of $M_{(0,\eps]}$  has a {\em virtually abelian} fundamental group.
The submanifold $M_{(0,\eps]}$ is called the {\em thin part} of $M$ and its complement the {\em thick part} of $M$.
The compact components of $M_{(0,\eps]}$ are called {\em tubes} and the noncompact components are called {\em cusps}.

Then the cusps of $\Ga$ are in bijective correspondence with the cusps in  $M_{(0,\eps]}$:

For every cusp $[\Pi]$ in $\Ga$, there exists a noncompact component $C\subset  M_{(0,\eps]}$, so that $\Pi=\pi_1(C)$.
Conversely, for each cusp $C\subset M$, there exists a maximal parabolic subgroup $\Pi\subset \Ga$ such that $\Pi=\pi_1(C)$.

\smallskip
\noindent Taking the $\Ga$-conjugacy class of $\Pi$ reflects the ambiguity in the choice of the base-point needed to
identify $\pi_1(M)$ and $\Ga$.

If $n\le 2$ and the manifold $M$ is oriented, then the components $C_i$ of   $M_{(0,\eps]}$ are convex:
The cusps in $M$ are quotients of horoballs in $\H^{n+1}$, while the compact components $T_i$ of $M_{(0,\eps]}$
are metric $R_i$-neighborhoods of closed geodesics $\ga_i\subset M$.
In higher dimensions ($n\ge 3$) convexity (in general) fails. However every tube $T_i$ in $M_{(0,\eps]}$ is a finite union
of convex sets containing a certain closed geodesic $\ga_i\subset T_i$.
In particular, every tube $T_i$ is homeomorphic to a disk bundle over $\S^1$.
A similar, although more complicated description, holds for the cusps, where one
has to consider (in general) a union of infinitely many convex subsets. See e.g. \cite{Kapovich07}.

 \medskip
{\bf M\"obius structures.} In this paper we shall also discuss the subject closely related
to the theory of Kleinian groups, namely {\em M\"obius structures}.
When $M$ is a smooth manifold of dimension $\ge 3$, a {\em M\"obius} (or {\em flat conformal}) structure $K$ on a
 $M$ is the conformal class of a conformally-Euclidean Riemannian metric on $M$.
Topologically, $K$ is a maximal {\em M\"obius atlas} on
$M$, i.e., an atlas  with M\"obius transition maps. Thus, for each Kleinian group $\Ga$ and $\Ga$-invariant subset $\Om\subset \Om(\Ga)$,
the standard M\"obius structure on $\Om\subset \S^n$ projects to a M\"obius structure
$K_{\Ga}$ on the manifold $\Omega/\Ga$. The M\"obius structures of this type are called
{\em uniformizable}.

\medskip
{\bf Complex-hyperbolic Kleinian groups.} Instead of considering
the isometry group of the hyperbolic space, one can consider other negatively curved symmetric spaces, for instance,
the {\em complex}--{\em hyperbolic} $n$-{\em space} $\C \H^n$ and its group of automorphisms
$PU(n,1)$. From the analytical viewpoint, $\C \H^n$ is the unit ball in $\C^n$ and $PU(n,1)$ is the group
of biholomorphic automorphisms of this ball.
The Bergman metric on $\C \H^n$ is a K\"ahler metric of negative sectional curvature. The discrete subgroups of $PU(n,1)$ are
{\em complex--hyperbolic Kleinian groups}. They share many properties with Kleinian groups. In fact,
{\em nearly all positive results} stated in this survey for Kleinian subgroups of $\mob$ ($n\ge 3$) are also valid
for the  complex--hyperbolic Kleinian groups! (One has to replace {\em virtually abelian} with {\em virtually nilpotent}
in the discussion of cusps.) There exists an isometric embedding $\H^n\to \C \H^n$ which induces an embedding of the isometry groups
and  therefore complex--hyperbolic Kleinian groups ($n\ge 4$) also inherit the pathologies of the higher--dimensional
Kleinian groups. We refer the reader to \cite{Goldman(1992), Goldman(1998), Schwartz(2002)} for detailed discussion.

\section{Ways and means of Kleinian groups}
\label{waysandmeans}

\subsection{Ways: Sources of Kleinian groups}

The following is a list of ways to construct Kleinian groups.

\medskip
(a) {\bf Poincar\'e fundamental polyhedron theorem} (see e.g. \cite{Ratcliffe(1994)} for a very detailed discussion, as well
as \cite{Maskit(1987)}). This source is, in principle, the most general. The Poincar\'e
fundamental polyhedron theorem  asserts that given a polyhedron $\Phi$ in $\H^{n+1}$
and a collection of elements $\ga_1 , \ga_2, ..., \ga_k ,...$ of $\mob$, pairing the faces of
$\Phi$, under certain conditions on this data,
the group $\Ga$ generated by $\ga_1 , \ga_2, ..., \ga_k ,...$ is Kleinian and
$\Phi$ is a {\em fundamental domain} for the action of the group $\Ga$ on $\H^{n+1}$.

Every Kleinian group has a convex fundamental
polyhedron (for example, the {\em Dirichlet fundamental domain}). However,
in practice, the Poincar\'e fundamental polyhedron theorem is not always easy to use, especially if
$\Phi$ has many faces and $n$ is large. This theorem was used, for instance,
to construct non-arithmetic lattices in $\mob$ (see \cite{Makarov66, Makarov68, Vinberg(1967)}),
as well as other interesting Kleinian groups, see e.g.  \cite{Davis(1985), IRT, Kapovich(1993b), Kuiper,   Ratcliffe-Tschantz}.

\smallskip
(b) {\bf Klein--Maskit Combination Theorems} (see e.g. \cite{KAG} and \cite{Maskit(1987)}).
Suppose that we are given two Kleinian groups $\Ga_1, \Ga_2\subset \mob$ which share a common subgroup $\Ga_3$,
or a single Kleinian group $\Ga_1$ and a M\"obius transformation $\tau\in \mob$ which conjugates subgroups
$\Ga_3, \Ga_3'\subset \Ga_1$. The Combination Theorems provide conditions
which guarantee that the group $\Ga\subset \mob$ generated by $\Ga_1$ and $\Ga_2$ (or by $\Ga_1$ and $\tau$)
is again Kleinian and is isomorphic to the amalgam
$$
\Ga\cong \Ga_1*_{\Ga_3} \Ga_2,
$$
or to the HNN extension
$$
\Ga\cong \Ga_1*_{\Ga_3}= HNN(\Ga_1, \tau).
$$
The proofs of the Combination Theorems generalize the classical ``ping-pong'' argument due to Schottky and Klein.
The Combination Theorems also show that the quotient manifold $M^n(\Ga)$
of the group $\Ga$ is obtained from  $M^n(\Ga_1)$, $M^n(\Ga_2)$ (or $M(\Ga_1)$)
via some ``cut-and-paste" operation. Moreover, Combination Theorems generalize to
graph of groups. There should be a generalization of Combination Theorems to
{\em complexes of groups} (see e.g. \cite{Bridson-Haefliger}
for the definition); however, to the best of my knowledge,
nobody worked out the general result, see \cite{Kapovich05} for a special case.

\smallskip
(c) {\bf Arithmetic groups and their subgroups} (see e.g. \cite{Maclachlan-Reid} and \cite{Vinberg-Shvartsman}).
A subgroup $\Ga\subset O(n,1)$ is called {\em arithmetic} if there exists an embedding
$$
\iota: O(n,1) \hook GL(N,\R),
$$
such that the image $\iota(\Ga)$ is {\em commensurable} with the intersection
$$
\iota(O(n,1))\cap GL(N,\Z).
$$
Recall that two subgroups $\Ga_1, \Ga_2\subset G$ are called {\em commensurable} if
$\Ga_1\cap \Ga_2$ has finite index in both $\Ga_1$ and $\Ga_2$.

\medskip
Below is a specific construction of arithmetic groups.
Let $f$ be a quadratic form of signature $(n,1)$ in $n+1$
variables with coefficients in a totally real algebraic number field
$K\subset \R$ satisfying the following condition:

(*) For every nontrivial (i.e., different from the identity)
embedding $\si: K\to \R$, the quadratic form $f^\si$ is
positive definite.

\smallskip
Without loss of generality one may assume that this quadratic form is diagonal. For instance, take
$$
f(x)= -\sqrt{2}x_0^2+ x_1^2+...+x_n^2.
$$

We now define discrete subgroups of $\Isom(\H^n)$ using the form $f$.
Let $A$ denote the ring of integers of $K$. We define the group $\Gamma:=O(f, A)$ consisting of matrices
with entries in $A$ preserving the form $f$. Then $\Gamma$ is
a discrete subgroup of $O(f, \R)$. Moreover, it is a lattice:
its index 2 subgroup
$$
\Gamma'=O'(f, A):=O(f, A)\cap O'(f, \R)$$
acts on $\H^n$ so that $\H^n/\Gamma'$ has finite volume.
Such groups $\Gamma$ (and subgroups of $\Isom(\H^n)$ commensurable to them) are called {\em  arithmetic
subgroups of the simplest type} in $O(n, 1)$, see \cite{Vinberg-Shvartsman}.

\begin{rem}
If $\Gamma\subset O(n, 1)$ is an arithmetic lattice so that
either $\Gamma$ is non-cocompact or $n$ is even, then it follows from
the classification of rational structures on $O(n,1)$ that $\Gamma$
is commensurable to an arithmetic lattice of the simplest type. For odd $n$ there
is another family of arithmetic lattices given as
the groups  of units of appropriate skew-Hermitian forms over
quaternionic algebras. Yet other families of arithmetic lattices exist for $n=3$ and $n=7$. See e.g.
\cite{Vinberg-Shvartsman}.
\end{rem}

We refer the reader to \cite{Maclachlan-Reid} for the detailed treatment of geometry and topology of arithmetic subgroups
of $\mbd$.

\medskip
(d) {\bf Small deformations of a given Kleinian group}. We discuss this construction in detail in Section \ref{local}.
The idea is to take a Kleinian group $\Ga\subset \mob$ and to ``perturb it a little bit'', by modifying the generators
slightly (within $\mob$) and preserving the relators. The result is a new group $\Ga'$ which may or may not be Kleinian
and even if it is, $\Ga'$ is not necessarily isomorphic to $\Ga$.
However if $\Ga$ is convex-cocompact, $\Ga'$ is again a convex-cocompact group isomorphic to $\Ga$, see
Theorem \ref{stabilitythm}.

\medskip (e) {\bf Limits  of  sequences  of  Kleinian groups}, see Section \ref{global}.
Take a sequence $\Ga_i$ of Kleinian subgroups of $\mob$ and assume that it has a limit $\Ga$:
It turns out that there are two ways to make sense of this procedure ({\em algebraic} and {\em geometric} limit).
In any case, $\Ga$ is again a Kleinian group. Even if the (algebraic) limit does not exist
as a subgroup of $\mob$, there is a way to make sense of the limiting group as a group of isometries of a metric tree.
This logic turns out to be useful for  proving compactness theorems for sequences of Kleinian  groups.

\medskip
(f) {\bf Differential-geometric constructions of hyperbolic metrics.} The only (but spectacular) example where it has been
used is Perelman's work on Ricci flow and proof of Thurston's geometrization conjecture.
See \cite{Kleiner-Lott, Morgan-Tian, Perelman1, Perelman2}.
However applicability of this tool at the moment appears to be limited to 3-manifolds.

\medskip
A beautiful example of application of (b) and (c) is
the construction of M.~Gromov and I.~Piatetski-Shapiro \cite{GP} of {\em non-arithmetic}
lattices in $\mob$.  Starting with two arithmetic groups $\Gamma_j$
($j= 1, 2$) they first ``cut these groups in half", take ``one half'' $\Del_j \subset \Gamma_j$ of each,
and then combine $\Del_1$ and $\Del_2$ via Maskit Combination. The construction of Kleinian groups in
\cite{Gromov-Thurston(1987)} (see also Section \ref{codimension1}) is an application of (b), (c) and (d).
Thurston's hyperbolic Dehn surgery theorem is an example of (e).
One of the most sophisticated  constructions of Kleinian groups is given by Thurston's hyperbolization theorem
(see e.g.  \cite{Kapovich00}, \cite{Otal(1996)}, \cite{Otal(1997)});  still, it is essentially a
combination (a very complicated one!) of (b), (d) and (e).

\begin{rem}
There is  potentially the sixth source of Kleinian groups in higher
dimensions: monodromy of linear ordinary differential equations.
However, to the best of my knowledge, the only example of its application relevant to Kleinian groups, is the construction
of lattices in $PU(n,1)$ (i.e., the isometry group of the complex-hyperbolic $n$--space) by Deligne and Mostow,
see \cite{Deligne-Mostow}.
\end{rem}

\subsection{Means: Tools of the theory of Kleinian groups in higher dimensions}

Several key tools of the ``classical'' theory of Kleinian subgroups of $\mbd$ (mainly, the Beltrami equation and pleated hypersurfaces)
are missing in higher dimensions. Below is the list of main tools that are currently available.

\medskip
(a) {\bf Dynamics, more specifically, the convergence property.}
Namely, every sequence of  M\"obius transformations $\ga_i\in \mob$ either contains a convergent subsequence
or contains a subsequence which converges to a constant map away from a point in $\S^n$. See e.g. \cite{Kapovich00}.

\medskip
(b) {\bf Kazhdan-Margulis lemma} and its corollaries.

\smallskip
It turns out that the lion share of the general results about higher-dimensional Kleinian groups is a combination of (a) and (b),
together with some hyperbolic geometry.

\medskip
(c) {\bf Group actions on trees and Rips theory}. This is a very potent tool for proving compactness results for
families of representations of Kleinian groups, see for instance Theorem \ref{morg}.

\medskip
(d) {\bf Barycentric maps}. These maps were originally introduced by A.~Douady and C.~Earle \cite{Douady-Earle}
as a tool of the Teichm\"uller theory of Riemann surfaces. In the hands of G.~Besson, G.~Courtois and S.~Gallot
 these maps became a powerful analytic tool of the theory
of Kleinian groups in higher dimensions, see e.g. \cite{Besson-Courtois-Gallot(1999), BCG03},
as well as Theorems \ref{inequality} and \ref{contract} in this survey.
In contrast, equivariant harmonic maps which proved so useful in the study of, say, K\"ahler groups,
seem at the moment to be only of a very limited use in the theory of Kleinian groups  in higher dimensions.

\medskip
(e) {\bf Ergodic theory of the actions of $\Ga$ on its limit set and Patterson--Sullivan measures.}
See for instance \cite{Nicholls, Sullivan(1981b)} and the survey of P.~Tukia \cite{Tukia(1995)}.

\medskip
(f) {\bf Conformal geometric analysis.} This is a branch of (conformal) differential geometry concerned with
the analysis of the conformally-flat Riemannian metrics on $M^n(\Ga)=\Om(\Ga)/\Ga$. This tool tends to work rather well
in the case when $M^n(\Ga)$ is compact. The most interesting examples of this technique are due to R.~Schoen and S-T. Yau
\cite{Schoen-Yau(1988)},  S.~Nayatani \cite{Nayatani},  A. Chang, J. Qing, J. and P. Yang, \cite{CQY}, and
H.~Izeki \cite{Izeki(1995), Izeki(2000), Izeki(2002)}.

\medskip
(g) {\bf Infinite-dimensional representation theory} of the group $\mob$. The only (but rather striking) example
of its application is Y.~Shalom's work \cite{Shalom}.

\medskip
(h) {\bf Topological rigidity theorems of Farrell and Jones}: See Section \ref{toprigidity}.

\section{A bit of homological algebra}
\label{homology}

{\bf Why does one need homological algebra in order
to study higher-dimen\-sional Kleinian groups?}

Essentially the only time one encounters group cohomology with
twisted coefficients in the study of Kleinian subgroups of $\mbd$,
is in the proof of Ahlfors' finiteness theorem, see \cite{Kra}.
Another (minor) encounter appears in the proof of the smoothness
theorem for deformation spaces of Kleinian groups, Section 8.8 in
\cite{Kapovich00}. Otherwise, homological algebra is hardly ever
needed. The main reason for this, I believe, is 3-fold:

1. Solvability of the 2-dimensional Beltrami equation, which
implies smoothness of the deformation spaces of Kleinian groups in
the most interesting situations.

2. Scott compact core theorem \cite{Scott(1973a), Scott(1973b)}
ensures that every finitely-generated Kleinian group $\Ga\subset
\mbd$ satisfies a very strong finiteness property: Not only it is
finitely-presented, it is also (canonically) isomorphic to the
fundamental group of a compact aspherical 3-manifold with boundary
(Scott compact core).

3. The separation between Kleinian groups of the cohomological
dimension 1, 2 and 3 comes rather easily: Free groups, ``generic''
Kleinian groups, and lattices. Moreover, every Kleinian group
$\Ga\subset \mbd$ which is not a lattice, splits as
\begin{equation}
\label{freepro}
 \Ga\cong \Ga_0 * \Ga_1 * ... * \Ga_k,
\end{equation}
where $\Ga_0$ is free and each $\Ga_i$, $i\ge 1$, is freely
indecomposable, 2-dimensional group. In the language of
homological algebra, the group $\Ga_0$ has cohomological dimension
1, while the groups $\Ga_i$, $i\ge 1$, are {\em two-dimensional
duality groups}.

\medskip
All this changes rather dramatically in higher dimensions:

1. Solvability of the Beltrami equation fails, which, in
particular, leads to non-smoothness of the deformation spaces of
Kleinian groups, Theorem  \ref{jm}. In order to study the local
structure of character varieties one then needs the first and the
second group cohomology with (finite-dimensional) twisted
coefficients.

2. Scott compact core theorem fails for Kleinian subgroups of
$\mbt$, for instance, they do not have to be finitely-presented,
see Section \ref{grail}. Therefore, it appears that one has to
reconsider the assumption that Kleinian groups are
finitely-generated. It is quite likely, that in higher dimensions,
in order to get good structural results, one has to restrict to
Kleinian groups of {\em finite type}, i.e., type $FP$, defined
below. This definition requires homological algebra.

3. One has to learn how to separate $k$-dimensional from
$m$-dimensional in the algebraic structure of Kleinian groups. For
the subgroups of $\mbd$ this separation comes in the form of the
free product decomposition (\ref{freepro}). It appears at the
moment that ``truly'' $m$-dimensional groups are the {\em
$m$-dimensional duality groups}. For instance, for Kleinian
subgroups $\Ga$ of $\mob$ which are $n$-dimensional duality
groups, one can prove a coarse form of the Scott compact core
theorem, \cite{Kapovich-Kleiner(2005)}. In particular, every such
group admits the structure of an $n+1$-dimensional {\em Poincar\'e
duality pair} $(\Ga, \Del)$. The latter is a homological analogue
of the fundamental group of a compact aspherical $n+1$-manifold
with boundary (where the boundary corresponds to the collection of
subgroups $\Del$ in $\Ga$). See Section \ref{core} for more
details.

4. The (co)homological dimension appears to be an integral part of
the discussion of the critical exponent of higher-dimensional
Kleinian groups, see Section \ref{sec:small} and Izeki's papers
\cite{Izeki(1995), Izeki(2000)}.

\medskip
{\bf (Co)homology of groups.}  Some of the above discussion was
rather speculative; we now return to the firm ground of
homological algebra. We refer the reader to \cite{Bieri(1976a),
Brown(1982)} for the comprehensive treatment of (co)homologies of
groups.

Throughout this section we let $R$ be a commutative ring with a
unit. The examples that the reader should have in mind
are $R=\Z$, $\Z/p\Z$ and $R=\R$.
The {\em group ring} $R\Ga$ of a group $\Ga$ consists of finite
linear combinations of the form
$$
\sum_{\ga\in \Ga} r_{\ga} \ga,
$$
with $r_{\ga}\in R$ equal to zero for all but finitely many $\ga\in \Ga$.  Let $V$ be a (left) $R\Ga$-module.
Basic examples include $V=R$ (with the trivial $R\Ga$-module structure)
and $V=R\Ga$. If $R$ is a field, then $V$ is nothing but a vector space over $R$ equipped with a linear
action of the group $\Ga$.  The very useful (for the theory of Kleinian groups) example is the following:

Let $G=\mob$, ${\mathfrak g}$ be the Lie algebra of $G$. Then $G$ acts on ${\mathfrak g}$ via the {\em adjoint representation}
$Ad=Ad_G$. Therefore ${\mathfrak g}$ becomes an $\R G$-module. For every abstract group $\Ga$ and a representation
$\rho: \Ga\to G$ we obtain the $\R\Ga$-module
$$
V={\mathfrak g}_{Ad(\rho)},
$$
where the action of $\Ga$ is given by the composition $Ad\circ \rho$.
We will abbreviate this module to $Ad(\rho)$.
From the theory of Kleinian groups viewpoint,
the most important example of this module is when $\Ga$ is a Kleinian subgroup of $G$ and $\rho$ is the identity embedding.

A {\em projective} $R\Ga$-module, is a module $P$, such that every
exact sequence of $R\Ga$-modules
$$
Q\to P \to 0
$$
splits. For instance, every free $R\Ga$-module is projective.

\medskip
Assume now that $V$ be an $R\Ga$-module.  A {\em resolution of
$V$} is an exact sequence of $R \Ga$-modules:
$$
\cdots\to P_n \to \cdots \to P_0 \to V\to 0.
$$
Every $R\Ga$-module has a unique projective resolution up to a chain
homotopy equivalence.

\begin{example}
Let $V=\Z$, be the trivial $\Z\Ga$-module. Let $K$ be a cell
complex which is $K(\Ga,1)$, i.e., $K$ is connected, $\pi_1(K)\cong
\Ga$ and $\pi_i(K)=0$ for $i\ge 2$. Let $X$ denote the universal
cover of $K$. Lift the cell complex structure from $K$ to $X$. The
group action $\Ga\acts X$, determines a natural structure of a
$\Z\Ga$-module on the cellular chain complex $C_*(X)$. Since the
latter is acyclic, we obtain a resolution of $\Z$ with
$$
P_i=C_i(X),
$$
and the homomorphism $P_0\to \Z$ given by the augmentation.
Moreover, as the group $\Ga$ acts freely on $X$, each module $P_i$
is a {\em free} $\Z\Ga$-module:
$$
P_i\cong \oplus_{j\in C_i} \Z\Ga,
$$
where $C_i$ is the set of $i$-cells in $K$.
\end{example}

A group $\Ga$ is said to be of {\em finite type}, or $FP$ (over
$R$), if there exists a resolution by finitely generated
projective $R\Ga$--modules
$$
0\to P_k\to P_{k-1}\to ... \to P_0\to R\to 0.
$$
For example, if there exists a finite cell complex $K=K(\Ga,1)$,
then $\Ga$ has finite type for an arbitrary ring $R$.  Every group
of finite type is finitely generated, although it does not have to
be finitely-presented, see \cite{Bestvina-Brady}.

The {\em cohomology} of $\Ga$ with coefficients in an
$R\Ga$-module $V$, $H^*(\Ga, V)$, is defined as the homology of
chain complex
$$
{\rm Hom}_{R \Ga}(P_*, M),
$$
where $P_*$ is a projective resolution of the trivial
$R\Ga$-module $R$.  The {\em homology} of $\Ga$ with coefficients
in $V$, $H_*(\Ga, V)$, is the homology of the chain complex
$$
P_*\otimes_{R \Ga} V.$$

An example to keep in mind is the following. Suppose that $K$ is a
manifold, or, more generally, a cell complex, which is an
Eilenberg-MacLane space $K(\Ga, 1)$. Then one can use the chain
complex $C_*(X, R)$ as the resolution $P_*$. Therefore, for the
trivial $\Ga$-module $R$ we have
$$
H^*(\Ga, R)\cong H^*(K, R), \quad H_*(\Ga, R)\cong H_*(K, R).
$$
For the more general modules $V$, in order to compute $H^*(\Ga,
V)$ and $H_*(\Ga, V)$, one uses the (co)homology of $K$ with
coefficients in an appropriate bundle over $K$.

Similarly, given a collection $\Pi$ of subgroups of $\Ga$, one
defines the relative (co)homolo\-gy groups $H^*(\Ga, \Pi; V)$ and $H_*(\Ga,
\Pi; V)$. Whenever discussing (co)homology with $R=\Z$ and $\Z$ as
the trivial $\Z\Ga$-module, we will use the notation $H^*(\Ga),
H_*(\Ga)$.

The (co)homology of groups behaves in a manner similar to the more familiar (co)homology of cell complexes. For instance,
if $\Ga$ admits an $n$-dimensional $K(\Ga,1)$, then $H^i(\Ga,V)=H_i(\Ga, V)=0$ for all $i>n$ and all $R\Ga$-modules.

\smallskip {\bf (Co)homological dimension.} For a group $\Ga$, let
$cd_R(\Ga)$ and $hd_R(\Ga)$ denote the {\em cohomological} and
{\em homological}  dimensions of $\Ga$ (over $R$):
$$
cd_R(\Ga)=\sup \{n: \exists \hbox{~an~} R\Ga\hbox{--module~} V \hbox{~so that~} H^n(\Ga, V)\ne 0\},
$$
$$
hd_R(\Ga)=\sup \{n: \exists \hbox{~an~} R\Ga\hbox{--module~} V \hbox{~so that~} H_n(\Ga, V)\ne 0\}.
$$
We will omit the subscript $\Z$ whenever $R=\Z$.

\medskip
Using the relative (co)homology one defines the relative
(co)homological dimension of $\Ga$ with respect to a collection
$\Pi$ of its subgroups, $cd_R(\Ga, \Pi)$ and $hd_R(\Ga, \Pi)$.

We will use this definition in the case when $\Ga$ is a Kleinian
group as follows. Let ${\mathcal P}$ denote the set of all maximal
(elementary) subgroups of $\Ga$ which contain $\Z^2$. For every
$\Ga$-conjugacy class $[\Pi_i]$ in ${\mathcal P}$, choose a
representative $\Pi_i\subset \Ga$. Then $\Pi$ will denote the set
of all these representatives $\Pi_i$. By abusing the notation, we
will refer to the set $\Pi$ as the set of cusps of virtual rank
$\ge 2$ in $\Ga$.

\medskip
If $\Ga$ is of type $FP$, then
$$
hd_{R}(\Ga)= cd_R (\Ga), \quad \forall \hbox{~~rings~~} R,
$$
see for instance \cite{Bieri(1976a)}. In general,
$$
hd_{R}(\Ga)\le cd_R (\Ga)\le hd_R(\Ga)+1.
$$

\begin{example}
Let $\Ga$ be a free group of finite rank $k>0$. Then $hd_{R}(\Ga)= cd_R (\Ga)$ for all rings $R$.
Indeed, $\Ga$ admits a finite $K(\Ga,1)$ which is the bouquet $B$ of $k$ circles. Since $B$ is 1-dimensional,
$$
hd_{R}(\Ga)= cd_R (\Ga)\le 1.
$$
On the other hand, by taking the trivial $R\Ga$-module $V=R$ we obtain
$$
H_1(B, R)= R^k,
$$
the direct sum of $k$ copies of $R$, and hence is nontrivial.
\end{example}

It turns out that the converse to this example is also true, which is an application  of the famous theorem of J.~Stallings
on the ends of groups:

\begin{thm}\label{stalling}
(J.~Stallings \cite{Stallings(1968)}.) If $\Ga$ is a finitely generated group with $cd(\Ga)=1$, then $\Ga$ is free.
\end{thm}

This result was generalized by M.~Dunwoody:

\begin{thm}
(M. Dunwoody \cite{Dunwoody(1979)}.) Let $R$ be an commutative
ring with a unit.

1. If $\Ga$ is a finitely generated torsion-free group with $cd_R(\Ga)=1$, then $\Ga$ is free.

2. If $\Ga$ is finitely-presented and $cd_R(\Ga)=1$ then $\Ga$ is a free product of finite and cyclic groups with
amalgamation over finite subgroups. In particular, $\Ga$ is virtually free.
\end{thm}

{\bf Duality groups.} A group $\Ga$ is said to be an
$m$-dimensional {\em duality group}, if $\Ga$ has type $FP$ and
$$
H^i(\Ga, R\Ga)\ne 0, \hbox{~for~} i=m \hbox{~and~}  H^i(\Ga, R\Ga)=0, \hbox{~for~} i\ne m.
$$
For instance, a finitely-presented group $\Ga$ is a 2-dimensional duality group (over $\Z$) if and only if
$cd(\Ga)=2$ and $\Ga$ does not split as a nontrivial free product.

\medskip
{\bf Poincar\'e duality groups.} Poincar\'e duality groups are
homological generalizations of the fundamental groups of closed
aspherical manifolds.

\begin{defn}
A group $\Ga$ is an (oriented) $m$-dimensional Poincar\'e
duality group over $R$ (a $PD(m)$-group for short) if $\Ga$ is of
type $FP$ and
$$
H^i(\Ga, R\Ga)\cong R, \hbox{~for~} i=m \hbox{~and~}  H^i(\Ga, R\Ga)=0, \hbox{~for~} i\ne m.
$$
\end{defn}

\noindent The basic examples are the fundamental groups of closed oriented aspherical $n$-manifolds.

This definition generalizes to (possibly non-oriented) $PD(n)$-groups, where we have to twist the module $V=R\Ga$
by an appropriate orientation character $\chi: R\Ga\to R$. The basic
examples are the fundamental groups of closed aspherical $n$-manifolds $M$. The character $\chi$ in this case corresponds to
the orientation character $\pi_1(M)\to R$.

We will need (in Section \ref{core}) the following relative version of the $PD(n)$ groups.

\begin{defn}
Let $\Ga$ be an $(n-1)$-dimensional group of type $FP$, and let
$$
\Del_1,\ldots, \Del_k\subset \Ga$$
be $PD(n-1)$ subgroups of $\Ga$. Set $\Del:=\{\Del_1,...,\Del_k\}$. Then,
the group pair $(\Ga, \Del)$ is an {\em $n$-dimensional Poincar\'e duality pair, or a $PD(n)$ pair},
if the double of $\Ga$ over the $\Del_i$'s is a $PD(n)$ group.
\end{defn}

We recall that the double of $\Ga$ over the $\Del_i$'s is the
fundamental group of the graph of groups ${\mathcal G}$, where
${\mathcal G}$ has two vertices labelled by $\Ga$, $k$ edges with
the $i$-th edge labelled by $\Del_i$, and edge monomorphisms are
the inclusions $\Del_i\ra \Ga$.

An alternate homological definition of $PD(n)$ pairs is the
following: A group pair $(\Ga, \Del)$ is a $PD(n)$ pair
if $\Ga$ and each $\Del_i$ has type $FP$, and
$$
H^*(\Ga, \Del; \Z\Ga)\simeq H^*_c(\R^n).
$$
If
$(\Ga, \Del)$ is a $PD(n)$ pair, where $\Ga$ and each $\Del_i$
admit a finite Eilenberg--MacLane space $X$ and $Y_i$
respectively, then the inclusions $\Del_i\ra \Ga$ induce a map
$$
\sqcup_i Y_i\ra X$$
whose mapping cylinder $C$ gives a {\em Poincar\'e pair}
$(C, \sqcup_i Y_i)$. The latter is a pair which satisfies Poincar\'e duality for manifolds with
boundary with local coefficients, where $\sqcup_i Y_i$ serves as
the boundary of $C$.  The most important example of a $PD(n)$ pair is the following. Let
$M$ be a compact manifold which is $K(\Ga,1)$. We suppose that the boundary of $M$ is the disjoint union
$$
\partial M= N_1\cup...\cup N_k,
$$
of $\pi_1$-injective components, each of which is a $K(\Del_i, 1)$, $i=1,...,k$. Then the pair
$$
(\Ga, \{\Del_1,...,\Del_k\})
$$
is a $PD(n)$ pair. See \cite{Bieri-Eckmann(1978)} for the details.

The following is one of the major problems in higher-dimensional topology:

\begin{conjecture}\label{wall}
(C.T.C. Wall, see a very detailed discussion in \cite{Kirby(1995)}.)
Suppose that $\Ga$ is an finitely-presented $n$-dimensional Poincar\'e duality group over $\Z$. Then there exits
a closed $n$-dimensional manifold $M$ which is $K(\Ga,1)$.
\end{conjecture}

This problem is open for all $n\ge 3$. The case $n=1$ is an easy corollary of the Stallings-Dunwoody theorem.
In the case $n=2$, the positive solution is due to Eckmann, Linnel and Muller, see
\cite{Eckmann-Linnel, Eckmann-Muller}. This result was extended to the case of fields $R$ by B.~Bowditch \cite{Bowditch(2004)}
and for general rings $R$ by M.~Kapovich and B.~Kleiner \cite{Kapovich-Kleiner(2004)} and B.~Kleiner \cite{Kleiner}.
Cannon's conjecture below is a special case (after Perelman's work) of Wall's problem for $n=3$:

\begin{conjecture}[J.\! Cannon]
\label{can}
Suppose that $\Ga$ is a Gromov-hyperbolic group whose ideal boundary is homeomorphic to $\S^2$.
Then $\Ga$ admits a cocompact properly discontinuous isometric  action on $\H^3$.
\end{conjecture}

\section{Topological rigidity and coarse compact core theorem}
\label{toprigidity}

First, few historical remarks. After the work of B.\! Maskit
\cite{Maskit(1970)} and A.~Marden \cite{Marden}, it became clear
that the major developments in the 3-dimensional topology
occurring at that time (in the 1960-s and the early 1970-s) were
of extreme importance to the theory of Kleinian groups. The key
topological results were:

1. Topological rigidity theorems of Stallings and Waldhausen.
Under appropriate assumptions they proved that  homotopy equivalence of
Haken manifolds implies homeomorphism,  see \cite{Hempel(1976)} for the detailed discussion.
In the context of Kleinian groups, it meant that the  (properly understood) algebraic structure
of a geometrically finite Kleinian group $\Ga\subset \mbd$ determines the topology of the associated
hyperbolic 3-manifold $\H^3/\Ga$.

2. Dehn Lemma, Loop Theorem and their consequences. The most important (for Kleinian groups)
of these consequences was the {\em Scott compact core theorem} \cite{Scott(1973a), Scott(1973b)}.
This theorem meant for (possibly geometrically infinite)
Kleinian groups, that the hyperbolic 3-manifold $M=\H^3/\Ga$ admits
a deformation retraction to an (essentially canonical) compact submanifold $M_c\subset M$ (the {\em compact core} of $M$).

\begin{rem}
Of course, after W.~Thurston entered the area of Kleinian groups, the theory experienced yet another radical change
and became the theory of hyperbolic 3-manifolds. However, this is another story.
\end{rem}

We now turn to the higher dimensions.

\subsection{Results of Farrell and Jones}

The following conjecture is a natural
generalization of the topological rigidity of 3-manifolds:

\begin{conjecture}
[A.~Borel]\label{borel}
 Let $M, N$ be closed aspherical $n$-manifolds and $f: M\to N$ is a homotopy-equivalence.
Then $f$ is homotopic to a homeomorphism. (There is also a relative version of this conjecture.)
\end{conjecture}

We refer the reader to \cite{Kirby(1995)} for a detailed discussion of Borel's Conjecture and its relation to Wall's Conjecture \ref{wall}.
Although, in full generality, Conjecture \ref{borel} is expected to be false, in the last 20 years there has been
a remarkable progress in proving this conjecture in the context to Kleinian groups. Most of these results appear in the works
of T.~Farrell and L.~Jones. We collect some of them below.

\begin{thm}
\label{fj2}
(T.\! Farrell and L.\! Jones \cite{Farrell-Jones(1989)}.) Suppose that $\Ga\subset \mob$ is a convex-cocompact
Kleinian group, $n\ge 4$, $N$ is a compact aspherical manifold (possibly with nonempty boundary $\D N$)
and $f: (\bar{M}^{n+1}(\Ga), M^{n}(\Ga)) \to (N, \D N)$ is a homotopy-equivalence which is a homeomorphism
on the boundary. Then $f$ is homotopic to a homeomorphism (rel. $ M^{n}(\Ga)$).
\end{thm}

\begin{thm}\label{fj3}
(T.\! Farrell and L.\! Jones, \cite[Theorem 0.1]{Farrell-Jones(1993c)})
Suppose that $X$ is a nonpositively curved closed Riemannian manifold,
$Y$ is a closed aspherical manifold of dimension $\ge 5$ and $f: X\to Y$ is a homotopy-equivalence.
Then $f$ is homotopic to a homeomorphism.
\end{thm}

\begin{thm}
\label{w=0}
(T.\! Farrell and L.\! Jones, \cite[Proposition 0.10]{Farrell-Jones(1998)})
For each Kleinian group $\Ga$ the Whitehead group $Wh(\Ga)$ is trivial.
\end{thm}

By combining Theorem \ref{w=0} with the $s$-cobordism theorem (see e.g. \cite{Kreck-Luck, Poenaru, Rourke-Sanderson(1982)}), one gets:

\begin{corollary}
\label{fj1} Suppose that $W^{n+1}$ is a topological (resp. PL,
smooth) $h$-cobordism so that $n\ge 4$ and $\pi_1(W^{n+1})$ is
isomorphic to a Kleinian group. Then $W$ is trivial in the
topological (resp. PL, smooth) category.
\end{corollary}

\subsection{Limit sets and homological algebra}

Let $\Ga\subset \mob$ be a convex-cocom\-pact subgroup with the limit set $\La$  and $R$ be a ring.
(One can also deal with geometrically finite groups by using relative cohomology, see \cite{Kapovich07}.)
The following theorem establishes a link between topology of the limit sets and the cohomology of $\Ga$:

\begin{thm}
(M. Bestina, G. Mess \cite{Bestvina-Mess(1991)}.)
$$
H^*(\Ga, R\Ga)\cong H^*_c(Hull(\La), R) \cong \tilde{H}^{*-1}(\La, R).
$$
Here we are using the Chech cohomology of the limit set.
\end{thm}

\noindent In particular, $\Ga$ is an $m+1$-dimensional duality group over $\Z$ (see Section \ref{homology}) iff $\t{H}^{*-1}(\La)$ vanishes
except in dimension $m$. In this case $\La$, homologically, looks like an infinite bouquet of $m$-spheres. Moreover
$$
cd(\Ga)=\dim(\La),
$$
which gives a geometric interpretation of the cohomological dimension of $\Ga$.

\subsection{Coarse compact core} \label{core}

In this section we state the best (presently) available higher-dimensional generalization of the
Scott compact core theorem.
The main drawback of this result is that it applies only to Kleinian groups
$\Ga\subset \mob$ which are $n$-dimensional duality groups.

We first need some definitions. We recall that every end $e$ of a manifold
$M$ admits a {\em basic system of neighborhoods}, which is a
decreasing sequence $(E_j)$ of nested connected subsets of
$M$ with compact frontier, so that
$$
\bigcap_{j\in \N} E_j=\emptyset.
$$
Given such $(E_j)$, we obtain the inverse system of the fundamental groups $(\pi_1(E_j))$. Consider the image $\Ga_j$
of $\pi_1(E_j)$ in $\pi_1(M)=\Ga$. Then
$$
\Ga_1\supset \Ga_2\supset ...
$$
The {\em fundamental group} of $e$ (or, rather, its image in $\Ga$) is defined as
$$
\pi_1(e)=\bigcap_{j\in \N} \Ga_j.
$$
An end $e$ is called {\em almost stable} if every sequence $(\Ga_j)$ as above is eventually constant,
in which case
$$
\pi_1(e)=\Ga_j
$$
for all sufficiently large $j$. (This notion is weaker than the notion of {\em semistabe} ends, which the reader might be familiar with.)
For instance, if $M$ is an open handlebody of finite genus,
then $M$ has unique end $e$, which is almost stable, whose fundamental group $\pi_1(e)$ is the free group $\pi_1(M)$.
On the other hand, if $S$ is a surface of infinite genus with a unique end $e$, then $e$ is not almost stable.
If $M$ is the complement to a Cantor set in $\S^2$, then no end of $M$ is almost stable.

The following is the {\em Coarse Compact Core Theorem} proved in \cite{Kapovich-Kleiner(2005)}
in the general context of {\em coarse Poincar\'e duality spaces}.

\begin{thm}\label{ad}
(M. Kapovich, B. Kleiner \cite{Kapovich-Kleiner(2005)}.)
Let $\Ga\subset \mob$ be a Kleinian subgroup which is an $n$-dimensional  duality group.
Then the manifold $M=\H^{n+1}/\Ga$ contains a compact submanifold $M_c$ (the {\em coarse compact core})
satisfying the following:

1. $\Ga$ contains a finite collection $\Del$ of $PD(n-1)$ subgroups $\Del_1,...,\Del_k$.

2. The pair $(\Ga, \Del)$ is an $n$-dimensional Poincar\'e duality
pair.

3. The group $\pi_1(M_c)$ maps onto $\pi_1(M)=\Ga$.

4. The manifold $M$ has exactly $k$ ends $e_1,...,e_k$, each of which is almost stable;
the components $E_1,..., E_k$ of $M\setminus M_c$ are basic neighborhoods of $e_1,...,e_k$.

5. For every $i=1,...,k$, $\pi_1(e_i)=\Del_i$ is the image of  $\pi_1(E_i)$ in $\Ga$.
\end{thm}

\begin{figure}[tbh]
\begin{center}
\input{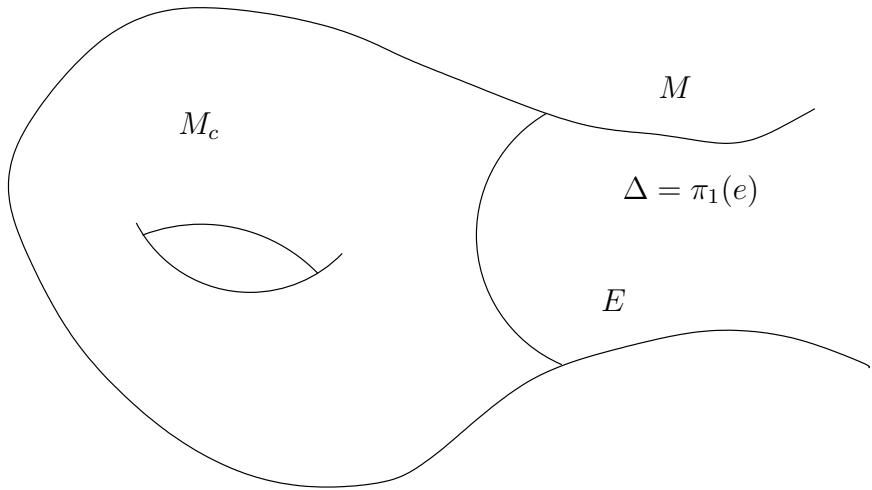}
\end{center}
\caption{\sl A coarse compact core.} \label{core.fig}
\end{figure}

See Figure \ref{core.fig}.  In the case when $n=3$, this theorem, of course, is a special case of the
{\em Scott compact core theorem} \cite{Scott(1973a), Scott(1973b)}.
More precisely, it covers the case when Scott compact core has incompressible boundary, for otherwise
$\Ga$ splits as a free product and is not a 2-dimensional duality group.
If $M$ is a tame manifold, e.g. $\Ga$ is geometrically finite, this theorem is also obvious.
At the moment, all known examples of Kleinian groups in $\mob$, $n\ge 3$, which are
$n$-dimensional duality groups, are geometrically finite.

\begin{problem}
Generalize Theorem \ref{ad}  to groups of type $FP$ which are not $n$-dimen\-sional duality groups. (Of course, the
conclusion of Part 1 of Theorem \ref{ad} would have to be suitably modified.)
\end{problem}

\section{Notions of equivalence for Kleinian groups}
\label{equivalences}

In this section we discuss various equivalence relations for Kleinian subgroups $\Ga_1, \Ga_2$ of $\mob$.
We start with the weakest one and end with the strongest.

\medskip
(0) {\bf Algebraic}: $\Ga_1$ is isomorphic to $\Ga_2$ as an abstract group.

\smallskip (1) {\bf Dynamical}: there exists a homeomorphism $f: \La(\Ga_1)\to \La(\Ga_2)$
such that $f \Ga_1 f^{-1} = \Ga_2$; i.e., the groups $\Ga_1$ and $\Ga_2$ have
the same topological dynamics on their limit sets. Thus, $\Ga_1$ is geometrically finite iff $\Ga_2$ is,
since geometric finiteness can be stated in terms of topological dynamics of a group on its limit set (Theorem \ref{BiMa}).

\smallskip (2) {\bf Topological conjugation}:
there exists a homeomorphism $f : \S ^n \to \S ^n$ such that $f \Ga_1 f^{-1} = \Ga_2$. (One can relax
this by assuming that $f$ is defined only on the domain of discontinuity of $\Ga_1$.)

\smallskip (3) {\bf Quasiconformal conjugation}: in (2) one can find a quasiconformal
homeomorphism. In the case $n=1$ one should replace {\em quasiconformal} with {\em quasisymmetric}.

\smallskip (4) {\bf Topological isotopy}: in (2) there exists a continuous
family of homeomorphisms $h_t: \S^n\to \S^n$ such that:
$h_0 = id$,~ $\forall t$,~ $h_t \Ga_1 h_t ^{-1} \subset \mob$ and
$h_1 \Ga_1 h_1 ^{-1} = \Ga_2$.

\smallskip (5) {\bf Quasiconformal isotopy}: in (4) all homeomorphisms are
quasiconformal (quasisymmetric).

\smallskip (6) {\bf M\"obius conjugation:} there is $f\in \mob$ such that $f\Ga_1 f^{-1}=\Ga_2$.

\smallskip
We refer the reader to \cite{Heinonen(2001), Iwaniec-Martin} for the definitions of quasisymmetric and quasiconformal  homeomorphisms.

\smallskip
Below is a collection of facts about the relation between different notions of equivalence of Kleinian groups.

Suppose that both groups $\Ga_j$ are {\em geometrically finite} and $\varphi : \Ga_1 \to \Ga_2$
is an isomorphism which preserves the {\em type} of elements, i.e., for $\ga\in \Ga_1$,
$\varphi(\ga)$ is hyperbolic if and only if $\ga$ is hyperbolic. It is clear that the above assumptions
are necessary for getting the equivalence (1). The following theorem shows that these assumptions are also sufficient.

\begin{theorem}\label{tukia}
(P.\ Tukia \cite{Tukia(1985a)}.) Under the above assumptions,
the isomorphism $\varphi$ can be realized by
the equivalence (1), i.e., there exists a (quasisymmetric) homeomorphism
$f$ of the limit sets, so that $f \ga f^{-1} = \varphi (\ga)$ for all $\ga \in \Ga_1$. Moreover, if $f: \Om(\Ga_1)\to \Om(\Ga_2)$
is a $\varphi$-equivariant quasiconformal (quasisymmetric) homeomorphism, then $f$ admits a $\varphi$-equivariant
quasiconformal (quasisymmetric) extension to the entire sphere.
\end{theorem}

\begin{question}\label{qctop}
{\bf (Quasiconformal vs. topological.)}
Suppose that two Kleinian groups $\Ga_{1}, \Ga_{2}\subset \mob$ are topologically
conjugate by a homeomorphism $f$ (defined either on $\Om(\Ga_1)$, or on $\La(\Ga_1)$, or on the entire $\S^n$),
which induces a type-preserving isomorphism  $\varphi: \Ga_1\to \Ga_2$.
Does it imply that $\varphi$ is induced by a quasiconformal (quasisymmetric) homeomorphism with the same domain as $f$?
\end{question}

Note that, for every $n$, the above question actually consists of 3 subquestions, depending on the domain of $f$.
Here is what is currently known about these questions:

\smallskip
1. If $n=1$ then all three questions have the affirmative answer and the proof is rather elementary.
It also follows for instance from Theorem \ref{tukia}.

\smallskip
2. If $n= 2$ then the answer to all three questions  is again positive, but the proof is highly nontrivial.
The easiest case is when the homeomorphism $f$ is defined on $\Om(\Ga_1)$. Then we get the induced
homeomorphism $\bar{f}$ of the quotient surfaces $S_1\to S_2$, where $S_i=\Om(\Ga_i)/\Ga_i$.
The existence of a diffeomorphism $S_1\to S_2$ homotopic to $\bar{f}$ follows from the uniqueness
of the smooth structure on surfaces. If $S_1$ is compact, then this diffeomorphism lifts to an equivariant
quasiconformal homeomorphism $\Om(\Ga_1)\to \Om(\Ga_2)$. Two noncompact surfaces can be
 diffeomorphic but not quasiconformally homeomorphic: For instance, the open disk is not quasiconformally
 equivalent to the complex plane. However, since $\varphi$ is type-preserving, Ahlfors Finiteness
 Theorem \cite{Ahlfors(1964)} in conjunction with a lemma of Bers and Maskit (see e.g. \cite[Corollary 4.85]{Kapovich00}),
 implies the existence of a quasiconformal homeomorphism $S_1\to S_2$.

If $f$ is defined on the limit set and $\Ga_1, \Ga_2$ are geometrically finite, then the positive answer is
a special case of Tukia's theorem \ref{tukia}. However, if the groups $\Ga_i$ are not geometrically finite, the proof
becomes very difficult and is a corollary of the solution of the Ending Lamination Conjecture
in the work of  J.~Brock, R.~Canary and Y.~Minsky in \cite{BCM,  Minsky03, Minsky(2003)},
and M.~Rees \cite{Rees04}.

Combination of the Ahlfors Finiteness Theorem with the Ending  Lamination Conjecture also gives the positive answer in
the case when $f$ is defined on $\S^2$.

3. If $f$ is defined on $\Om(\Ga_1)$, then the answer is positive provided that $n \neq  4$ and $M^n(\Ga_{1})$ is compact. This is a
consequence of the theorem of  D.\! Sullivan \cite{Sullivan77}, who proved uniqueness of the quasiconformal structure
on compact $n$-manifolds ($n\ne 4$): Apply Sullivan's theorem to the manifolds $M^n(\Ga_{i})$, $i=1, 2$, and lift
the quasiconformal homeomorphism to the domain of discontinuity.

\begin{rem}
An alternative proof of Sullivan's theorem and its generalization was given by  J.~Luukkainen in \cite{Luukkainen(2001)},
see also \cite{Tukia-Vaisala}.
\end{rem}

%Solution of the following appears to be needed for Sullivan's arguments in \cite{Sullivan77}:
%\begin{problem}
%For every $n$ construct a compact stably parallelizable hyperbolic $n$-manifold.
%\end{problem}
%Sullivan solves this by reference to his work with Deligne on complex vector bundles: Complex vector bundles with
%discrete structure group are virtually trivial.
%However the bundle in question here is not complex.

If $n= 4$, $f$ is defined on $\Om(\Ga_1)$, and $M^4(\Ga_{1})$ is compact, then the
situation is unclear but one probably should expect the negative answer since the uniqueness of quasiconformal structures
in dimension $4$ was disproved by S.~Donaldson and D.~Sullivan \cite{Donaldson-Sullivan}.

\begin{question}
Is there a pair of Kleinian groups $\Ga_1, \Ga_2\subset \mob$ so that the manifolds $M^n(\Ga_1), M^n(\Ga_2)$
are homeomorphic but not diffeomorphic?
\end{question}

Note that in view of the examples in \cite{Farrell-Jones89}, the positive answer to the above question would not be too surprising.

\medskip
If $f$ is defined on $\Om(\Ga_1)$ and we do not assume compactness of $M^n(\Ga_1)$, then the answer  to Question \ref{qctop}
is negative in a variety of ways.

(a) For instance, take {\em singly degenerate groups} $\Ga_1, \Ga_2\subset \mbd$, which are both isomorphic to the fundamental group
of a closed oriented surface $S$, contain no parabolic elements and have distinct ending laminations.
Then  $\Om(\Ga_i)\subset \S^2$ are open disks $D_i$ for both $i$. There exists an equivariant homeomorphism $h: D_1\to D_2$,
which induces an isomorphism $\varphi:\Ga_1\to \Ga_2$. However, since the ending laminations are different,
there is no equivariant homeomorphism $\La(\Ga_1)\to \La(\Ga_2)$.

Now extend both groups to the 3-sphere so that
$\Ga_i\subset \mbt$, $i=1 ,2$. Then the 3-dimensional domains  of discontinuity $B_i$ of both groups are diffeomorphic to
the open 3-ball, $i=1, 2$;  the quotient manifolds are
$$
M^3(\Ga_i)= B_i/\Ga_i\cong S\times \R, \quad i=1, 2.$$
Therefore there exists an equivariant diffeomorphism $f: B_1\to B_2$. We claim that this map
cannot be quasiconformal. Indeed, otherwise it would extend to an equivariant homeomorphism of the limit sets
(which are planar subsets of $\R^3$). This is a contradiction.

\medskip
(b) One can construct geometrically finite examples as well. The reason is that even though all (orientation-preserving)
parabolic elements of $\mbd$ are quasiconformally conjugate, the analogous assertion is
false for the parabolic elements of $\mbt$. Suppose that
$\tau$ is the translation in $\R^3$ by a nonzero vector $v$. Let $R_{\theta_i}$,  $i=1, 2$,
denote the rotations around $v$ by the angles $\theta_1, \theta_2\in [0, \pi]$. Then the skew motions
$$
\ga_i=R_{\theta_i} \circ \tau_i, \quad i=1, 2$$
are parabolic elements of $\mbt$. One can show that

\begin{prop}\label{nonqi}
The M\"obius transformations $\ga_1$ and $\ga_2$ are quasiconformally conjugate in $\S^3$ if and only if $\theta_1=\theta_2$.
\end{prop}

The proof is based on a calculation of the extremal length of a certain family of curves in $\R^3$ and  we will not present it here.

Note that the cyclic groups $\Ga_i=\<\ga_i\>$ are geometrically finite, the isomorphism $\varphi: \Ga_1\to \Ga_2$ sending
$\ga_1$ to $\ga_2$ is type-preserving. The quotient manifolds $M^3(\Ga_i)$ are both diffeomorphic to $\R^2\times \S^1$,
therefore there exists a $\varphi$-equivariant diffeomorphism $f: \Om(\Ga_1)\to \Om(\Ga_2)$ which, of course, extends to a
homeomorphism $\S^3\to \S^3$. However, according to Proposition \ref{nonqi},  this homeomorphism cannot be made quasiconformal.

These examples do not resolve the following:

\begin{question}
Suppose that $\Ga_1, \Ga_2\subset \mob$, $n\ge 3$, are Kleinian groups and
$f: \La(\Ga_1)\to \La(\Ga_2)$ is a homeomorphism which induces an isomorphism
$\Ga_1\to \Ga_2$. Does it follow that $f$ is quasisymmetric?
\end{question}

\medskip
If $n\le 2$, then (in the list of equivalences between Kleinian groups) we have the implication
$$
(3) \Rightarrow (5).
$$
Indeed, consider a quasiconformal homeomorphism $f$ conjugating Kleinian groups $\Ga_1$ and $\Ga_2$ and let $\mu$ denote the
Beltrami differential of $f$. Then for $t\in [0, 1]$ the solutions of the Beltrami equation
$$
\frac{\D f_t}{\D \bar{z}}= t\mu \frac{\D f}{\D {z}}
$$
also conjugate $\Ga_1$ to Kleinian subgroups of $\mbd$, see e.g.  \cite{Bers(1972)}. This gives the required quasiconformal isotopy.
Since (2) is equivalent to (3) for $n\le 2$, it follows that for $n\le 2$ we have
$$
(2) \iff (3) \iff (4) \iff (5)
$$

This argument however fails completely in higher dimensions, since the Beltrami equation in
$\R^n$ for $n\ge 3$ is overdetermined.
%The following is open for $n\ge 3$ even for convex-cocompact groups:

\begin{question}\label{deform}
In the list of equivalences between Kleinian groups:

(a)  Does (2) $\Rightarrow$ (4) ?

(b)  Does (3) $\Rightarrow$ (5) ?
\end{question}

One can show (using quasiconformal stability, see Section \ref{stability},
cf. \cite[Theorem 7.2]{Martin(1989a)})
%(cf. \cite[Theorem 7.2]{Martin(1989a)})
%proved in \cite{Epstein-Canary-Green}
that for convex-cocompact groups parts (a) and (b) of the above question are equivalent.
In Theorem \ref{nonisotopy} we give examples of convex-cocompact Kleinian groups in $\mob$, $n\ge 5$,
for which the answer to Question \ref{deform} is negative. The situation in dimensions 3 and 4
at the moment is unclear, but we expect in these dimensions the answer to be  negative as well.

\medskip
The implications (i)$\Rightarrow$(6) for i$\le 5$ are, of course, extremely rare. The most celebrated example is provided by the
Mostow rigidity theorem:

\begin{thm}
Suppose that $\Ga_1, \Ga_2\subset \mob$ are lattices and $n\ge 2$. Then (0)$\Rightarrow$(6) for these groups.
\end{thm}

See \cite{Mostow} for G.~Mostow's original proof or \cite{Kapovich06} for a more elementary argument along the same lines
which uses only the analytical properties of quasiconformal mappings.
A completely different argument due to M.~Gromov
can be found in \cite{Benedetti-Petronio}.  Yet another proof is an application of the barycentric maps
\cite{Besson-Courtois-Gallot(1999)}. Note that, presently, there are no proofs using equivariant harmonic maps.

\medskip
Mostow's ergodic arguments were greatly generalized by D.~Sullivan in \cite{Sullivan(1981b)}, see also \cite{Ahlfors(1980)}:

\begin{thm}
(D.~Sullivan \cite{Sullivan(1981b)}.) Suppose that $\Ga_1, \Ga_2\subset \mob$ are Kleinian groups whose limit set is the entire $\S^n$
and so that the action of $\Ga_1$ on $\S^n$ is {\em recurrent}. Then (3)$\Rightarrow$(6)  for these groups.
\end{thm}

The action of $\Ga\subset \mob$ on $\S^n$ is called {\em recurrent} if for every measurable subset $E\subset \S^n$
of positive Lebesgue measure, the measure of the intersection $\ga(E)\cap E$ is positive for some $\ga\in \Ga\setminus \{1\}$.

\section{Groups with zero-dimensional limit sets}
\label{0}

In what follows, we let $\dim$ denote the covering dimension of topological spaces, see for instance \cite{Hurewitz-Wallman}.
Suppose that $\Ga\subset \mob$ is a non-elementary Kleinian subgroup of $\mob$ and $\dim(\La(\Ga))=0$;
hence  $\La(\Ga)$ is totally disconnected (its only connected components are points). Recall that a {\em discontinuum}
is a nonempty perfect totally disconnected compact topological space, see e.g. \cite{Alexandroff}. Hence $\La(\Ga)$ is a discontinuum.
It follows (see e.g. \cite{Alexandroff}) that $\La(\Ga)$ is homeomorphic to the standard Cantor set $K\subset [0,1]$.
Below is a couple of examples of Kleinian groups whose limit sets are totally disconnected.

\begin{example}
(A {\bf Schottky group}, see e.g. \cite{KAG, Maskit(1987)}.) Let $n, k\ge 1$.
Suppose that we are given a collection of  disjoint  closed  topological
$n$-balls
$$
B^+_{1},\ldots, B^+_{k} , B^-_1 , \ldots, B^-_k\subset \S^n$$
and M\"obius transformations $\ga_{j}\in \mob$ so that $\ga_j( \int(B_{j}^+) ) =  \ext(B^-_j)$.
We assume that for each pair $B_j^+, B_j^-$ there exists a diffeomorphism of $\S^n$ which carries these balls to the
round balls.\footnote{By the smooth Schoenflies theorem, for $n\ne 4$ it suffices to assume that the balls
$B_j^+$ have smooth boundary.} Then
$$
\Phi:= \S^{n} -\bigcup^{k}_{j=1}( B^+_{j} \cup \int( B^-_j) )$$
is  a  fundamental domain for the group $\Ga$ generated by  $\ga_{1},\ldots, \ga_{k}$.
The group $\Ga$ is called a {\em Schottky group}. It is isomorphic to a free group of rank $k$,
and the limit set of $\Ga$ is a discontinuum provided that $k\ge 2$. Every nontrivial element of $\Ga$ is hyperbolic.
\end{example}

\begin{figure}[tbh]
\begin{center}
\input{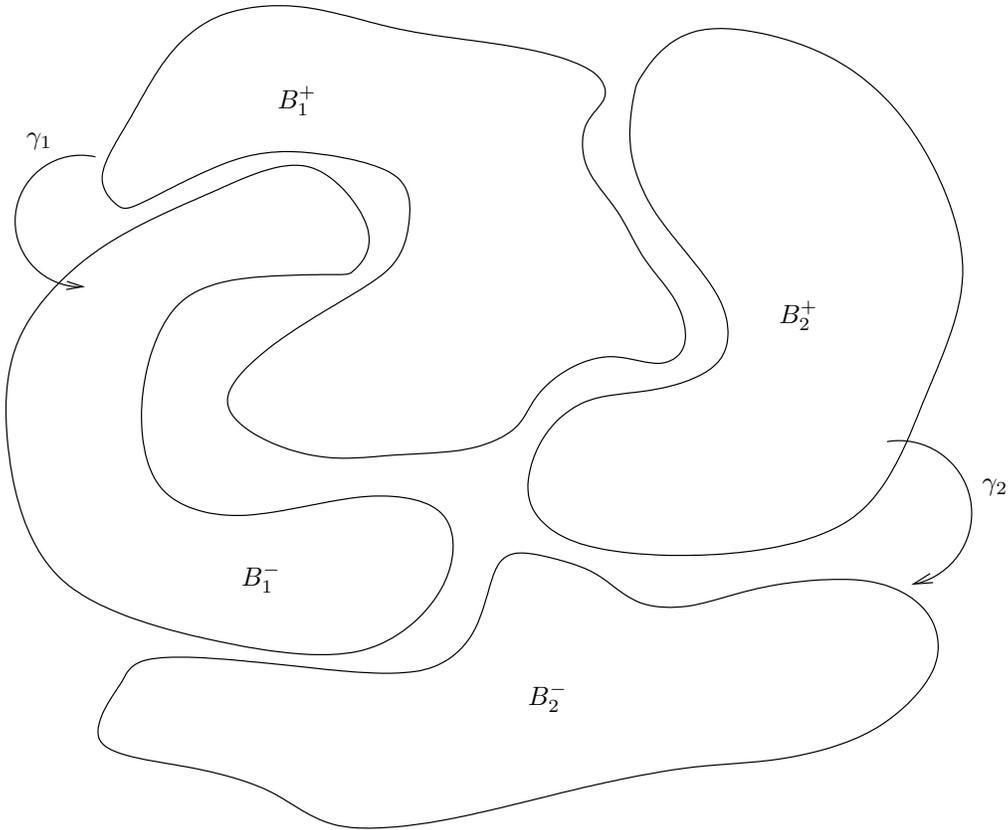}
\end{center}
\caption{\sl A Schottky group.} \label{sch.fig}
\end{figure}

Before giving the next example we need a definition.
Suppose that $\Ga\subset \mob$ is a nontrivial elementary subgroup.
Then, after conjugating $\Ga$ if necessary, we can assume that either:

\smallskip
1. $\Ga$ fixes $0, \infty\in \overline{\R^n}=\S^n$ and therefore is generated by $\ga(x)=Ax$, where $A$ is the product of a scalar $c>1$
by an orthogonal matrix.

2. Or $\Ga$ acts on $\R^n\subset \S^n$ by Euclidean isometries.

\smallskip
In the first case we take the fundamental domain $\Phi$ for the action of $\Ga$ on $\S^n$
to be an annulus bounded by two disjoint round spheres.
In the second case we take a Dirichlet fundamental domain $\Phi\subset \R^n$ for $\Ga$.

We refer to the fundamental domains $\Phi$ as {\em standard} fundamental domains. A fundamental domain for $\Ga$
is {\em topologically standard} if it is the image of a standard fundamental domain of $\Ga$ under a diffeomorphism of $\Om(\Ga)$
commuting with $\Ga$. For instance, the fundamental domain for a rank 1 Schottky group is topologically standard.
Therefore, the fundamental domain $\Phi$ for the Schottky group satisfies the property that it is the intersection of
topologically standard fundamental domains
$$
\Phi_j= \S^{n} -( B^+_{j} \cup \int( B^-_j) )
$$
for the groups $\Ga_j=\<\ga_j\>$.

Given a domain $\Phi\subset \S^n$, we let $\Phi^c\subset \S^n$ denote the closure of the complement
of $\Phi$. We are now ready for the second example which is a generalization of the first.

\begin{example}\label{extended}
({\bf Schottky-type groups, see e.g. \cite{KAG, Maskit(1987)}.})
Start with a collection of elementary Kleinian groups $\Ga_i\subset \mob$, $i=1,...,k$. Let $\Phi_i\subset \S^n$
be topologically standard fundamental domains for these groups. Assume that
$$
\Phi_i^c\cap \Phi_j^c=\emptyset
$$
for all $i\ne j$.
Let $\Ga\subset \mob$ be the group generated by $\Ga_1,...,\Ga_k$. Then:

1. As an abstract group, $\Ga$ is isomorphic to the free product $\Ga_1*...*\Ga_k$.

2. $\Phi:= \Phi_1\cap ... \cap \Phi_k$ is a fundamental domain for the group $\Ga$.

3.  The limit set of $\Ga$ is totally disconnected.

\smallskip
\noindent The groups $\Ga$ obtained in this fashion are called {\em Schottky-type groups}.
\end{example}

\begin{figure}[tbh]
\begin{center}
\input{type.pstex_t}
\end{center}
\caption{\sl A classical Schottky-type group isomorphic to $\Z * \Z^2$ with $\Z=\<\ga_1\>$, $\Z^2=\<\al_2, \be_2\>$.}
\label{type.fig}
\end{figure}

A Schottky-type group is called {\em classical} if it admits a fundamental domain $\Phi:= \Phi_1\cap ... \cap \Phi_k$,
so that each $\Phi_i$ is geometrically standard. It is not difficult to see
that Schottky-type groups are geometrically finite. For instance, consider the case of a Schottky group
$\Ga$ of rank $k$, for $n\ge 2$. We have  the map
$$
j: \Z= H_{n}(M^n(\Ga))\to H_n( \bar{M}^{n+1}(\Ga))
$$
induced by the inclusion of manifolds. Since the manifold $\bar{M}=\bar{M}^{n+1}(\Ga)$ is $K(\Ga,1)$, it follows that
$$
H^n(\bar{M})=H^n(\Ga)= H^n(B)=0,
$$
where $B$ is the bouquet of $k$ circles. Therefore $j=0$. Hence $\bar{M}$ is compact and hence $\Ga$ is convex-cocompact.
A similar argument works for Schottky-type groups, provided that one uses cohomology relative to the cusps.

The quotient manifolds of the Schottky-type groups $\Ga$ have a rather simple topology, as it follows from the
explicit description of their fundamental domains. Namely, let $M_i=M^n(\Ga_i)$, $i=1,...,k$. Then
we get the smooth connected sum decomposition
$$
M^n(\Ga)= M_1 \# ... \# M_k.$$
By combining this with Theorem \ref{tukia} we obtain

\begin{prop}
1. Suppose that $\Ga, \Ga'$ are Schottky groups of the same rank.
Then there exists a quasiconformal homeomorphism $f: \S^n\to \S^n$ which conjugates $\Ga$ to $\Ga'$, i.e., $f \Ga f^{-1}=\Ga'$.

2. Suppose that $\Ga, \Ga'$ are Schottky-type groups and $\varphi: \Ga\to \Ga'$ is a type-preserving isomorphism, so
that for every free factor $\Ga_i$ in $\Ga$, the restriction $\varphi: \Ga_i\to \Ga_i'\subset \Ga'$ is induced by
a quasiconformal homeomorphism of $\S^n$.
Then there exists a quasiconformal homeomorphism $f: \S^n\to \S^n$ which induces the isomorphism $\varphi$.
\end{prop}

\begin{question}
Let $n\ge 3$. Is there a quasiconformal isotopy between $\Ga$ and $\Ga'$ in the above theorem (either part 1 or part 2)?
\end{question}

In the case when $\Ga$ and $\Ga'$ are both classical, the positive answer follows rather easily.
In the non-classical case the above question is open even if $n=3$ and $\Ga$ is a Schottky group.

\medskip
Schottky subgroups of $\mbd$ can be characterized as follows:

\begin{thm}
( B.\! Maskit \cite{Maskit(1967)}.) A Kleinian subgroup $\Ga\subset \mbd$ is a Schottky
group if and only if $\Ga$ is free, has nonempty domain of discontinuity in $\S^2$ and
consists only of hyperbolic elements.
\end{thm}

This result was generalized by N.\! Gusevskii and N.\! Zindinova \cite{GZ}:

\begin{thm}
Let $\Ga\subset \mbd$ be a Kleinian subgroup, which has nonempty domain of discontinuity in $\S^2$
and is isomorphic to a Schottky-type group $\Ga'$ via a type-preserving isomorphism $\Ga\to \Ga'$.
Then $\Ga$ is a Schottky-type group.
\end{thm}

Both theorems are easy under the assumption that $\Ga$ is geometrically finite,
the key point here is that (in dimension 2) one can prove geometric finiteness under the above mild assumptions.

If $\Ga$ is a Kleinian subgroup of $\mbt$, then the above results are
 not longer true, moreover, $\Ga$ can be geometrically infinite.
For instance, take a free finitely generated purely hyperbolic discrete subgroup of
$PSL(2,\C)$, whose limit set is the  2-sphere (the existence of such groups
was first established by V.\! Chuckrow \cite{Chuckrow}). The M\"obius extension of this group to $\S^3$
has nonempty domain of discontinuity, but is not geometrically finite.

\medskip
{\bf Tameness of limit sets.} Below we address the following:

\begin{question}
Suppose that $\Ga\subset \mob$ is a Kleinian group, whose limit set is a discontinuum. What can be said about the
embedding $\La(\Ga)\subset \S^{n}$?
\end{question}

A discontinuum $D\subset \S^n$ is called {\em tame} if there exists a homeomorphism $f: \S^n\to \S^n$ which carries
$D$ to the Cantor set $K\subset [0, 1]$ and is called {\em wild} otherwise.
It is a classical (and easy) result that every discontinuum in $\S^2$ is tame, see e.g. \cite{Bing}.
The (historically) first example of a wild discontinuum was the {\em Antoine's necklace} $A\subset \S^3$:
$$
\pi_1(\S^3\setminus A)\ne \{1\},
$$
which explains why $A$ is wild, see \cite{Bing}.
D.~DeGryse and R.~Osborne \cite{DeGryse-Osborne} constructed for every $n\ge 3$
examples of wild discontinua $D_n\subset \S^n$, such that
$$
\pi_1(\S^n\setminus D_n)=\{1\}.
$$
See also \cite{GRZ} for infinitely many inequivalent 3-dimensional examples of this type.
%We refer the reader to \cite{Bing} for the comprehensive discussion of wild
%embeddings  of discontinua and other compacta in $S^3$.

The algebraic structure of Kleinian groups with totally disconnected limit sets is given by

\begin{thm}
(R.\! Kulkarni \cite{Kulkarni2}.)
Suppose that a Kleinian group $\Ga\subset \mob$  has a
totally disconnected limit set. Then $\Ga$ is isomorphic to a Schottky-type group.
\end{thm}

One can even describe (to some extent) fundamental domains of such groups:

\begin{thm}
\label{gsg}
(N.\! Gusevskii \cite{Gusevskii}.) Suppose that the limit set of $\Ga\subset \mob$  is
totally disconnected.  Then $\Ga$  admits a fundamental domain $\Phi$
of the same shape as in Example \ref{extended}, only
the fundamental domains $\Phi_i$ for $\Ga_i$'s are not required to be
topologically standard.
\end{thm}

The proof of Theorem \ref{gsg} is based on the following

\begin{thm}
(M.\! Brin \cite{Brin}.)
Let $\tilde M$ be a smooth oriented $n$-manifold of dimension $>2$, so that $H^{1}(\tilde{M})=0$. Let  $\Ga\acts \tilde M$
is a smooth properly discontinuous free action.

Then, for every smooth oriented compact hypersurface $\Sigma$ in $\tilde M$ and an
open neighborhood $U $ of $\Ga\cdot \Sigma$, there exists a
smooth compact connected oriented hypersurface $\Sigma^* \subset U$ such that for every $\ga \in \Ga$ either
$\ga\Sigma^* \cap \Sigma^* = \emptyset$ or $\ga\Sigma^*=\Sigma^*$.
\end{thm}

This theorem allows one to split (inductively) the Kleinian group $\Ga$
as a free product in a ``geometric fashion": Start with a compact hypersurface in $\Om(\Ga)$
which separates components of $\La(\Ga)$. Find $\Sigma^*$ as in Brin's theorem
which still separates. Then cut open the manifold
$M^n(\Ga)$ along the projection of $\Si^*$. This decomposition yields a free product decomposition
$\Ga=\Ga'*\Ga''$ so that $\Ga$ is a Klein combination of the groups $\Ga', \Ga''$. Continue inductively.
Finite generation of $\Ga$ implies that the
decomposition process will terminate and the terminal groups must be elementary.
Note that if all $\Sigma^*$ were spheres, then this decomposition would be of Schottky-type.

\begin{corollary}
Every Kleinian group with a totally disconnected limit set is
geometrically finite.
\end{corollary}
\proof Repeat the arguments which we used to establish geometric
finiteness of Schottky groups. \qed

\begin{problem}\label{stg}
Suppose that $\Ga\subset \mob$ is such that $\La(\Ga)$ is a tame discontinuum. Does is
follow that $\Ga$ is a Schottky-type group?
\end{problem}

If $n=2$ then the affirmative answer to this question follows for instance from Maskit's theorem.
If $n=3$ then the answer is again positive; moreover,

\begin{prop}
Suppose that $\Ga\subset \mbt$ is such that $\La(\Ga)$ is totally disconnected and
$\pi_1(\Om(\Ga))=1$. Then $\Ga$ is a Schottky-type group
\end{prop}
\proof Under the above assumptions, $\pi_2(\Om(\Ga))\ne \{0\}$;
hence, by the Sphere Theorem (see e.g. \cite{Hempel(1976)}), we can find a smooth hypersurface $\Si^*$ as in Brin's theorem,
so that $\Si^*$ is diffeomorphic to $\S^2$. Therefore, as we saw above, it follows that  $\Ga$ is a Schottky-type group. \qed

\smallskip
This argument however fails for $n\ge 4$, where Problem \ref{stg} is still open.
On the other hand, there are Kleinian subgroups of $\mbt$ with wild discontinua as limit sets.
The first such example was given by M.\! Bestvina and D.\! Cooper:

\begin{thm}
(M.\! Bestvina, D.\! Cooper \cite{Bestvina-Cooper}.)
There exists a  Kleinian group
$\Ga \subset \Mob(\S^3)$  which contains parabolic
elements, so that $\La(\Ga)$ is a wild discontinuum.
\end{thm}

The proof that $\pi_{1}(\Omega (\Ga)) \neq  1$ presented in \cite{Bestvina-Cooper}
was incomplete; however the gap was filled several years later by S.\! Matsumoto:

\begin{thm}
(S.\! Matsumoto \cite{Matsumoto89, Matsumoto92}, see also
 \cite{Gusevskii96}.) There are
Kleinian groups $\Ga$  in  $\Mob(\S^3)$ without parabolic elements whose
limit sets are wild discontinua.
\end{thm}

\section{Groups with one-dimensional limit sets}
\label{kk}

%We recall the construction of two standard examples of 1-dimensional compact spaces.
The simplest examples of 1-dimensional limit sets of Kleinian groups are topological circles.
For instance, the limit set of a lattice $\Ga\subset \Isom(\H^2)$ is the round circle.
Of course, even if the limit set of $\Ga\subset \mob$ is a topological circle, its embedding in $\S^n$ can be complicated.
We will discuss this issue later on. For now, we are only interested in the topology of the limit set itself.

\medskip
Given convex-cocompact Kleinian groups $\Ga_1, \Ga_2\subset \mob$ with 1-dimensional limit sets, one can
use Klein-Maskit Combination theorems in order to get convex-cocompact
Kleinian groups $\Ga\subset \mob$ isomorphic to
$$
\Ga_1*_{\Del}\Ga_2,$$
where $\Del$ is either trivial or infinite cyclic.
The limit sets of the resulting groups are again 1-dimensional. For instance, if $\La(\Ga_i)$ is a topological circle
for $i=1, 2$ then the limit set of  $\Ga=\Ga_1*\Ga_2$ will be disconnected: The connected components of $\La(\Ga)$
are topological circles and points. Similarly, if $\Del=\Z$,
then the limit set of $\Ga=\Ga_1*_{\Del} \Ga_2$ will have {\em cut pairs}:
The complement to the 2-point set $\La(\Del)$ in $\La(\Ga)$ is disconnected.

\begin{figure}[tbh]
\begin{center}
\input{comb.pstex_t}
\end{center}
\caption{\em Combination of two 1-\qf groups: $\Ga= {\Ga_1} *_{\Z}{\Ga_2}$.} \label{comb.fig}
\end{figure}

\medskip
These constructions are, of course, not very illuminating.
Therefore we are interested in examples of 1-dimensional limit sets which are connected and which do not contain
cut-pairs. It turns out that there are only two such examples:

\medskip
{\bf 1. The Sierpinski carpet ${\mathcal S}$.} Start with the unit square $S=I\times I$.
Subdivide this square into 9 squares of the size $\frac{1}{3}\times \frac{1}{3}$ and
then remove from $S$ the open middle square
$(\frac{1}{3}, \frac{2}{3})\times (\frac{1}{3}, \frac{2}{3})$. Repeat this
for each of the remaining $\frac{1}{3}\times \frac{1}{3}$ sub-squares in $S$ and continue inductively.
After removing a countable collection of open squares we are left with a compact
subset  ${\mathcal S}\subset \R^2$, called the {\em Sierpinski carpet}.

\medskip
{\bf 2. The Menger curve ${\mathcal M}$.} Start with the unit cube $Q=I\times I\times I$.
Each face $F_i$ of $Q$ contains a copy of the Sierpinski carpet ${\mathcal S}_i$. Let $p_i: Q\to F_i$
denote the orthogonal projection. Then
$$
{\mathcal M}:= \bigcap_{i} p^{-1}_i({\mathcal S}_i)
$$
is called the {\em Menger curve}.

\begin{example}
There exists a convex-cocompact subgroup $G\subset \Mob(\S^2)$
whose limit set is homeomorphic to the Sierpinski carpet ${\mathcal S}$.
\end{example}

\begin{figure}[h]
\leavevmode
\centerline{\epsfxsize=4in \epsfbox{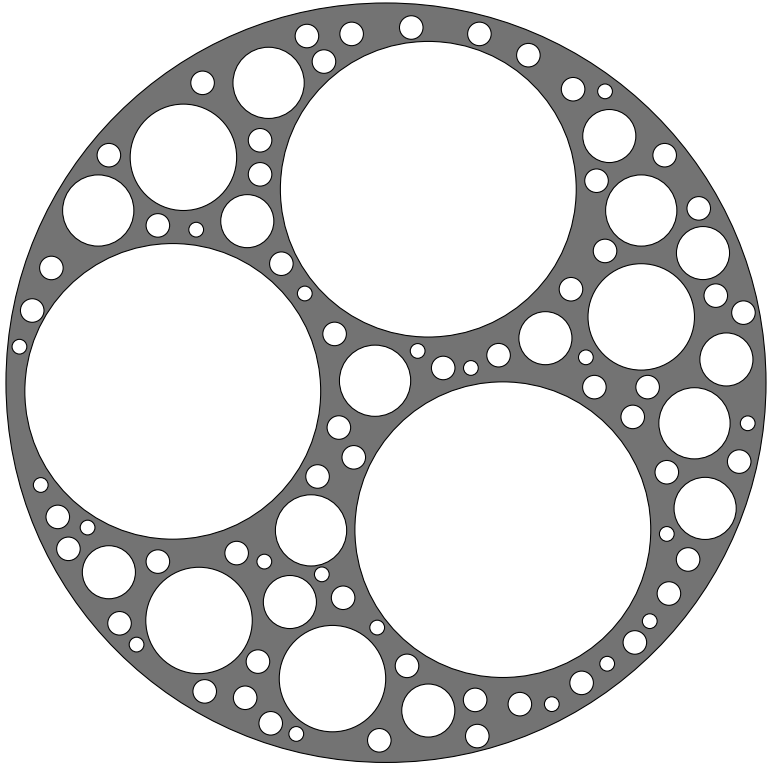}}
\caption{\sl A limit set homeomorphic to the  Sierpinski carpet.}
\label{carpet}
\end{figure}

To construct such an example start with a compact
hyperbolic manifold $M^3$ with nonempty totally-geodesic boundary. Thus we get an embedding
of $\Ga=\pi_1(M^3)$ into $\mbd$ as a convex-cocompact Kleinian subgroup.
The limit set of $\Ga$ is homeomorphic to the Sierpinski carpet.
To see this note that the convex hull
$Hull(\La(\Ga))$ in $\H^3$ is obtained by removing from $\H^3$ a countable
collection of disjoint open half-spaces $H_j\subset \H^3$.
The ideal boundary of each $H_j$ is the open round disk $D_j\subset \S^2$.
Thus
$$
\La(\Ga)=\S^2\setminus \bigcup_{j} \int(D_j).
$$
Clearly, $D_j\cap D_i=\emptyset$, unless $i=j$; since $\La(\Ga)$ has empty interior. See Figure \ref{carpet}.
According to Claytor's theorem \cite{Claytor}, it follows that
$\La(\Ga)$ is homeomorphic to ${\mathcal S}$. Moreover, it is easy to see
that this homeomorphism extends to the 2-sphere, since it sends the boundary circles of $\La(\Ga)$ to
the boundary squares of ${\mathcal S}$.

\medskip
The construction of Kleinian groups whose limit sets are homeomorphic to ${\mathcal M}$ is more complicated:

\begin{example}
\label{bou}
(M.\! Bourdon \cite{Bourdon97}; see also \cite{Kapovich05}.) There exists a
convex-cocompact subgroup $\Ga\subset \Mob(\S^3)$
whose limit set is homeomorphic to the Menger curve ${\mathcal M}$.
\end{example}

The following theorem is proved in \cite{Kapovich-Kleiner00} in the more general context of Gromov-hyperbolic groups:

\begin{thm}
(M.\! Kapovich, B.\! Kleiner \cite{Kapovich-Kleiner00}.) Suppose that $\Ga\subset \mob$ is a
(torsion-free) nonelementary convex-cocompact subgroup such that:

(a) $\Ga$ does not split as a free product,

(b) $\Ga$ does not split as an amalgam over $\Z$,

(c) $\dim(\La(\Ga))=1$.

Then $\La(\Ga)$ is either homeomorphic to the Sierpinski carpet or to the Menger curve.
\end{thm}

\begin{conjecture}
(M.\! Kapovich, B.\! Kleiner \cite{Kapovich-Kleiner00}) If $\Ga\subset \mob$ is a (torsion-free)
convex-cocompact Kleinian group whose limit set is homeomorphic to the Sierpinski carpet,
then $\Ga$ is isomorphic to a convex-cocompact subgroup in $\Mob(\S^2)$.
\end{conjecture}

It was proved in \cite{Kapovich-Kleiner00} that this conjecture would follow either from the positive solution
of the 3-dimensional Wall's problem (Problem \ref{wall}) or from Cannon's conjecture
(Conjecture \ref{can}).

\medskip
Topology of the limit sets of geometrically infinite Kleinian groups can be more complicated.
A {\em dendroid} is a compact locally connected simply-connected 1-dimensional topological space.

\begin{thm}
(J. Cannon and W. Thurston \cite{Cannon-Thurston}, see also \cite{Abikoff(1987)} and \cite{Cooper-Long-Reid(1994)}.)
There exist singly-degenerate Kleinian groups whose limit sets are {\em dendroids}.
\end{thm}

Conjecturally, limit sets of all singly-degenerate Kleinian groups are dendroids and
the following problem is open even for $n=2$:

\begin{problem}
Suppose that $\Ga\subset \mob$ is a Kleinian group whose limit set is connected and 1-dimensional. Is it true
that $\La(\Ga)$ is locally connected?
\end{problem}

See \cite{McMullen(2001)}, \cite{Minsky(1994)}   and \cite{Mitra05a, Mitra06a} for partial results in dimension $2$.

\section{Groups whose limit sets are topological spheres}
\label{spheres}

\begin{definition}
A Kleinian group $\Ga\subset \mob$ is called {\em $i$-fuchsian}\footnote{Our definition is
somewhat different from the classical: fuchsian subgroups of $PSL(2,\C)$
are usually required to preserves a round disk in $\S^2$.}
if $\La(\Ga)$ is a round $i$-dimensional  sphere in $\S^n$.
\end{definition}

To construct examples of $i$-fuchsian groups start with a lattice $\Ga\subset \Mob(\S^i)$.
The limit set of $\Ga$ is the round sphere $\S^i$. Define the {\em canonical embedding}
$$
\iota: \Mob(\S^i) \hook \mob
$$
induced by the embedding of the Lorentz groups
$$
O(i+1, 1)\hook O(n+1, 1)
$$
$$
A\mapsto \left[\begin{array}{cc}
A& 0\\
0& I
\end{array}\right],
$$
where $I$ is the identity matrix. Therefore we get the {\em canonical} embedding
$$
\iota: \Ga\hook \mob.
$$
One can modify this construction as follows.
Note that the stabilizer of $\S^i$ in $\mob$ contains $\Mob(\S^i)\times SO(n-i)$.
Choose a homomorphism $\phi: \Ga\to SO(n-i)$. For instance, being residually finite, the
group $\Ga$ will have many epimorphisms to finite groups, which then can be embedded in  $SO(n-i)$ if $n-i$ is sufficiently large.
Alternatively, in many cases the group $\Ga$ will have infinite abelianization $\Ga^{ab}$. The abelian group  $\Ga^{ab}$ admits
many embeddings into $SO(n-i)$ provided that $n-i\ge 2$. Then the image of
$$
\rho=\iota\times \phi: \Ga\to \Mob(\S^i)\times SO(n-i)\subset \mob
$$
is also an $i$-fuchsian group, since $\rho(\Ga)$ preserves $\S^i$ and the
action of $\rho(\Ga)$ on $\S^i$ is the same as the action of $\Ga$.

\begin{definition}
A Kleinian group $\Ga\subset \mob$ is called {\em $i$-quasifuchsian} if its limit set
is a topological $i$-dimensional sphere.
\end{definition}

We will refer to the number $n-i$ as the {\em codimension} of a (quasi)fuchsian group $\Ga$.

\begin{example}
Suppose that $\Ga$ is an $i$-fuchsian subgroup of $\mob$ and $\Ga'\subset \mob$ is another
group which is topologically conjugate to $\Ga$ (with a homeomorphism $f$ defined on the entire
$n$-sphere). Then $\Ga'$ is $i$-quasifuchsian. However, as we will see, there are $i$-quasifuchsian groups
(for $n\ge 3$) which cannot be obtained in this fashion.
\end{example}

The (quasiconformal) homeomorphisms $f$ as in the previous example exist in abundance if
$i=1, n=2$, due to the solvability of the Beltrami equation.
If $i\ge 2$, the situation is very different and it is not so easy
to construct nontrivial examples of $i$-quasifuchsian groups which are not fuchsian.
Some of these examples will be discussed below.

The following result was proved by M.~Bestvina and G.~Mess \cite{Bestvina-Mess(1991)} in the context of Gromov-hyperbolic groups:

\begin{thm}
Each convex-cocompact $i$-\qf group is a {\em Poincar\'e duality group} (over $\Z$) of dimension $i+1$. Conversely, if
$\Ga\subset \mob$ is a convex-cocompact  Poincar\'e duality group, then $\La(\Ga)$ is a homology manifold which is
a homology sphere.
\end{thm}

\begin{question}
Is it true that each convex-cocompact \qf group is isomorphic to the fundamental group of a
closed aspherical manifold?
\end{question}

This is, of course, a special case of Wall's problem (Problem \ref{wall}).

\begin{question}
Is there a convex-cocompact group $\Ga\subset \mob$ whose limit set is a homology  manifold which is homology sphere, so
that $\La(\Ga)$ is not homeomorphic to a sphere?
\end{question}

\subsection{Quasifuchsian groups of codimension 1}\label{codimension1}

The situation in the case of $n=2$ is completely understood due to the following:

\begin{thm}
(B.\! Maskit \cite{Maskit(1970)}, see also \cite{Marden}.)
Let $\Ga \subset \Mob(\S^2)$  be a Kleinian group whose domain of discontinuity
$\Om(\Ga)$ consists of precisely two components.
Then:

1. $\Ga$ is 1-quasifuchsian and geometrically finite.

2. $\Ga$ is quasiconformally conjugate to a 1-fuchsian group.

3. $\bar{M}^3(\Ga)=({\H}^{3}\cup\Omega(\Ga))/\Ga$ is homeomorphic to an interval bundle over a surface $S$, which
is 2-fold covered by $\Omega(\Ga)/\Ga$.
\end{thm}

Our goal is to compare the higher-dimensional situation with this theorem.
Suppose that $\Ga\subset \mob$ is a codimension $1$ \qf group. Then
$\Omega(\Ga)$ consists of two components, $\Om_1, \Om_2$. After
replacing $\Ga$ by an appropriate index 2 subgroup, we can assume that each $\Om_i$ is $\Ga$-invariant; hence
$M^n(\Ga) = M_{1} \cup  M_{2}$, where $M_i:= \Om_i/\Ga$. Then, by the Alexander duality,
$H_{*}(\Omega _{i}) \cong H_{*}(point)$, $i=1, 2$.
Therefore, if $\Om_i$ is simply-connected, then $\Om_i$ is contractible. Below we discuss what is currently known
about such \qf groups for $n\ne 4$.

\begin{thm}
\label{fjt}
Suppose that both $M_i$ are compact and both $\Om_i$ are simply-connected. Then
$\bar{M}^{n+1}(\Ga)$ is diffeomorphic to $M_{1} \times  [0, 1]$ provided that $n\ge 5$.
\end{thm}
\proof Note that, for homological reasons, $W=\bar{M}^{n+1}(\Ga)$ is compact, hence $\Ga$ is
convex-cocompact, see Theorem \ref{T2}. Since both $\Om_1, \Om_2, \H^{n+1}$ are contractible, the inclusions
$$
M_i\hook W, \quad i=1, 2
$$
are homotopy-equivalences. Therefore $W$ defines a smooth h-cobordism between the aspherical manifolds
$M_{1}$ and $M_{2}$. According to Corollary \ref{fj1}, this h-cobordism is smoothly trivial.
\qed

Suppose that $n=3$ and both $\Om_1, \Om_2$ are contractible. Then
$$
\pi_1(M_1)\cong \pi_1(M_2)\cong \Ga,
$$
the manifolds $M_1$ and $M_2$ are both irreducible and have infinite fundamental groups. If $\Ga$ were to contain
a subgroup $\Pi$ isomorphic to $\Z^2$, the subgroup $\Pi$ would be parabolic. This would contradict compactness of $W$.
Therefore, according to Perelman's solution of Thurston's hyperbolization conjecture, there exists a closed
hyperbolic 3-manifold $M_0$ which is homeomorphic to $M_1$ and $M_2$. Since $M_0$ is hyperbolic, its fundamental group $\Ga_0$
acts as a $2$-fuchsian group on $\S^3$. Therefore, according to our discussion in Section \ref{equivalences},
the group $\Ga$ is quasiconformally conjugate to $\Ga_0$. It is not difficult to see that passage to the index
2 subgroup which we used above does no harm and we obtain:

\begin{prop}
\label{precor}
Suppose that $\Ga$ is a codimension 1 \qf subgroup of $\mbt$, so that both components of $\Om(\Ga)$ are simply-connected and
$M^3(\Ga)$ is compact. Then $\Ga$ is quasiconformally conjugate to a $2$-fuchsian group $\Ga_0\subset \mbt$.
\end{prop}

On the other hand, we do not know if the 4-dimensional manifold $\bar{M}(\Ga)$ is
homeomorphic (or diffeomorphic) to an interval bundle over a 3-manifold.

Proposition \ref{precor} fails for $n\ge 4$:

\begin{thm}
For every $n\ge 4$ there are codimension 1 convex-cocompact \qf subgroups $\Ga\subset \mob$, so that both components of $\Om(\Ga)$
are simply-connected, but $\Ga$ is not isomorphic to a fuchsian group.
\end{thm}
\proof We give only a sketch of the proof, the details will appear
elsewhere. Fix $n\ge 4$. M.~Gromov and W.~Thurston in
\cite{Gromov-Thurston(1987)} construct examples of negatively
curved compact conformally-flat $n$-manifolds $M^n$, so that $M^n$
is not homotopy-equivalent to any closed hyperbolic $n$-manifold
$N^n$.   (See also \cite{Kapovich07b} for a review of the
Gromov--Thurston examples and for a construction of a convex
projective structure on $M^n$.)

By choosing parameters in the construction of \cite{Gromov-Thurston(1987)} more carefully, one
can construct an example of a {\em uniformizable} flat conformal manifold $M^n$ with the same properties.
Moreover, $M^n = \Omega_1/\Ga$, $\Ga\subset \mob$ is convex-cocompact, and
$\Omega(\Ga) = \Omega_1 \cup \Omega_2$ is the union of two simply-connected components.
Then $\Lambda(\Ga)$ is homeomorphic to $\S^{n-1}$, since the limit set of
$\Ga$ is homeomorphic to the ideal boundary of the universal cover of the negatively curved manifold $M^n$.
If $\Ga$ were isomorphic to an $n-1$-fuchsian group $\Ga'$, then $\Ga'$ would be isomorphic to
the fundamental group of a closed hyperbolic $n$-manifold $N^n$, which is a contradiction. \qed

\medskip
The above examples have another interesting property.
Let $\Om^{n+1}$ denote the domain of discontinuity of the group $\Ga$
(regarded as a subgroup of $\MOB$). Note that
$\Om^{n+1}$ is connected and $\pi_1(\Om^{n+1})\cong \Z$.
Since both $\Om_1, \Om_2$ are contractible, it follows that
$$
\pi_i(\Om^{n+1})=0, \quad i\ge 2.
$$
Set $M^{n+1}= \Om^{n+1}/\Ga$. We have the short exact sequence
$$
1\to \Z=\pi_1(\Om^{n+1}) \to \pi_1(M^{n+1}) \to \Ga\to 1.
$$
The embedding $M_1\to M^{n+1}$ determines a splitting of this sequence. Hence
the manifolds $M^{n+1}$ and $\S^1 \times M^n$ are homotopy-equivalent.
Given the existence of a metric of negative curvature on $M^n$, we obtain a
metric of nonpositive curvature on $\S^1 \times M^n$. Therefore, by Theorem \ref{fj3},
the manifolds $M^{n+1}$ and $\S^1 \times M^n$ are homeomorphic.

Let  $k\Z\subset \Z\subset \pi_1(M^{n+1})$ be the index $k$ subgroup in the center of $\pi_1(M^{n+1})$.
Then we obtain the $k$-fold covering $X_k\to X_1=M^{n+1}$, where
$$
\pi_1(X_k)= k\Z\times \Ga\subset \pi_1(X_1).
$$
Since the manifolds $X_k$ have isomorphic fundamental groups and $\pi_i(X_k)=0$ for all $i\ge 2, k\in \N$,
these manifolds are all homeomorphic to the smooth manifold
$X_1$ by Theorem \ref{fj3}. By \cite[Essay IV]{Kirby-Siebenmann(1977)}, there only finitely many smooth structures
on the manifold $X_1$. Therefore we obtain an infinite family of diffeomorphic manifolds
$$
M_j^{n+1}:=X_{k_j}, j\in \N
$$
and smooth covering maps $p_j: M_j^{n+1}\to M^{n+1}$.

The $(n+1)$-manifold $M^{n+1}=\Om^{n+1}/\Ga$ has the flat conformal structure $K_1$ uniformized
by the group $\Ga$. Let $K_j$ denote the flat conformal structure  on   $M^{n+1}$, which is the
lift  of $K_1$ via $p_j$.\footnote{The structures $K_j$
are obtained via {\em grafting} of $(M^{n+1},K_1)$ along the hypersurface $M^n$.}
We thus obtain an infinite family of diffeomorphic
flat conformal manifolds
$$
(M_j^{n+1}, K_j), j= 1, 2, ...
$$

\begin{question}
Suppose that $M$ is a closed hyperbolic $n$-manifold. Is there a finite cover $f:M'\to M$ such
that the pull-back map $f^*: H^3(M,\Z/2)\to H^3(M',\Z/2)$ is trivial?
(Recall \cite{Kirby-Siebenmann(1977)} that the group
$H^3(M,\Z/2)$ classifies PL structures on $M$ if $n\ge 5$.)
\end{question}

We regard the structures $K_j$ as elements of ${\mathfrak M}(X)$, the moduli space of the flat conformal structures
on a fixed smooth manifold $X$. The proof of the following claim is similar to \cite{Kapovich(1990)},
where it was proved in the context of 3-manifolds:

\begin{claim}
\label{disconnect}
For different $i, j$ the structures $K_i, K_j$ lie in
different connected components of the moduli space ${\mathfrak M}(X)$.
Thus ${\mathfrak M}(X)$ consists of infinitely many connected components.
\end{claim}

We note that K.\! Scannell in \cite{Scannell02}  constructed examples of hyperbolic
3-manifolds $X$ for which ${\mathfrak M}(X)$ consists of infinitely many
connected components.

To get the same phenomenon in dimension 4 consider one of the hyperbolic manifolds $M^3$
obtained by Dehn surgery on a 2-bridge knot, so that the natural embedding
$$
\Ga=\pi_1(M)\hook \mbt
$$
is (locally) rigid (see Example \ref{rigid}). Then the natural embedding of $\Ga$ to $\mbc$ is also  rigid and hence
the manifold $M^4=M^4(\Ga)\cong M\times \S^1$ has rigid  flat conformal structure. Taking $k$-fold covers of this
manifold we obtain infinitely many rigid flat conformal structures on $M^4$. By combining these results we obtain

\begin{thm}
For every $n\ge 3$ there exists a smooth compact $n$-manifold $X^n$ such that ${\mathfrak M}(X^n)$ consists of infinitely many
connected components.
\end{thm}

We now return to our discussion of Kleinian groups, restricting to $n=3$. Suppose that
$\Ga$ is a convex-cocompact $2$-\qf group, such that both components of $\Om(\Ga)$ are simply-connected.
Then, by proposition \ref{precor}, the limit set of $\Ga$ is {\em tame}, i.e., there
is a  homeomorphism of ${\S}^{3}$ which maps $\Lambda(\Ga)$  to the round sphere.

\begin{thm}
(B.\! Apanasov and A.\! Tetenov \cite{Apanasov-Tetenov}.)
There exists  a convex-cocompact 2-\qf group $\Ga\subset \Mob(\S^3)$ whose limit set is a
{\em wild} 2-sphere, i.e., there is no homeomorphism of ${\S}^{3}$
which maps $\Lambda(\Ga)$  to the round sphere. Moreover, one component of $\Om(\Ga)$
is simply-connected.
\end{thm}

\subsection{$1$-quasifuchsian subgroups of $\Mob(\S^3)$}

\medskip
Given a Kleinian subgroup $\Ga\subset \Mob(\S^3)$ whose limit set is a topological
circle $C$, we would like to analyze the embedding $C\hook \S^3$. It is clear
that $C$ could be an unknot in $\S^3$ (i.e., there exists a homeomorphism of
$\S^{3}$ which maps $C$ to a round circle), take for instance any 1-fuchsian subgroup of $\mbt$.

A topological circle $C$ in $\S^3$ is called {\em tame} if it is isotopic to
a polygonal knot in $\S^3$; if $C$ is not tame, it is called {\em wild}.

\begin{prop}
1. If $\Ga$ is a 1-\qf subgroup of $\Mob(\S^3)$, then either $\La(\Ga)$ is an unknot or
it is a wild knot $K$ such that $\pi_1(\S^3\setminus K)$ is infinitely generated.

 2. Each 1-\qf group is geometrically finite.
\end{prop}
\proof Since $\Ga$ is a 1-\qf subgroup of $\Mob(\S^3)$, this group is
nonelementary. The fundamental group of $M=M^3(\Ga)$ is finitely generated
(since $\Ga$ is) and we have the exact sequence:
$$
1\to \pi_1(\Om(\Ga))\to \pi_1(M)\to \Ga\to 1.
$$
Suppose that $\pi_1(\Om(\Ga))$ is finitely generated. Then, according to Jaco-Hempel's Theorem
\cite{Hempel(1976)},  $\pi_1(\Om(\Ga))\cong \Z$. This immediately excludes
tame nontrivial knots (the result proved by R.~Kulkarni \cite{Kulkarni1}). It remains to exclude wild knots
with
$$
\Del:=\pi_1(\Om(\Ga))\cong \Z.$$
Without loss of generality (after passing to an index 2 subgroup in $\Ga$), we can assume that $\Del$ is
contained in the center of $\pi_1(M)$. Note that $M$ is a Seifert manifold since
its fundamental group has infinite center.
Hence $M$ admits an $\S^1$-action. Lift this action to $\Om(\Ga)$ and then extend it to the
entire 3-sphere (so that the fixed point set is the limit set).
Raymond's classification \cite{Raymond} of topological $\S^1$--actions on $\S^3$ implies that
this $\S^1$--action is topologically conjugate to the orthogonal action, hence
$\La(\Ga)$ is an unknot.  This proves (1).

To prove (2) note that the group $\Ga$ acts as a {\em convergence group} on $\S^1=\La(\Ga)$
(see \cite{Tukia(1988)}).
Hence, according to \cite{Tukia(1988)}, there exists a homeomorphism $f: \La(\Ga)\to \S^1$
such that
$$
f \Ga f^{-1}=\Ga'\subset PSL(2,\R).
$$
Since finitely generated discrete subgroups of $PSL(2,\R)$ are geometrically finite, it follows that
$\Ga'$ is geometrically finite. As geometric finiteness is an invariant of the topological
dynamics on the limit set (see Theorem \ref{BiMa}), the group $\Ga$ is geometrically finite as well. \qed

\medskip
On the other hand, even if $\La(\Ga)$ is an unknot, the 3-manifold $\Om(\Ga)/\Ga$ is not necessarily a product:

\begin{thm}
\label{gltk}
(M.\! Gromov, B.\! Lawson, W.\! Thurston \cite{Gromov-Lawson-Thurston(1988)}, N.\! Kuiper \cite{Kuiper},
and M.~Kapovich \cite{Kapovich(1989b), Kapovich(1993a)}.)
There are 1-\qf groups $\Ga\subset \Mob(\S^3)$  such that $\Lambda(\Ga)$ are
unknotted but $\Ga$ are not topologically conjugate to 1-fuchsian groups.
\end{thm}

In the examples constructed in this theorem, the manifolds $M^3(\Ga)$ are
 nontrivial oriented circle bundles over orientable surfaces. On the other hand, for every
1-fuchsian group $\Ga_0\subset \mbt$, the manifold $M^3(\Ga_0)$ is a 3-dimensional
Seifert manifold with the zero Euler number, since it admits a natural $\H^2\times \R$--structure.
Hence, in this case, $M^3(\Ga_0)$ admits a finite cover which is homeomorphic to
the product of $\S^1$ and a surface.

\begin{figure}[tbh]
\begin{center}
\input{knot.pstex_t}
\end{center}
 \caption{}

 \label{knot.fig}
\end{figure}

\begin{thm}
(B.\! Apanasov \cite{Apanasov},  B.\! Maskit \cite{Maskit(1987)}, see also \cite{Gromov-Lawson-Thurston(1988)}.)
There are 1-\qf groups $\Ga \subset \Mob(\S^3)$  such that $\Lambda (\Ga)$  are wild knots.
\end{thm}
\sproof ~Start with a necklace of round balls
$$
B_0, B_1 ,..., B_{m-1}\subset \S^3,$$
so that $B_i$ is tangent to $B_{j}$, if $j=i+1\in \Z/m\Z$ and is disjoint otherwise.
Assume that this necklace is {\em knotted}, i.e.,
the polygonal knot obtained by connecting the consecutive points of tangency is a nontrivial knot $K\subset \S^3$.
See Figure \ref{knot.fig}.

Let $\ga_i\in \mbt$ denote the inversion in the sphere $\partial B_i$, $i=0, 1,...,m-1$.
Let $\Ga\subset \mob$ be the group generated by these inversions. By the Poincar\'e fundamental polyhedron theorem,
$$
\Phi=\S^3 \setminus \bigcup_{i=0}^{m-1} B_i
$$
is a fundamental domain for $\Ga$, the group $\Ga$ is geometrically finite and is isomorphic to a 1-fuchsian group $\Ga'$.
Therefore,  Tukia's theorem \ref{tukia} applied to the isomorphism $\Ga\to \Ga'$,  implies that $\Ga$ is
1-quasifuchsian. By Seifert--van Kampen Theorem, $\pi_1(\S^3 \setminus K)$ embeds in $\pi_1(\Om(\Ga))$.
Therefore  $\pi_1(\Om(\Ga))$ is not isomorphic to $\Z$ and, hence, the limit set of $\Ga$
is a wild knot. \qed

\medskip
By modifying the above construction, S.~Hwang  proved

\begin{thm}\label{H1}
(S.~Hwang \cite{Hwang}.) Let $L$ be a polygonal link in $\S^3$. Then there exists a (torsion-free)
convex-cocompact Kleinian group $\Ga \subset \Mob(\S^3)$ with a fundamental domain $\Phi\subset \S^3$
such that the complement $\S^3\setminus \Phi$ is isotopic to a regular neighborhood of $L$.
\end{thm}

The above theorem is the key for proving

\begin{thm}\label{H2}
(S.~Hwang \cite{Hwang}.) Let $M^3$ be a closed oriented 3-manifold. Then there
exists a closed oriented 3-manifold $N^3$ such that the connected sum $M^3\# N^3$ admits a M\"obius structure.
\end{thm}

Results similar to Theorems \ref{H1} and \ref{H2} hold in dimension 4, see  \cite{Kapovich04}, although
one has to assume that $M^4$ is a Spin-manifold.
One of the key ingredients in \cite{Kapovich04} is the following:

\begin{thm}
Let $Q\subset \R^4\subset \S^4$ be a finite 2-dimensional subcomplex in the standard cubulation of $\R^4$.
Then there exists a convex-cocompact Kleinian subgroup $\Ga \subset \mbc$ (generated by reflections)
with a fundamental domain $\Phi\subset \S^4$,  such that the complement $\S^4\setminus \Phi$ is isotopic
to a regular neighborhood of $Q$.
\end{thm}

Very little is known about \qf groups in $\mob$ whose limit sets have dimension between $2$ and $n-2$.
Perhaps the most interesting result here is obtained by I.\! Belegradek \cite{Belegradek(1997)}
 who used the construction from \cite{Gromov-Lawson-Thurston(1988)} to get

\begin{thm}
There exist convex-cocompact 2-\qf subgroups $\Ga_1, \Ga_2 \subset \Mob(\S^4)$ so that:

1) $\Lambda(\Ga_1)$ is a wild 2-sphere in $\S^4$.

2) $\Lambda(\Ga_2)$ is tame but the group $\Ga_2$ is not topologically
conjugate to a 2-fuchsian group: $M^4(\Ga_2)$ is a
nontrivial circle bundle over a hyperbolic 3-manifold.
\end{thm}

\noindent Similar results probably hold for codimension 2 quasifuchsian subgroups in $\mob$, $n\ge 5$.

\section{Ahlfors finiteness theorem in higher dimensions:\\
Quest for the holy grail}\label{grail}

\subsection{The holy grail}

\medskip
One of the most fundamental results of the theory of Kleinian subgroups of $\Mob(\S^2)$ is
the {\em Ahlfors Finiteness Theorem} (the ``Holy Grail''), which we state here together with its companions:

\begin{thm}
\label{aft}
Suppose that $\Ga\subset PSL(2, \C)$ is a
Kleinian group\footnote{Recall that all Kleinian groups
are assumed to be finitely generated.} which may have torsion. Then the following hold:

1. (L.\! Ahlfors \cite{Ahlfors(1964)}, L.\! Greenberg \cite{Greenberg(1967)}.)
The group $\Ga$ is analytically finite, i.e., the quotient $O:=\Om(\Ga)/\Ga$ is a complex
orbifold of finite conformal type\footnote{I.e., as a Riemann surface it
biholomorphic to a compact Riemann surface with a finite subset removed;
as an orbifold it has only finitely many singular cone-points.}. In particular, $O$
is homotopy-equivalent to a finite CW complex.

2. (D. Sullivan \cite{Sullivan(1981a)}.) $\Ga$ has only finitely many cusps.

3. (M.\! Feighn and G.\! Mess \cite{Feighn-Mess}.) $\Ga$ has only finitely many
$\Ga$-conjugacy classes of finite order elements.

4. (P.\! Scott \cite{Scott(1973a), Scott(1973b)}.)
$\Ga$ is finitely presentable and the orbifold $\H^3/\Ga$ is
finitely covered by a manifold $\H^3/\Ga'$, which is homotopy-equivalent
to a compact 3-manifold.

5. (L.\! Ahlfors \cite{Ahlfors(1971)}.) The action of $\Ga$ on $\La(\Ga)$ is {\em recurrent} with
respect to the Lebesgue measure $\mu$.
\end{thm}

Alternative analytical proofs of Part 1 (i.e., the original Ahlfors'  finiteness theorem) are given for instance in
\cite[Section 8.14]{Kapovich00} and \cite{Marden(2006)}.  A geometric proof (valid even in the context of manifolds of negative curvature)
follows from the solution of Tameness Conjecture, see \cite{Agol}.

\begin{corollary}
\label{aftc}
If $\Ga$ is as above then:

a. For each component $\Om_0$ of $\Om(\Ga)$, the limit set of the stabilizer of $\Om_0$ in $\Ga$
equals $\D \Om_0$ (follows directly from Part 1 of Theorem \ref{aft}). In particular, no component of
$\Om(\Ga)$ has trivial stabilizer.

b. Kleinian subgroups $\Ga$ of $\Mob(\S^2)$ are {\em coherent}, i.e., each finitely generated
subgroup of $\Ga$ is also finitely presented (follows from Part 4 of Theorem \ref{aft}).

c. (W.\! Thurston, see \cite{Morgan(1984)}.) If $\Ga\subset \Mob(\S^2)$ is geometrically finite with
$\Om(\Ga)\ne \emptyset$ then each finitely generated subgroup
$\Del\subset \Ga$ is geometrically finite as well.
\end{corollary}

We also now fully understand the topology of the manifold (orbifold) $\H^3/\Ga$:

\begin{thm}\label{tamenessconjecture}
(Former tameness conjecture.) The quotient $\H^3/\Ga$ is {\em tame}, i.e., it
is homeomorphic to the interior of a compact manifold (orbifold) with boundary.
\end{thm}

The above theorem was proved for freely indecomposable groups $\Ga$ by F.~Bonahon \cite{Bonahon(1986)}
and by I.~Agol \cite{Agol}, D.~Calegari and D.~Gabai \cite{Calegari-Gabai} in the general case.

\medskip
The next theorem is a combination of a result by Thurston \cite{Thurston(1978-81)},
who proved ergodicity for tame Kleinian subgroups of $\mbd$,
and the proof of the tameness conjecture:

\begin{thm}\label{ergodic}
If $\Ga$ is as above, then the action of $\Ga$ on $\La(\Ga)$ is ergodic with respect to the
Lebesgue measure: each measurable $\Ga$-invariant function on $\La(\Ga)$ is constant a.e..
\end{thm}

Note, that the conglomerate of assertions presented above contains statements
of different nature: algebraic, topological, dynamical.
For a while it was hoped that a theorem analogous to Theorem \ref{aft}
can be proved for Kleinian groups in higher dimensions; an attempt
to develop  analytical technique to achieve this was made by Ahlfors
in \cite{Ahlfors(1981)} (see also \cite{Ohtake}).

\smallskip
Nearly all algebraic and topological assertions of Theorem \ref{aft} and the Corollary \ref{aftc}
have been disproved in the case of Kleinian
groups acting in higher dimensions
(M.\ Kapovich and L.\ Potyagailo \cite{Kapovich(1995a)}, \cite{KP2}, \cite{KP1},
\cite{Potyagailo(1990)}, \cite{Potyagailo(1994)}):

\begin{thm}\label{fail}
There exist Kleinian subgroups $\Ga_1,...,\Ga_5\subset \mbt$ so that:

1. The group $\Ga_1$ is not finitely presentable.

2. For each $i$, the manifold $M(\Ga_i)=\Om(\Ga_i)/\Ga_i$  contains a component with infinitely generated
fundamental group.

3. $\Ga_2$ is free and has infinitely many cusps (of rank 1).

4. $\Ga_{3}$ is not torsion-free and has infinitely many conjugacy classes of finite order elements.

5. $\Ga_4$ is a normal subgroup of a convex-cocompact group
$\widehat\Ga_4 \subset \Mob({\S^3})$ and satisfies (1), (2) and (4).

6. (B.\! Bowditch, G.\! Mess \cite{Bowditch-Mess}.) The group $\Ga_5$ satisfies
(1) and (2) and is contained in a cocompact lattice $\widehat\Ga_5\subset \mbt$.

7. Groups $\Ga_i$, $i=1,...,4$ are normal subgroups of geometrically finite groups $\widehat\Ga_i$ so
that $\widehat\Ga_i/\Ga_i\cong \Z$.
\end{thm}

\begin{remark}
By modifying $\Ga_3$ one can also construct an example $\Ga_{6} \subset \mbt$
such that $\Omega(\Ga_6)/\Ga_6$ has infinitely many connected components.
\end{remark}

At the time when the above examples were constructed, they were regarded as a ``rare pathology''. It appears however
that such examples are rather common:

\begin{conjecture}
(M.~Kapovich, L.~Potyagailo, E.~Vinberg \cite{KPV}.)
Suppose that $\Ga\subset \mob$ is an arithmetic lattice, where $n\ge 3$. Then $\Ga$ is {\em noncoherent}, i.e.,
it contains a finitely generated subgroup $\Del$ which is not finitely presentable.
\end{conjecture}

This conjecture was proved in  \cite{KPV} in a number of special cases, e.g. for all non-cocompact arithmetic lattices provided that
$n\ge 5$.

\medskip
All the examples $\Ga_i$ in the above theorem are based upon the existence of hyperbolic
3-manifolds $M^3$ of finite volume which fiber over the circle:
the groups $\Ga_i$ are obtained by manipulating with the normal surface subgroups
in $\pi_1(M^3)$.

\begin{problem}
Find examples similar to $\Ga_i$'s without using hyperbolic 3-manifolds fibering over the circle.
\end{problem}

\begin{problem}
Construct a finitely generated Kleinian subgroup $\Ga\subset \mob$ such that Part (a) of Corollary \ref{aftc}
fails for $\Ga$.
\end{problem}

\begin{problem}
(G.~Mess.) Construct a finitely-presented Kleinian subgroup of $\mob$ ($n\ge 3$) which contains no parabolic elements
and for which any of the assertions of Theorem \ref{aft} fail.
(In Part (a) one would need to replace  analytical finiteness
with finiteness of the homotopy type.)
\end{problem}

\begin{problem}
Construct a finitely generated Kleinian subgroup $\Ga\subset \mbt$ such that $\Om(\Ga)$ contains a contractible
component $\Om_0$ so that:

The stabilizer $\Ga_0$ of $\Om_0$ in $\Ga$  is  finitely-generated, but the manifold $\Om_0/\Ga_0$ is not tame.
\end{problem}

Note however that although algebra and topology fail,
the assertions of dynamical nature
(part 5 of Theorem \ref{aft}, part (a) of Corollary \ref{aftc}, and Theorem \ref{ergodic})
remain open in higher dimensions.
Moreover, an attempt  to construct a higher-dimensional counter-example to Theorem \ref{aft} (part 5) along the
lines of the examples $\Ga_i$, is doomed to failure:

\begin{thm}
(K.\ Matsuzaki \cite{Matsuzaki}.) Let $\widehat\Ga$ be a geometrically finite
subgroup\footnote{Actually, the proof also works for subgroups of any rank 1 Lie group.}
in $\mob$. Suppose that $\Ga\subset \widehat\Ga$ is a normal subgroup (which does not have to be finitely
generated). Then the action of $\Ga$ on its limit set is recurrent.
\end{thm}

Ergodicity fails however for discrete subgroups of $PU(2,1)$ (it probably also fails for Kleinian
groups in higher dimensions but an example would be difficult to construct):

\begin{thm}
There exists a finitely generated (but not finitely presentable!) discrete group $\Ga$ of isometries
of complex-hyperbolic 2-plane $\C\H^2$ so that the limit set of $\Ga$ is the 3-sphere and the action
of $\Ga$ on $\S^3$ is not ergodic.
\end{thm}
\proof There are examples (the first was constructed by R.\! Livne in his thesis \cite{Livne}, see also
\cite{Deligne-Mostow}) of cocompact torsion-free discrete subgroups $\widehat\Ga\subset PU(2,1)$ such that
the complex 2-manifold $M=\C\H^2/\widehat\Ga$ admits a nonconstant holomorphic map $f: M\to S$
to a Riemann surface $S$ of genus $\ge 2$. The fundamental group of the generic fiber of $f$ maps onto
a normal subgroup $\Ga$ in $\widehat\Ga$, so that $\Ga$ is finitely generated but is not
finitely presentable \cite{Kapovich(1998a)}. By lifting $f$ to the universal covers we get a
nonconstant holomorphic $\Ga$-invariant function
$$
\tilde{f}: \C\H^2\to \H^2.
$$
Then the bounded harmonic function $Re(\tilde f)$ is also $\Ga$-invariant and
 nonconstant. This harmonic function admits a measurable extension $h$  to $\S^3$,
 the boundary of the complex ball $\C\H^2$, so that  $h$  is $\Ga$-invariant and not a.e. constant. \qed

\subsection{Groups with small limit sets}
\label{sec:small}

So far, our quest for the holy grail mostly resembled Monty Python's:
We are not sure what to look for in higher dimensions. Nevertheless,
there is a glimmer of hope.

Recall that the Hausdorff dimension $\dim_H$ of a subset
$E\subset \R^n$ is defined as follows.
For each $\al>0$ consider the $\al$-Hausdorff measure of $E$:
$$
mes_\al(E)= \lim_{\rho\to 0} \inf\{ \sum_i r_i^{\al}:\quad r_i \le \rho, E
\hbox{~is contained in the union of $r_i$-balls} \}.
$$
The Hausdorff dimension of $E$ is
$$
\dim_H(E)=\inf_{\al} \{ \al : mes_\al(E)=0\}.
$$
According to \cite{Hurewitz-Wallman}, for every bounded subset
$E\subset \R^n$ one has the inequality
$$
\dim(E)\le \dim_H(E)
$$
between topological and Hausdorff dimensions.
In particular, if $\Ga$ is a Kleinian group with $\dim_H(\La(\Ga))< 1$, then $\Ga$ is
geometrically finite and is isomorphic to a Schottky-type group, see Theorem \ref{gsg}.

\begin{conjecture}
If $\dim_H(\La(\Ga))< 1$, then $\Ga$ is a Schottky-type group. Moreover, $\Ga$ is classical.
\end{conjecture}

The {\em critical exponent}  of a Kleinian group $\Ga\subset \mob$ is
$$
\del(\Ga):=\inf \{ s>0 : \sum_{\ga\in \Ga} e^{-s d(x, \ga x)} <\infty\},
$$
where $d$ is the hyperbolic metric in $\H^{n+1}$. The following theorem is the result of combined efforts of a
large number of mathematicians, including P.~Tukia, D.~Sullivan and P.~Nicholls, we refer to
\cite{Nicholls}, \cite{Bishop-Jones} for the proofs:

\begin{thm}
For every Kleinian subgroup $\Ga\subset \mob$,
$$
\delta(\Ga)= \dim_H(\La_c(\Ga)).
$$
\end{thm}

Recall that $\La_c(\Ga)$ is the conical limit set of $\Ga$.

The critical exponent relates to $\la_0$, the bottom of the spectrum of Laplacian  on the hyperbolic manifold
$\H^{n+1}/\Ga$, by the following

\begin{thm}
(D. Sullivan \cite{Sullivan(1987)})
$$
\la_0 =\left(\frac{n}{2}\right)^2, \hbox{~~if~~}\quad 0\le \del(\Ga)\le \frac{n}{2},
$$
$$
\la_0 = \del(\Ga)(n-  \del(\Ga)), \hbox{~~if~~}\quad \frac{n}{2}< \del(\Ga)\le n.
$$
\end{thm}

The expectation is that Kleinian groups in $\mob$ with {\em small} limit sets behave analogously to the Kleinian subgroups
of $\mbd$.

\begin{conjecture}
\label{C2}
Suppose that $\Ga$ is a (finitely generated) subgroup of $\mob$ so that $\La(\Ga)$
has Hausdorff dimension $< 2$. Then $\Ga$ is geometrically finite.
\end{conjecture}

For $n=2$, this conjecture is a theorem of C.~Bishop and P.~Jones \cite{Bishop-Jones}.
A partial generalization of \cite{Bishop-Jones} was proved by A.~Chang, J.~Qing, J. and P.~Yang \cite{CQY}:

\begin{thm}
Suppose that $\Ga$ is a (finitely generated) {\em conformally finite}\footnote{I.e.,
$M^n(\Ga)=\Om(\Ga)/\Ga$ is compact modulo cusps.}
subgroup of $\mob$ such that $\dim_H(\La(\Ga))<n$. Then $\Ga$ is geometrically finite.
\end{thm}

The converse to the above theorem was proved earlier by P.~Tukia \cite{Tukia(1984)}.

\begin{thm}
(Y.\! Shalom \cite{Shalom}.)
Suppose that $\Ga$ is a geometrically finite subgroup of $\mob$ such that $\dim_H(\La(\Ga))<2$
and $\Del\subset \Ga$ is a finitely generated normal subgroup.
Then $\Del$ has finite index in $\Ga$.
In particular, $\Del$ is geometrically finite as well.
\end{thm}

Thus, attempts to construct geometrically infinite groups using normal subgrooups in geometrically finite
Kleinian groups with {\em small} limit sets, are doomed to failure.
On the other hand, the assumption that $\delta(\Ga)$ is small should impose strong
restrictions on the algebraic properties of the group $\Ga$.

\begin{conjecture}\label{C3}
Suppose that $\Ga$ is a Kleinian group in $\mob$ which does not contain
parabolic elements. Then:

1. $cd(\Ga)-1\le \delta(\Ga)$.

2. In the case of equality, $\Ga$ is an $i$-fuchsian convex-cocompact group, $i=\delta(\Ga)$.
\end{conjecture}

Recall that $cd$ and $hd$ stand for the cohomological and homological dimensions of a group.
A partial confirmation of Part 1 of this conjecture is obtained in

\begin{thm}\label{inequality}
(M. Kapovich \cite{Kapovich07}.) Suppose that  $\Ga\subset \mob$
is a Kleinian group. Then for every ring $R$,
$$
hd_R(\Ga, \Pi) -1\le \delta(\Ga),
$$
where $\Pi\subset \Ga$ is the set of virtually abelian
subgroups of $\Ga$ of (virtual) rank $\ge 2$ and
$hd_R(\Ga, \Pi)$ is the relative homological dimension.
\end{thm}

We refer the reader to the series of papers by H.~Izeki \cite{Izeki(1995), Izeki(2000), Izeki(2002)}
for the related results.

\begin{corollary}
(M. Kapovich \cite{Kapovich07}.) Suppose that the group $\Ga$ is
finitely-presented and $\del(\Ga)<1$. Then $\Ga$ is free.
\end{corollary}
\proof Since $\del(\Ga)<1$, it follows that $\Ga$ contains no rank 2 abelian subgroups.
Then we have the inequalities
$$
cd(\Ga)\le 1+ hd(\Ga)\le \delta(\Ga)+2<3.
$$
Combined with finite presentability of $\Ga$, the inequality $cd(\Ga)\le 2$ implies that
$\Ga$ has {\em finite type}; therefore
$$
cd(\Ga)=hd(\Ga)\le  \delta(\Ga)+1<2.
$$
Hence $\Ga$ is free by Stallings' Theorem \ref{stalling}. \qed

\medskip
An  inequality similar to Conjecture \ref{C3} was proved by A.~Reznikov: For a (finitely-generated) group $\Ga$ define
$$
\al(\Ga):= \inf\{p\in [1,\infty] : \ell_p H^1(\Ga)\ne 0\}.
$$
Here $\ell_p H^1$ is the first $\ell_p$-cohomology of the group $\Ga$, see \cite{Bourdon-Martin-Valette} for
the precise definition. Then

\begin{thm}
(A.~Reznikov \cite{Reznikov(1999)}, see also \cite{Bourdon-Martin-Valette} for the detailed proof in the case
of isometries of $CAT(-1)$ spaces.) For every Kleinian group $\Ga\subset \mob$,
$$
\al(\Ga)\le \max(\del(\Ga), 1).
$$
\end{thm}

\begin{question}
What can be said about $\Ga$ in the case of equality in Reznikov's theorem?
\end{question}

In the case of geometrically finite groups, Part 2 of Conjecture \ref{C3} holds:

\begin{thm}
1. (Chenbo Yue \cite{Yue}, see also \cite{BCG03} and \cite{Bonk-Kleiner01}.)
Suppose that $\Ga$ is convex-cocompact and
$i=\del(\Ga)=cd(\Ga)-1$. Then $\Ga$ is $i$-fuchsian.

2. (M. Kapovich \cite{Kapovich07}.) Suppose that $\Ga$ is geometrically finite.
Then the following three conditions are equivalent:

$i=\del(\Ga)=cd(\Ga, \Pi)-1$,

$i=\dim(\La(\Ga))=\dim_H(\La(\Ga))$,

$\Ga$ is $i$-fuchsian.
\end{thm}

%\begin{rem}
%Chenbo Yue states his theorem in the context of $i$-\qf groups, but his proof actually does
%not need this assumption.
%\end{rem}

\begin{conjecture}
Suppose that $\Ga$ is a Kleinian group in $\mob$ whose limit set is not totally disconnected and
has Hausdorff dimension $1$. Then $\Ga$ is $1$-fuchsian.
\end{conjecture}

This conjecture is known to be true for $n=2$, see \cite{Canary-Taylor}.

\begin{problem}
(The {\em gap} problem, L.~Bowen, cf. \cite{Stratmann}.)
1. Is there a number $d_n<n$ such that for every Schottky subgroup $\Ga\subset \mob$, $n\ge 3$, we
have:
$$
\del(\Ga)< d_n.
$$
2. More generally, consider a sequence $\Ga_j\subset \mob$ of convex-cocompact groups isomorphic to a fixed group $\Ga$ so that:
$\La(\Ga_j)\ne \S^n$ for each $j$. Is it true that
$$
\lim \sup_{j\to\infty} \del(\Ga_j)<n~ ?
$$
\end{problem}

By the work of R.~Phillips and P.~Sarnak \cite{Phillips-Sarnak}, the answer to the Part 1 of this question is positive in the class of
classical Schottky groups.

\section{Representation varieties of Kleinian groups}
\label{varieties}

For a finitely-generated $\Ga$ consider the {\em representation variety} of $\Ga$:
$$
R_n(\Ga):=Hom (\Ga, \mob).$$
If $\Ga$ has the presentation
$$
\Ga=\<x_1,...,x_m| r_1,..., r_k,...\>,
$$
the representation variety is given by
$$
\{ (g_1,...,g_m)\in (\mob)^m: r_1(g_1,...,g_m)=1, ... r_k(g_1,...,g_m)=1,....\}.
$$

The group $\mob$ acts on $R_n(\Ga)$ via conjugation:
$$
\theta\cdot \rho(\ga)=\theta \rho(\ga)\theta^{-1}, \quad \theta\in \mob.
$$
Given this action, one can form the quotient variety
$$
X_n(\Ga):= R_n(\Ga)/\!/\mob,
$$
called the {\em character variety}. Roughly speaking, the elements of $X_n(\Ga)$ are represented
by conjugacy classes of representations $\rho:\Ga\to \mob$. This is literally true for
``most'' representations, the ones for which $\rho(\Ga)$ does not contain a normal parabolic
subgroup, see \cite{Johnson-Millson}.
%\footnote{Such $\rho$'s are called {\em reductive}.}),
In general, the representations $\rho_1, \rho_2$ project to the same point in $X_n(\Ga)$
iff the closures of their $\mob$-orbits have nonempty intersection.  We let $[\rho]$ denote the projection of
$\rho\in R_n(\Ga)$ to $X_n(\Ga)$.

A {\em trivial deformation} of a representation $\rho_0\in R_n(\Ga)$ is a connected curve
$\rho_t\in R_n(\Ga)$ which projects to a point in $X_n(\Ga)$. A representation $\rho_0$
is called {\em rigid} if it admits no nontrivial deformations.

We will be mostly interested in representations $\rho\in R_n(\Ga)$ which have discrete, nonelementary image,
however much of our discussion is more general.

\medskip
In this section we address the following issues related to the character varieties:

i. Local structure of $X_n(\Ga)$ and existence of small deformations of a given Kleinian group (rigidity vs. flexibility).

ii. Connectedness of the subspace ${\mathcal D}_n(\Ga)$ of discrete and faithful representations in $X_n(\Ga)$.

iii. Structural stability: What happens to a Kleinian group in $\mob$ if we deform it a little bit? Does it stay Kleinian?

iv. Compactness of  ${\mathcal D}_n(\Ga)$ and estimates on various natural continuous functionals on ${\mathcal D}_n(\Ga)$.

v. Difficulties in constructing ``truly higher-dimensional'' geometrically infinite Kleinian groups.

\subsection{Local theory}
 \label{local}

We start by considering the {\em local} structure of $X_n(\Ga)$. Given an abstract group $\Ga$ and a representation
$\rho\in R_n(\Ga)$, we  have the adjoint action of $\rho(\Ga)$
on the Lie algebra ${\mathfrak g}$ of $\mob$ and  the associated first cohomology group
$$
H^1(\Ga, Ad(\rho))= Z^1(\Ga, Ad(\rho))/B^1(\Ga, Ad(\rho)),
$$
see Section \ref{homology}. It was first observed by A.~Weil \cite{Weil} (in the general context of representations to Lie groups)
that if $X_n(\Ga)$ is smooth at $[\rho]\in X_n(\Ga)$ then $H^1(\Ga, Ad(\rho))$ is
isomorphic to the tangent space to $X_n(\Ga)$ at $[\rho]$. Moreover, Weil proved that if $H^1(\Ga, Ad(\rho))=0$
then $[\rho]$ is an isolated point on $X_n(\Ga)$, i.e.,
$\rho$ is rigid.

Therefore, the elements of $H^1(\Ga, Ad(\rho))$ can be regarded as {\em infinitesimal deformations} of
the representation $\rho$. An infinitesimal deformation $[\xi]\in H^1(\Ga, Ad(\rho))$ is called {\em integrable} if it is
tangent to a smooth curve in $X_n(\Ga)$. The obstructions to integrability are cohomological in nature, they are certain elements of
$H^2(\Ga, Ad(\rho))$, called {\em Massey products}. However, in practice,
these cohomology classes are very difficult to compute. The first such obstruction is the {\em cup-product}:
$$
\phi([\xi])= [\xi] \cup [\xi] \in H^2(\Ga, Ad(\rho)),
$$
see for instance \cite{Goldman-Millson(1988)}. Here $\phi([\xi])$ is represented by
the 2-cocycle
$$
\tau(x,y)= [\xi(x), Ad\circ\rho(x)\xi(y)],
$$
where $[,]$ is the Lie bracket on the Lie algebra ${\mathfrak g}$.
If the first obstruction vanishes and $\Ga$ is the fundamental group of
a surface, then $X_n(\Ga)$ is smooth at $\rho$, see \cite{Goldman-Millson(1988)},
where a much more general result is proved.

\medskip
We will be mostly interested in the case where $\rho: \Ga\hook \mob$ is a discrete embedding, whose image we will
identify with $\Ga$. Then, by abusing the terminology, we will talk of small deformations of $\rho$ in $R_n(\Ga)$ as
{\em small deformations of} $\Ga$ itself.

\subsubsection{Small deformations of $1$-quasifuchsian groups}

Recall that the group $\mob$ has dimension $d=(n+2)(n + 1)/2$. Suppose that $\Ga\subset \mob$ is
$1$-quasifuchsian; in this subsection we allow $\Ga$  to have nontrivial finite order elements. We assume however that $\Ga$
contains no elements fixing the circle $C=\La(\Ga)$ pointwise. Therefore we obtain the injective map
$$
\Ga\to \Isom(\H^2)
$$
given by the restriction of the elements of $\Ga$ to the round circle $C$. To simplify the discussion we assume that $\Ga$
preserves the orientation on $C$. Then $\Ga$ embeds as a lattice in $PSL(2,\R)$.

If $\Ga$ contains no parabolic elements then it has the presentation:
$$
\< a_1 , b_1 ,...., a_q , b_q , c_1 ,..., c_k |
[a_1 , b_1 ]\cdot ... [a_q , b_q ] \cdot c_1 \cdot...\cdot c_k = 1 ,
c_j^{r_j} = 1 , j= 1 ,..., k \>.
$$
For a representation $\rho: \Ga\to \mob$ we let
$$
e_j:=d-\dim\{ \xi \in {\mathfrak g} : Ad\circ\rho(c_j)(\xi )= \xi \};$$
in other words, $e_j$ is the codimension of the centralizer of $\rho(c_j)$ in $\mob$. Let $s$ denote
the dimension of the centralizer of $\rho(\Gamma)$ in $\mob$.

\begin{thm}
(A.\! Weil \cite{Weil}.)
\begin{equation}
\label{wei}
h = \dim H^1(\Gamma, Ad(\rho)) = (2q-2)d + 2s + e_1 +...+ e_k.
\end{equation}
Moreover, if $s = 0$ then $X_n(\Ga)$ near $[\rho]$ is a smooth $h$-dimensional
 manifold.
\end{thm}

For instance, if $n=1$, $\Ga\subset PSL(2,\R)\subset \Mob(\S^1)$;
therefore $d=3$, we get $e_i=1$ for each $i=1,...,k$, $s=0$.  Hence
$$
h= 6q-6 +k,
$$
which is the familiar formula for the dimension of the Teichm\"uller space of the orbifold $O=\H^2/\Ga$.
If $n=2$, we, of course, obtain $h=2(6q-6 +k)$ which is the (real) dimension of the
space of the quasifuchsian deformations of $\Ga$ in $PSL(2,\C)$.

To better understand the difficulties which one encounters in
the case of $i$-fuchsian groups for $i\ge 2$, we consider the {\em hyperbolic triangle groups} $\Ga$.
The reason for considering these groups is that they are rigid in $PSL(2, \R)$
(similarly to rigidity of lattices in $\mob$, $n\ge 2$).

The triangle groups are the $1$-fuchsian groups with $q=0$, $k=3$; every such $\Ga$
has the presentation
$$
\< c_1 , c_2 , c_3 | c_1 \cdot c_2 \cdot c_3 = 1, c_j^{r_j} = 1 , j= 1,2, 3 \>,$$
where $r_1^{-1}+r_2^{-1}+r_3^{-1}<1$. Such group embeds discretely into $PSL(2,\R)$
and we will denote the image of this embedding by $\Delta=\Delta(r_1,r_2,r_3)$

As a subgroup of $\Mob(\S^2)$, the group $\Delta$ is rigid
(which follows from vanishing of $H^1$). Moreover, triangle groups are ``strongly rigid" in $\Mob(\S^2)$,
i.e., every discrete embedding of $\Gamma$ into $\Mob(\S^2)$ is induced by conjugation of the
identity embedding, see \cite{Greenberg(1981)} for the complete description of $X_2(\Ga)$.

The situation changes somewhat if we consider representations into
$\Mob(\S^3)$.  First, let $\rho_0: \Delta\to \Gamma' \subset \Mob(\S^3)$ be
the embedding obtained as the composition
$$
\Delta\subset \Mob(\S^1)\hook  \Mob(\S^2)\hook \Mob(\S^3).
$$
of natural embeddings. Then $\dim H^1(\Delta, Ad(\rho_0))=0$ and hence $\rho_0$ is still rigid in $\Mob(\S^3)$.
The easiest way to  see this is to use Weil's formula (\ref{wei}):
$$
d = 10, s = 1, e_j = 6, \quad \hbox{~~for~~} j= 1, 2, 3$$
and hence
$$
h= -20+ 2+6+6+6=0.
$$
However, instead of $\rho_0$ we can take a {\em twisted} extension.
 Suppose that we can find numbers $m_j \in \Z, 1<|m_j|<r_j-1$ ($j= 1, 2, 3 $) such that:
$$
m_1^{-1} + m_2^{-1} + m_3^{-1} = 0, \hbox{~and~} \forall j, ~~m_j \hbox{~~divides~~} r_j.
$$
(This is satisfied for instance by $m_1=m_2=4, m_3=-2$ and $r_j=8$ for all $j$.)

Define a homomorphism $\theta:\Delta\to SO(2)$ by sending $c_i$ to the rotation by $2\pi/m_i$.
Then define $\rho: \Delta\to \Mob(\S^1)\times SO(2) \subset \Mob(\S^3)$ by twisting $\rho_0$
via $\theta$:
$$
\rho(\ga)=\rho_0(\ga)\times \theta(\ga), \quad \ga\in \Del.
$$

It is clear that $\rho(\Delta)$ is again a 1-fuchsian subgroup in $\Mob(\S^3)$.
If $r_j > 3$ for each $j$, then $e_j = 8$, $s=1$ and the
formula (\ref{wei}) gives the dimension $h=6$ for $H^1(\Delta, Ad(\rho))$.
I do not know if any of these {\em infinitesimal} deformations is integrable.
To decide this one has to analyze the quadratic form
$$
\phi: H^1(\Delta, Ad(\rho)) \to  H^2(\Delta, Ad(\rho))=\R,
$$
given by the cup-product.   According to \cite{Goldman-Millson(1988)}, the quadratic cone $\{\phi=0\}$ is analytically isomorphic to
a neighborhood of $[\rho]$ in $X_3(\Delta)$; hence it suffices to find a nontrivial 1-cocycle $\xi$ for
which $\phi([\xi])=0$ to get nontrivial deformations of the representation $\rho$.

On the other hand, one can use (\ref{wei}) to show that every representation $\rho$ of the group
$\Delta=\Delta(2,3,r_3)$ into $\Mob(\S^3)$ has zero cohomology $H^1(\Delta, Ad(\rho))$. Therefore $X_3(\Del)$ is
a zero-dimensional algebraic variety and, hence, is a finite set.

This situation is somewhat typical for representations of lattices in $\mob$ ($n\ge 2$) into $\MOB$:
In a number of cases we can prove rigidity by making cohomological computations; in some cases we
can only conclude that $H^1$ is nonzero, without being able to make a definitive conclusion about existence of
nontrivial deformations.

\subsubsection{Small deformations of $i$-\qf groups for $i\ge 2$}
\label{small}

In the case of $(n-1)$-quasi\-fuchsian groups $\Ga$ ($n\ge 2$),
the existence of nontrivial deformations of $\Ga$ in $\mob$ is not
at all clear. Suppose that $\Ga\subset \Isom(\H^{n})$ is a
cocompact lattice. Then  the identity embedding
$$
\iota: \Ga\hook \Isom(\H^{n}),
$$
is rigid by Mostow's theorem.

\begin{rem}
Actually, (local) rigidity of $\iota$ was known prior to the work of Mostow; it was first established by E.~Calabi
\cite{Calabi}, whose proof was later generalized by A.~Weil \cite{Weil(1960), Weil(1962)}.  These arguments were
based on proving that $H^1(\Ga, Ad(\iota))=0$.
\end{rem}

Consider now the composition of $\iota$ with the natural embedding:
$$
\rho: \Ga \to \Isom(\H^n)\hook \Isom(\H^{n+1}).
$$
Then
$$
H^1(\Ga, Ad(\iota))\cong H^1(\Ga, V_n),
$$
where $V_n=\R^{n,1}$ and $\Ga$ acts on the Lorentz space $\R^{n,1}$ via the usual embedding $\Ga\hook O(n,1)$.

It turns out that $\rho$ may or may not be rigid, even if $\Ga$ is torsion-free: Rigidity depends on the
lattice $\Ga$. One has the following list ((a) through (e))
of constructions of deformations and infinitesimal deformations of $[\rho]$ in $X_n(\Ga)$:

\smallskip
(a) {\em Bending}, see \cite{Johnson-Millson}, \cite{Kourouniotis}. Given a connected properly embedded totally-geodesic
hypersurface $S\subset M=\H^n/\Ga$, one associates with $S$ a smooth curve through $[\rho]$ in $X_n(\Ga)$,
called the {\em bending deformation} of $[\rho]$. More generally, given a disjoint collections of such
hypersurfaces $S_1,...,S_k$, one obtains a $k$-dimensional smooth submanifold in $X_n(\Ga)$ containing $[\rho]$. This
construction is completely analogous to bending deformations of 1-fuchsian  subgroups in $\mbd$ defined by W.~Thurston
in \cite{Thurston(1978-81)}.  We let $[\xi_S]$ denote the element of $H^1(\Ga, Ad(\rho))$ corresponding to
a connected totally-geodesic hypersurface $S\subset M$.

\medskip
There are numerous  groups $\Ga$ satisfying assumptions of the bending construction.
Namely, start with an arithmetic group $O'(f,A)$ of
the simplest type (see Section \ref{basics}), where
$$
f= a_0 x_0^2+ a_1 x_1^2+ ... + a_n x_n^2,
$$
and $a_0<0, a_i>0, i= 1,...,n$. Identify $\H^{n}$ with a component of the hyperboloid $\{f(x)=-1\}$.
Then the stabilizer of the hyperplane $P=\{x_n=0\}$ in $O'(f,A)$ is an arithmetic lattice in
$\Isom(\H^{n-1})$. The intersection $H=P\cap \H^{n}$ is a hyperplane in $\H^n$.
After passing to an appropriate  finite-index subgroup $\Ga$ in $O'(f,A)$, one obtains  a totally-geodesic embedding
of the hypersurface $S=H/\Ga'$ into $H/\Ga$, where $\Ga'=\Ga\cap O'(f,A)$. We refer the reader to \cite{Millson(1976)} for the
details.

\begin{problem}
[I. Rivin] Construct examples of hyperbolic  $n$-manifolds $M$ of finite volume ($n\ge 4$) such that $M$ contains a
{\em separating} properly embedded totally-geodesic hypersurface $S\subset M$. Note that the main objective of \cite{Millson(1976)}
was to construct {\em nonseparating} hypersurfaces.
\end{problem}

The idea of bending deformations of representations is quite simple and has nothing to do with the hyperbolic
space. Below is a general description of bending as defined by D.~Johnson and J.~Millson in \cite{Johnson-Millson}.
Suppose that we are given a graph of groups ${\mathcal G}$ with the vertex groups $\Ga_v$ and the edge groups $\Ga_e$.
Let $\Ga=\pi_1({\mathcal G})$ be the fundamental group of  ${\mathcal G}$. For instance, the amalgam
\begin{equation}
\label{amal}
\Ga= \Ga_{v_1} * _{\Ga_e} \Ga_{v_2}
\end{equation}
is the fundamental group of a graph of groups which is a single edge $e$ with two vertices $v_1, v_2$.
Let $\rho_0: \Ga \to G$ be a representation of $\Ga$ to a Lie group $G$. A {\em bending deformation} of $\rho_0$
is a curve of representations $\rho_t: \Ga\to G, t\in [-1, 1]$, such that for each vertex group $\Ga_v$, we have
$$
\rho_t|\Ga_{v} = g_{v,t} (\rho_0|\Ga_{v}) g_{v,t}^{-1},
$$
for some curve $g_{v,t}$ of elements of $G$.

Therefore, the restriction of $\rho_t$ to each vertex group
determines a trivial deformation of the representation of this group. The trick is that the deformation of the representation of the
entire group $\Ga$ may be still nontrivial. For instance, in the case of the amalgam
(\ref{amal}), let $g_t\in G, g_0=1,$ be a curve of elements centralizing $\rho(\Ga_e)$, but not $\rho(\Ga_{v_1}),  \rho( \Ga_{v_2})$.
Define the family of representations
$$
\rho_t: \Ga\to G, \quad \rho_t|\Ga_{v_1}=\rho_0|\Ga_{v_1}, \quad \rho_t|\Ga_{v_2}= g_t(\rho_0|\Ga_{v_1})g_t^{-1}.
$$
In the case of the HNN extension
$$
\Ga= \Ga_{v_1} * _{\Ga_e}
$$
generated by $\Ga_{v_1}$ and $\tau\in \Ga$ such $\tau \Ga_e
\tau^{-1}=\Ga_e'\subset \Ga_{v_1}$, we take
$$
\rho_t|\Ga_{v_1}=\rho_0|\Ga_{v_1}, \quad  \rho_t(\tau) = \rho_0(\tau) g_{t}.
$$
This is a nontrivial deformation of the representation $\rho_0$.

We now return to the case when $\Ga=\pi_1(M)$, $M$ is a hyperbolic $n$-manifold containing pairwise
disjoint totally geodesic hypersurfaces $S_i, i=1,...,k$. Then
the group $\Ga$ splits as the graph of groups ${\mathcal G}$, so that the vertex subgroups $\Ga_{v_j}$
are the fundamental groups of the components of $M\setminus (S_1\cup... \cup S_k)$ and the edge groups $\Ga_{e_i}$
are the fundamental groups $\pi_1(S_i)$. Therefore

1.  The centralizer of each $\Ga_{e_i}=\pi_1(S_i)$ in $\mob$ is 1-dimensional (the group of elliptic
rotations around the limit set of $\Ga_{e_i}$).

2. The centralizer in $\mob$ of the fundamental group of each $\Ga_{v_j}$ is
zero-dimensional (i.e., $\Z_2$).

Hence one obtains nontrivial bending deformations $\rho_t$ of the
identity embedding of $\rho: \Ga\hook \mob$. The set of bending
parameters $t=(t_1,...,t_k)$ can be identified with $(\S^1)^k$, as
the centralizer of each $\Ga_{e_i}$ in $\mob$ is the circle
$SO(2)$.

\begin{thm}\label{jm}
(D.~Johnson and J.~Millson \cite{Johnson-Millson}.) For every
$n\ge 4$ there exists a uniform lattice $\Ga\subset \Isom(\H^n)$
and intersecting hypersurfaces $S_1, S_2\subset \H^n/\Ga$, so that
$$
[\xi_{S_1}]\cup [\xi_{S_2}]\in H^2(\Ga, Ad(\rho))\ne 0.
$$
In particular, $X_n(\Ga)$ is not smooth at $[\rho]$.
\end{thm}

In contrast with this result, the character varieties $X_2(\Ga)$
tend to be smooth:

\begin{thm}
(M. Kapovich \cite[Theorem 8.44]{Kapovich00}.) Let $\Ga\subset \mbd$
be a discrete subgroup. Then the identity embedding $\iota: \Ga\to
\mbd$ determines a smooth point on $X_2(\Ga)$.
\end{thm}

On the other hand, there are cocompact lattices $\Ga\subset \mbd$ and (nondiscrete) representations
$\rho: \Ga\to \mbd$ for which $X_2(\Ga)$ has nonquadratic singularity at $[\rho]$, see \cite{Kapovich-Millson(1996b)}.

\smallskip
(b) {\em Generalized bending} associated with a collection of compact
totally-geodesic submanifolds with boundary in $M^n$,  see \cite{Apanasov-Tetenov94}\footnote{Some of the theorems
stated in this paper are probably
incorrect since they do not take into account the restrictions on the angles between
the totally-geodesic submanifolds.}, \cite{Kapovich-Millson(1996)}, \cite{Maubon}, \cite{Bart-Scannell}.

The idea of the generalized bending is that instead of considering fundamental groups of
graphs of groups, one looks at the more general
{\em complexes of groups}. The only examples which had been worked out are 2-dimensional complexes of groups. Let
${\mathcal X}$ be such a complex with the vertex groups $\Ga_v$. Let $\pi_1({\mathcal X})=\Ga$
and $\rho_0: \Ga\to G$   be a representation to a Lie group.
Then, as in the definition of bending, a {\em generalized bending} of $\rho_0$ is a curve of representations
$\rho_t: \Ga\to G, t\in [-1, 1]$, whose restrictions to each vertex subgroup $\Ga_v$ are trivial deformations of $\rho_0|\Ga_v$.

\smallskip
(c) Suppose that a lattice $\Ga\subset \Isom(\H^n)$ is a {\em reflection group}, i.e., it is generated by reflections in the faces
of a convex acute polyhedron $\Phi\subset \H^n$ of finite volume (the fundamental domain of $\Ga$).
If $f$ is the number of facets of $\Phi$, then one can show that
$$
\dim H^1(\Ga, Ad(\rho))=f-n-1,
$$
see \cite{Kapovich(1994)}. The facets of $\Phi$ correspond to vectors spanning $H^1$.
If $n\ge 4$, it is unclear which (if any) of these infinitesimal deformations are integrable.
Of course, in {\em some} examples {\em some} of these infinitesimal
deformations are integrable, since they appear as infinitesimal bending
deformations. If $n=3$, then $X_3(\Ga)$ is smooth near $[\rho]$ and has dimension $f-4$, see  \cite{Kapovich(1994)}.

\smallskip
And that's all for $n\ge 3$.

\begin{problem}
[P. Storm]\label{storm}
 Let $M$ be a compact hyperbolic $(n+1)$-dimensional manifold with nonempty totally-geodesic boundary, $n\ge 3$.
Let $\Ga:=\pi_1(M)\subset \mob$.
Is it true that $\Ga$ is rigid in $\mob$?
\end{problem}

Note that (by Mostow rigidity) rigidity of $\Ga$ would follow if we knew that for each component
$S$ of $\partial M$, the fundamental group $\Ga_{S}:= \pi_1(S)$ is  rigid in $\mob$.
At the moment, we do not have results in either direction of this problem:

1. It is unclear if any of the  rigid hyperbolic 3-manifolds, or their disjoint union (see Example \ref{rigid}),
bounds a compact  hyperbolic 4-manifold.

2. Even if some $\Ga_{S}$ is not  rigid, it is unclear if its deformations extend
to deformations of $\Ga$.

\medskip
The only known example of a rigid group $\Ga$ (as in Problem
\ref{storm}) is the fundamental group of a 4-dimensional
hyperbolic {\em orbifold}. Moreover, in this example the group
$\Ga_S$ is
 {\em not} rigid:

\smallskip
Consider the 120-cell $D^4\subset \H^4$ which appears in \cite{Davis(1985)}.
Pick a facet $F\subset D^4$.
Let $\Ga\subset \mbt$ be the Kleinian group generated by reflections in all facets of $D^4$ except for
$F$. Then $\Ga$ is the fundamental group of a right-angled 4-dimensional reflection orbifold ${\mathcal O}$
with boundary (the boundary $S=\partial {\mathcal O}$  corresponds to the facet $F$).
The subgroup $\Ga_S:=\pi_1(S)$ is the Coxeter group generated by reflections in the facets
of the regular right-angled hyperbolic dodecahedron. In particular, $X_3(\Ga)$
is a smooth 8-dimensional manifold near $[\iota]$, where $\iota: \Ga_S\hook \mbd\hook \mbt$
is the identity embedding.

\begin{thm}
[M. Kapovich]
$\Ga$ is  rigid in $\mbt$.
\end{thm}

Assume now that $n=3$.

\smallskip
(d) The first obstruction to integrability of infinitesimal deformations
is always zero, see \cite{Kapovich-Millson(1996b)}.

\begin{question}
Suppose that $\Ga\subset \mbd$ is a cocompact lattice. It it true that the character variety
$X_3(\Ga)$ is smooth at the point $[\rho]$?
\end{question}

\smallskip
(e) Finally, there are several constructions which work for specific examples
of lattices $\Ga\subset \mbd$, e.g.  {\em stumping deformations}
\cite{Apanasov90}, generalized in \cite{Tan}.

\medskip
We recall the following

\begin{conjecture}\label{b1}
Suppose that $\Ga\subset \mob$ is a lattice. Then $\Ga$ contains a finite-index subgroup $\Ga'$
such that $\Ga'$ has infinite abelianization, i.e., $H^1(\Ga', \R)\ne 0$.
\end{conjecture}

We refer the reader to \cite{Li-Millson, Lubotzky(1996), Millson(1976), Raghunathan-Venkataramana,
 Schwermer(2004)} for various results towards this conjecture in the case of arithmetic lattices in $\Isom(\H^n)$.
The methods used in these papers for proving virtual nonvanishing of the first cohomology group usually
also apply to the cohomology groups $H^1(\Ga, Ad(\rho))$, where $\rho: \Ga\to \Isom(\H^{n+1})$ is the natural embedding.
On the other hand, the proofs of special cases of Conjecture \ref{b1} for hyperbolic 3-manifolds which use the methods of 3-dimensional
topology (see e.g. \cite{Lackenby(2006)}), usually provide no information about rigidity of $\Ga$
in $\Isom(\H^4)$.

\begin{conjecture}
Suppose that $\Ga\subset \Isom(\H^n)$ is a lattice.
Then there exists a finite-index subgroup $\Ga'\subset \Ga$ so that
$H^1(\Ga', Ad(\rho))\ne 0$.
\end{conjecture}

On the other hand, some uniform torsion-free lattices in $\mbd$ are  rigid in $\mbt$:

\begin{example}
\label{rigid} In \cite{Kapovich(1994)} we constructed examples of (torsion-free) cocompact lattices $\Ga$ in $\mbd$
for which $H^1(\Ga, Ad(\rho))= 0$, where $\rho: \Ga\to \mbt$
is the natural embedding. The quotient manifolds $\H^3/\Ga$ in these examples are
non-Haken. K.~Scannell  \cite{Scannell02} constructed analogous examples with Haken quotients  $\H^3/\Ga$.
\end{example}

More specifically, it was proved in  \cite{Kapovich(1994)} that for
every hyperbolic 2-bridge knot $K\subset \S^3$, there are infinitely many (hyperbolic) Dehn surgeries on $K$,
so that for the resulting manifolds $M_j$, $j\in \N$, we have
$$
H^1(\Ga_j, Ad(\rho))= 0, \hbox{~where~} \Ga_j=\pi_1(M_j).
$$

\subsubsection{Failure of quasiconformal isotopy}

The goal of this section is to construct examples of Kleinian groups which are quasiconformally conjugate, but
cannot be deformed to each other. As the reader will see, the tools for constructing such examples were available 12 years ago.
I realized that such examples exist only recently, while working on this survey.

\begin{thm}\label{nonisotopy}
There exists a pair of convex-cocompact Kleinian groups $\Del_1, \Del_2\subset \mbp$ and a quasiconformal homeomorphism
$f: \S^5\to \S^5$ conjugating $\Del_1$ to $\Del_2$, which is not isotopic to the identity through homeomorphisms
$h_t: \S^5\to \S^5$ such that
$$
h_t \Del_1 h_t^{-1}\subset \mbp.
$$
\end{thm}
\proof We begin with a lattice $\Ga=\pi_1(N)$, where $N=M_j$ is as in the discussion of Example \ref{rigid}
and $K\subset \S^3$ is the figure 8 knot.   Consider the representation
$$
\rho_1: \Ga\hook \mbd \hook \mbp$$
obtained by the composition of natural embeddings. Then
$$
H^1(\Ga, Ad(\rho_1))= H^1(\Ga, V_3\oplus V_3\oplus V_3\oplus \R^3),
$$
where $V_3=\R^{3,1}$ and $\R^3$ is the trivial $3$-dimensional $\R\Ga$-module. Since $H^1(\Ga, V_3)=0$ by \cite{Kapovich(1994)} and
 $H^1(\Ga, \R^3)=0$ since $N$ is a rational homology sphere, we obtain
$$
H^1(\Ga, Ad(\rho_1))= 0.
$$
Therefore $\rho_1$ is  rigid. If $N$ is an integer homology 3-sphere,
then nonvanishing of the Casson invariant of $K$ implies that $\Ga$ admits a nontrivial homomorphism
$$
\theta: \Ga\to SO(3),
$$
which lifts to $SU(2)$, see \cite{Akbulut-McCarthy(1990)}. If $M_j$ is not an integer homology sphere, then
$\Ga$ has nontrivial abelianization and hence we also obtain a nontrivial homomorphism $\theta: \Ga\to SO(3)$
with cyclic image. In any case, we twist the representation $\rho_1$ by $\theta$:
$$
\rho_2=\rho_1\times \theta: \Ga\to \mbd \times SO(3)\subset \mbp.
$$
It is clear that $[\rho_1], [\rho_2]$ are distinct points of $X_5(\Ga)$.
The images of $\rho_1$ and $\rho_2$ are $2$-fuchsian, convex-cocompact groups $\Del_1, \Del_2\subset \mbp$.
We obtain the isomorphism
$$
\rho:=\rho_2\circ \rho_1^{-1}: \Del_1\to \Del_2
$$
Clearly,
$$
M^5(\Del_1)= N\times \S^2,
$$
while $M^5(\Del_2)$ is the 2-sphere bundle over $N$ associated with $\theta$. It is easy to see that the latter bundle
is (smoothly) trivial. Therefore we obtain a diffeomorphism
$$
h: M^5(\Del_1)\to M^5(\Del_2)
$$
which lifts to a $\rho$-equivariant diffeomorphism $f: \Om(\Del_1)\to \Om(\Del_2)$. The latter extends
to a quasiconformal homeomorphism $f: \S^5\to \S^5$ by Theorem \ref{tukia}. If there was a continuous family of
homomorphisms $\rho_t$ connecting $\rho$ to the identity embedding $\Del_1\to \Del_2$, then the representation $\rho_1$
would not be  rigid in $X_5(\Ga)$. Contradiction. \qed

By embedding naturally the groups $\Del_1, \Del_2$   to $\mob$ for $n\ge 6$ one obtains higher-dimensional examples.

\subsection{Stability theorem}\label{stability}

Let $\Ga\subset \mob$ be a geometrically finite Kleinian group. Consider the set of cusps in $\Ga$:
%conjugacy classes
%of maximal parabolic subgroups
$$
 [\Pi_1],...,[\Pi_m],
$$
where $\Pi_i$ are maximal parabolic subgroups of $\Ga$. We define the {\em (topologically) relative} representation variety
$$
R_n^{top}(\Ga)=\{\rho: \Ga\to \mob: \rho(\Pi_i) \hbox{~is topologically conjugate to~} \Pi_i \hbox{~~in~~} \S^n, \forall i\}
$$
and the {\em (quasiconformally) relative} representation variety
$$
R_n^{qc}(\Ga)=\{\rho: \Ga\to \mob: \rho(\Pi_i) \hbox{~is quasiconformally conjugate to~} \Pi_i \hbox{~~in~~} \S^n, \forall i\}.
$$

Let $Homeo(\S^n)$ and $QC(\S^n)$ be the groups of homeomorphisms and quasiconformal homeomorphisms of $\S^n$ with
the topology of uniform convergence.
Let $X^{top}_n(\Ga), X^{qc}_n(\Ga)$ be the projections of $R_n^{top}(\Ga), R_n^{qc}(\Ga)$ to $X_n(\Ga)$.
Let $\iota: \Ga\to \mob$ be the identity embedding.
Then the {\em Stability Theorem} for geometrically finite groups states
that every homomorphism $\rho$ of $\Ga$ sufficiently close to
$\iota$ is induced by a (quasiconformal) homeomorphism $h_{\rho}$ close to the identity and depending continuously on $\rho$.
More precisely:

\begin{thm}\label{stabilitythm}
(Stability theorem, see \cite{Epstein-Canary-Green, Goldman(1987b), Kapovich(1990), Marden, Sullivan(1985)}.)
There exist neighborhoods $U^{top}, U^{qc}$ of $\iota$ in $R^{top}_n(\Ga), R^{qc}_n(\Ga)$ respectively, and continuous maps
$$
L^{top}: U^{top}\to Homeo(\S^n),\quad L^{qc}: U^{qc}\to QC(\S^n)
$$
so that
$$
L^{top}(\iota)= L^{qc}(\iota)=Id,
$$
and for every $\rho\in U^{top}$, resp. $\rho\in U^{qc}$, the homeomorphism $L^{top}(\rho)$, resp. $L^{qc}(\rho)$
is $\rho$-equivariant.
\end{thm}

This theorem was first proved by A.~Marden in \cite{Marden} in the case $n=3$. Marden was using
convex finitely-sided fundamental domains with simplicial links of vertices:
Such polyhedra are generic among the Dirichlet fundamental domains, see \cite{Jorgensen-Marden(1986)}.
Marden then argued that a small perturbation of such fundamental domain is again a fundamental domain (by the Poincar\'e
fundamental polyhedron theorem).
Moreover, the {\em simplicial} assumption implies that the combinatorics of the fundamental domain does not change
under a small perturbation.
This allowed Marden to construct an equivariant quasiconformal homeomorphism close to the identity.
This argument does not readily generalize to higher dimensions, mainly because finiteness of the number of
faces is not equivalent to geometric finiteness. (Otherwise, the same argument goes through.)

D.~Sullivan \cite{Sullivan(1985)} considered the case of general $n$, but assumed that
$\Ga$ is convex-cocompact. Then he proved the existence of a homeomorphism $h_{\rho}$
defined on the limit set of $\Ga$ and the fact that it depends continuously on $\rho$.
The fact that
$$
h_\rho: \La(\Ga)\to \La(\rho(\Ga))
$$
is necessarily  quasi-symmetric,  then follows from Tukia's theorem \ref{tukia}.
One then has to show existence of a $\rho$-equivariant diffeomorphism of the domains of discontinuity
$$
f_{\rho}: \Om(\Ga)\to \Om(\rho(\Ga))
$$
smoothly depending on $\rho$. This is achieved by appealing
to Thurston's holonomy theorem (see \cite{Epstein-Canary-Green, Goldman(1987b)})
for the M\"obius structures on the manifold $M^n(\Ga)$, as it is done in  \cite{Izeki(2000), Kapovich(1990)}.
The homeomorphisms $h_\rho$ and $f_\rho$ yield a $\rho$-equivariant quasiconformal homeomorphism of the $n$-sphere by
Theorem \ref{tukia}.

The proof in \cite{Epstein-Canary-Green} is a good alternative to the above argument; it is also
sufficiently flexible to handle  the case of geometrically
finite Kleinian groups with cusps. Namely, instead of  working with the $n$-dimensional manifold $M^n(\Ga)$
one works with the convex hyperbolic $(n+1)$-manifold
$$
H(\Ga):= Hull_{\eps}(\La(\Ga))/\Ga.
$$
An analogue of Thurston's holonomy theorem for manifolds with boundary applies in this case. Thus, for $\rho\in U^{top}$,
there exists a hyperbolic structure $s(\rho)$ (with the holonomy $\rho$) on the thick part
$$
H(\Ga)_{[\mu, \infty)}
$$
of the manifold $H(\Ga)$. Moreover, convexity of the boundary for the new hyperbolic structures (away from the cusps)
persists under small perturbations of the hyperbolic structure. Therefore, if $\Ga$ is convex-cocompact, $\Ga':=\rho(\Ga)$
is again convex-cocompact and $\rho: \Ga\to \Ga'$ is an isomorphism. If $\Ga$ is merely geometrically finite,
because $\rho$ belongs to the relative representation variety, it follows that the hyperbolic structure $s(\rho)$
extends to a convex complete hyperbolic structure on the cusps. This argument also yields
a $\rho$-equivariant diffeomorphism
$$
Hull_{\eps}(\La(\Ga))  \to Hull_{\eps}(\La(\Ga'))
$$
depending continuously on $\rho$. To get from the convex hulls to the domain of discontinuity one uses the
existence of the canonical equivariant diffeomorphisms (``the nearest-point projections'')
$$
\Om(\Ga)\to \partial Hull_{\eps}(\La(\Ga)),  \quad \Om(\Ga')\to \partial Hull_{\eps}(\La(\Ga')).
$$
We refer the reader to \cite{Epstein-Canary-Green} for the details.

\medskip
Sullivan also had a converse to the Stability Theorem for (finitely-generated) subgroups on $\mbd$:

\begin{thm}
\label{converse}
(D. Sullivan, \cite[Theorem A$^\prime$]{Sullivan(1985)})
If  a (finitely-generated) Kleinian subgroup of $\mbd$ is stable in the sense of Theorem \ref{stabilitythm}, then it is geometrically finite
or its identity embedding in $\mbd$ is  rigid in $X_2(\Ga)$.
\end{thm}

It was proved in \cite{Kapovich00} that every  rigid $\Ga$ in the above theorem has to be geometrically finite.
Now it, of course, follows from the positive solution of the Bers--Thurston density conjecture (geometrically finite groups
are dense among Kleinian subgroups of $\mbd$).

\begin{question}
Does Theorem \ref{converse} hold for subgroups of $\mob$, $n\ge 3$?
\end{question}

We expect the answer to be negative.

\subsection{Space of discrete and faithful representations}
\label{global}

Let ${\mathcal D}_n(\Ga)\subset X_n(\Ga)$ denote the subset corresponding to discrete,
injective and nonelementary representations of $\Ga$.%, which map parabolic elements to parabolic elements.

\begin{thm}
[Chuckrow--J{\o}rgensen--Wielenberg] ${\mathcal D}_n(\Ga)\subset X_n(\Ga)$ is closed. See for instance \cite{Wielenberg, Martin(1989a)}.
\end{thm}

It turns out that there exists another way to construct limits of sequences of Kleinian groups,
by regarding them as {\em closed subsets} of $\mob$. This leads to the topology of {\em geometric convergence}
of Kleinian groups. With few exceptions, the space of Kleinian groups is again closed in this topology (see e.g.
\cite{Thurston(1997),  Kapovich00}). In general, ${\mathcal D}_n(\Ga)$ is not compact. Nevertheless,
this space can be compactified by {\em projective classes} of nontrivial $\Ga$-actions on real trees. This compactification
generalizes Thurston's compactification of the Teichm\"uller space. The compactification by actions on trees was first
defined by J.~Morgan and P.~Shalen \cite{Morgan-Shalen(1988)} and J.~Morgan \cite{Morgan(1986)} using algebraic geometry.
More flexible, geometric, definitions of this compactification were introduced by M.~Bestvina
\cite{Bestvina(1988)} and F.~Paulin \cite{Paulin(1988)}. See also \cite{Kapovich00} for
the construction of this compactification using ultralimits of metric spaces.

This viewpoint provides a powerful tool for proving compactness of ${\mathcal D}_n(\Ga)$ for certain classes
of groups: If ${\mathcal D}_n(\Ga)$ is non-compact then $\Ga$ admits nontrivial action on a certain $\R$-tree. One
then proves that such action cannot exist. The tools for proving such non-existence theorems are originally due to
Morgan and Shalen (but limited to the fundamental groups of 3-manifolds, see \cite{Morgan-Shalen(1988)}); a
much more general method is due to E.~Rips (Rips theory), see \cite{Bestvina-Feighn(1995)}.  One then obtains
the following (see e.g. \cite{Kapovich00}):

\begin{thm}
[Rips--Thurston Compactness theorem]\label{morg}
Suppose that $\Ga$ is a finitely-presented group which
does not split as an amalgam over a virtually abelian group.
Then ${\mathcal D}_n(\Ga)$ is compact.
\end{thm}

\begin{rem}
W.~Thurston \cite{Thurston(1986)} proved this theorem for a certain class of 3-manifold groups in the case $n=2$.
\end{rem}

Unfortunately, none of the known proofs of Theorem \ref{morg} gives an explicit bound on the ``size''
of  ${\mathcal D}_n(\Ga)$.

\begin{problem}
\label{constructive}
Find a ``constructive'' proof of Theorem \ref{morg}. More precisely, consider a group
$\Ga$ with a finite presentation $\<g_1,..,g_k|R_1,..,R_m\>$.
Given $[\rho]\in {\mathcal D}_{n-1}(\Ga)$
define
$$
B_{n-1}([\rho]):= \inf_{x\in \H^{n}} \max_{i=1,...,k} d(x, \rho(g_i)(x)).
$$
Find an explicit constant $C$, which depends on $n$, $k, m$
and the lengths of the words $R_i$, so that the function $B_{n-1}: {\mathcal D}_{n-1}(\Ga)\to \R$ is bounded from
above by $C$.
\end{problem}

\noindent Theorem \ref{morg} suggests that one should also look for {\em geometric bounds} on
$[\rho]\in {\mathcal D}_n(\Ga)$: Even if ${\mathcal D}_n(\Ga)$ is noncompact (or its ``size'' is unknown),
one can still try to find some natural functionals on  ${\mathcal D}_n(\Ga)$  and obtain explicit bounds
(from below and from above) on these functionals.

\begin{defn}
[Diameter of a representation] Given a discrete embedding $\rho: \Ga\to \Ga'=\rho(\Ga)\subset \mob$,
consider the set $S$ of connected subgraphs $\si\subset \H^{n+1}/\Ga'$ with the property:
The map $\pi_1(\si)\to \pi_1(M)$ is surjective.

\smallskip
Then the {\em diameter} of $\rho$ is
$$
\diam(\rho):=\inf \{ \length(\si): \si\in S\}.
$$
\end{defn}

\begin{problem}
Given a group $\Ga$ as in Theorem \ref{morg}, find explicit bounds on $\diam(\rho)$  (in terms of the presentation of $\Ga$)
for representations $[\rho]\in  {\mathcal D}_n(\Ga)$.
\end{problem}

Note that the positive {\em lower bound} on $\diam(\rho)$ is an easy corollary of the Kazhdan-Margulis lemma.

\begin{defn}
(Volumes of a representation) Fix a homology class $[\zeta]\in Z_p(\Ga)$, $2\le p\le cd (\Ga)$. For a representation $\rho\in
{\mathcal D}_n(\Ga)$ consider the quotient manifold $M=\H^n/\rho(\Ga)$. Define the set $E(\zeta)$ of singular
$p$-cycles $\zeta'\in Z_p(M)$ which represent the homology class $[\zeta]$ under the isomorphism
$$
H_p(\Ga)\to H_p(M)
$$
induced by the isomorphism $\rho: \Ga\to \pi_1(M)$. Lastly, define
the {\em $\rho$-volume} of the class $[\zeta]$ by
$$
Vol_\rho(\zeta):=\inf \{ Vol(\zeta'): \zeta'\in E(\zeta)\}.
$$
\end{defn}

Let $||\zeta||$ denote the Gromov-norm of the class $[\zeta]$ and let $c_p$ denote the volume
of the regular ideal geodesic $p$-simplex in $\H^p$. Then an easy application of Thurston's ``chain-straightening''
is the inequality
$$
Vol_\rho(\zeta)\le c_p||\zeta||
$$
for all $\rho$, $[\zeta]$ and $p\ge 2$. However good lower bounds on the volume are considerably more difficult to get.

Given a hyperbolic manifold $M$ define $H^{par}_p(M)$ to be the image in $H_p(M)$ of the $p$-th homology group
of the union of all cusps of $M$. Then for every {\em parabolic} class $[\zeta]\in H^{par}_p(M)$ and every $\rho$, we clearly
have
$$
Vol_\rho(\zeta)=0.
$$
However there exists a positive constant $\eps=\eps(p,n)$ such that for every $p>0$, every non-cuspidal
class $[\zeta]$ and every $\rho$, we obtain
$$
Vol_\rho(\zeta)\ge \eps,
$$
see \cite{Kapovich07}. Below are some more interesting lower bounds on the volume:

\begin{thm}
(Follows directly from \cite[Theorem 5.38]{Gromov00}\footnote{I am grateful to A. Nabutovsky for this reference.}).
Let $\Ga$ be isomorphic to the fundamental group
of a compact aspherical $k$-manifold $N$ and $[\zeta]=[N]$ be the fundamental class of $M$. Then
there exists a universal (explicit) constant $c(p,n)>0$ depending only on $p$ and $n$, such that
$$
Vol_\rho(\zeta)\ge c(p,n)||N||.$$
\end{thm}

One gets better estimates using the work of Besson, Courtois and Gallot \cite{Besson-Courtois-Gallot(1999)}
\footnote{I am grateful to J.\! Souto for pointing this out.}:

\begin{thm}
\label{contract}
Fix a representation $[\phi]\in {\mathcal D}_n(\Ga)$. Then for every $[\rho]\in {\mathcal D}_n(\Ga)$ and $p\ge 3$ we obtain
$$
\left( \frac{p}{n} \right) ^{p+1}  Vol_{\phi}(\zeta)\le Vol_{\rho}(\zeta).
$$
\end{thm}

For instance, if $\Ga':=\phi(\Ga)$ happens to be a uniform lattice in $\Isom(\H^p)$, we obtain

\begin{corollary}
\label{BCGcor}
For every $[\rho]\in {\mathcal D}_n(\Ga)$ and $p\ge 3$ we have
$$
\left( \frac{p}{n} \right) ^{p+1}  Vol(M') \le Vol_{\rho}(\zeta),
$$
where $M'=\H^p/\Ga'$ and $[\zeta]$ is the fundamental class.
\end{corollary}

\subsection{Why is it so difficult to construct higher-dimensional geometrically infinite Kleinian groups?}

$~$

(i). The oldest trick for proving existence of geometrically infinite groups is due to L.~Bers \cite{Bers(1970)}:

Start with  (say) a convex-cocompact subgroup $\Ga\subset \mob$. Let $Q(\Ga)\subset {\mathcal D}_n(\Ga)$ be the
(open) subset of representations induced by quasiconformal conjugation. Let $Q_0(\Ga)$ denote the component of
$Q(\Ga)$ containing the (conjugacy class of) identity representation $[\rho_0]$.
We assume that the closure of $Q_0(\Ga)$ is not
open in $X_n(\Ga)$. Then there exists a curve
$[\rho_t]\in X_n(\Ga), t\in [0, 1]$, so that $\rho_1$ is either nondiscrete or non-injective.
Since ${\mathcal D}_n(\Ga)$ is closed, it follows that there exists $s\in (0,1)$ such that
$[\rho_s(\Ga)]$ belongs to ${\mathcal D}_n(\Ga)$ but $\rho_s(\Ga)$  is not convex-cocompact.
If $\Ga'=\rho_s(\Ga)$ contains no parabolic elements,
it would follow that $\Ga'$ is isomorphic to $\Ga$ and is not geometrically finite. However, it could happen that
the frontier of $Q_0(\Ga)$ consists entirely of the classes $[\rho]$ for which $\rho(\Ga)$ contains parabolic elements.

The latter cannot occur if $n=2$ for dimension reasons: The set of parabolic elements of $PSL(2, \C)$ has real codimension 2
and, hence, does not separate. However for all $n\ne 2$,
the set of parabolic elements has real codimension 1 and this argument is inconclusive.

\medskip
One can try to apply the above argument in the case of a codimension 1 fuchsian group $\Ga\subset \mob$
which acts as a cocompact lattice on $\H^n\subset \H^{n+1}$. Suppose that $M=\H^n/\Ga$ contains a totally-geodesic
compact hypersurface $S$. Then we have the circle $\S^1$ worth of bending deformations $\rho_t$ along $S$.
As $t=\pi$, the image of $\rho_t$ is again contained in $\Mob(\S^{n-1})$. Therefore $\rho_\pi$ is either nondiscrete or non-injective.
However, conceivably, in all such cases, for $[\rho_s]\in \partial Q_0(\Ga)$
the representation $\rho_s$ is geometrically finite (because its image may contain parabolic elements).
It happens, for instance, if $\Ga$ is a reflection group.

Note that even when $n=2$ and we are bending a 1-fuchsian group $\Ga$, it is hard to predict
which simple closed geodesics $\al\subset \H^2/\Ga$ yield geometrically infinite groups (via bending along $\al$).

\medskip
(ii). One can try to construct explicit examples of fundamental
domains, following, say, T.~J{\o}rgensen \cite{Jorgensen(1977)} or
A.~Marden and T.~J{\o}rgensen \cite{Jorgensen-Marden(1979)}.

The trouble is that constructing fundamental polyhedra with infinitely many faces in $\H^4$ is quite a bit harder than in $\H^3$.
One can try to find a lattice $\widehat\Ga\subset \mbt$ which contains a nontrivial finitely-generated normal subgroup $\Ga$
of infinite index.\footnote{If $\widehat\Ga$ is a Kleinian group containing a nontrivial normal subgroup $\Ga$ of infinite index, then
$\Ga$ is necessarily geometrically infinite.}
This is, probably, the most promising approach, since it works for complex-hyperbolic lattices in $PU(2,1)$,
cf. \cite{Kapovich(1998a)}. One can try to imitate Livne's examples, by constructing $\Ga\subset \widehat\Ga$ such that
$\widehat\Ga/\Ga$ is isomorphic to a surface group. This would require coming up with a specific compact
convex polyhedron in $\H^4$ such that the associated 4-manifold appears as a (singular) fibration over a surface.

\medskip
(iii). One can try to use combinatorial group theory. Note that there are plenty of examples of (mostly 2-dimensional)
Gromov-hyperbolic groups $\widehat\Ga$ which contain nontrivial finitely-generated normal subgroups $\Ga$ of infinite index.
See e.g. \cite{Bestvina-Feighn(1992), Brady(1999),  Mosher, Rips(1982)} for the examples which are not 3-manifold groups.
However embedding a given hyperbolic group $\widehat\Ga$ in $\mob$ is a nontrivial task, cf.
\cite{Bourdon97, Kapovich05} and discussion in Section \ref{constrains}.
The groups considered in \cite{Kapovich05}, probably provide the best opportunity here, since most of
them do not pass the {\em perimeter test} of J.~McCammond and D.~Wise \cite{McCammond-Wise}.
(If a geometrically finite group $\widehat\Ga$ satisfies the perimeter test,
then every finitely-generated subgroup of $\widehat\Ga$ is geometrically finite.)

\medskip
(iv). {\em What would geometrically infinite examples look like?} Let $\Ga\subset \mbd$ be a
singly-degenerate group; assume for simplicity that the injectivity radius of $\H^3/\Ga$ is bounded away from zero.
Let $S$ denote the boundary of
$$
Hull(\La(\Ga))/\Ga
$$
and $\la\subset S$ be the ending lamination of $\Ga$. Then every leaf of $\la$ lifts to an exponentially distorted curve $\kappa$
in $\H^3$: Given points $x, y\in  \kappa$, their extrinsic distance $d(x, y)$ in $\H^3$ is roughly the
logarithm of their intrinsic distance along $\kappa$.

One would like to imitate this behavior in dimension 4.
Let $M$ be a closed hyperbolic 3-manifold containing an embedded  compact
totally-geodesic surface $S\subset M$. Let $\la\subset S$ be an ending lamination from the above example.
One would like to construct a complete hyperbolic 4-manifold $N$ homotopy-equivalent to $M$, so that under the (smooth)
homotopy-equivalence $f: M\to N$ we have:

For every leaf $L$ of $\la$, $f(L)$ lifts to an exponentially distorted curve in $\H^4$.

Then $\pi_1(N)$ will necessarily be a geometrically infinite subgroup $\Ga$ of $\mbc$.
At the moment it is not even clear how to make this work with a hyperbolic metric on $N$ replaced by a complete
Riemannian metric of negatively pinched sectional curvature, although constructing a Gromov-hyperbolic metric with this behavior
is not that difficult. (Recall that a Riemannian metric is said to be {\em negatively pinched} if its sectional curvature
varies between two negative numbers.)
An example $\Ga$ of this type is likely to have two components of $\Om(\Ga)$: One contractible and one not.

More ambitiously, one can try to get a singly degenerate group $\Ga\subset \mbt$ (so that $\Om(\Ga)$ is contractible and
$M^3(\Ga)$ is compact). How would such an example  look like? One can imagine taking a 1-dimensional quasi-geodesic
foliation $\la$ of the 3-manifold $M$ as above and then requiring that for every leaf $L\subset \la$,
the curve $f(L)$ lifts to an exponentially distorted curve in $\H^4$. At the moment I do not see even a
Gromov-hyperbolic model of this behavior. Another option would be to work with 2-dimensional laminations $\nu$
(with simply-connected leaves) in $M$ and require every leaf $L\subset \nu$ to correspond to an exponentially
distorted surface in $\H^4$ (or a negatively-curved simply-connected 4-manifold),
which limits to a single point in $\S^3$.

\begin{problem}
Construct a complete negatively pinched 4-dimensional Riemannian manifold $N$ homotopy-equivalent to a hyperbolic
3-manifold $M$, so that the convex core of $N$ either has exactly one boundary component or equals $N$ itself.
\end{problem}

\begin{question}
Is there a geometrically infinite Kleinian subgroup of $\mob$ whose limit is homeomorphic to the  Menger curve?
Is there a geometrically infinite Kleinian subgroup of $\mob$ which is isomorphic to the fundamental
group of a closed aspherical manifold of dimension $\ge 3$? Are there examples of such groups acting
isometrically on complete negatively pinched manifolds? Are there examples of hyperbolic (or even negatively curved)
4-manifolds $M$ such that $\pi_1(M)=\Ga$ fits into a short exact sequence
$$
1\to \pi_1(S)\to \Ga\to \pi_1(F)\to 1,
$$
where $S, F$ are closed hyperbolic surfaces? Note that there are no complex-hyperbolic examples of this type, see
\cite{Kapovich(1998a)}.
\end{question}

\section{Algebraic and topological constraints on Kleinian groups}
\label{constrains}

Sadly, there are only few known algebraic and topological restrictions on Kleinian subgroups in $\mob$
that do not follow from the {\em elementary} restrictions, which come from the restrictions
on geometry of complete negatively curved Riemannian manifolds. Examples of the elementary restrictions on a
Kleinian group $\Ga$ are:

\smallskip
1. Every solvable subgroup of a Kleinian group is virtually
abelian.

\smallskip
2. The normalizer (in $\Ga$) of an infinite cyclic subgroup of
$\Ga$ is virtually abelian.

\smallskip
3. Every elementary (i.e., virtually abelian) subgroup $\Del\subset \Ga$ is contained in a unique maximal
elementary  subgroup $\tilde{\Del}\subset \Ga$.

\smallskip
4. Every Kleinian group has finite (virtual) cohomological dimension.

\smallskip
\noindent In this section we review known {\em nonelementary} algebraic and topological restrictions on Kleinian groups.

\subsection{Algebraic constraints}

\begin{defn}
An {\em abstract Kleinian group} is a group $\Ga$ which admits a discrete embedding in $\mob$ for some $n$.
Such a group is called {\em elementary} if it is virtually abelian.
\end{defn}

In order to eliminate trivial restrictions on abstract Kleinian groups one can restrict attention
to  Gromov-hyperbolic Kleinian groups. Below is the list of known algebraic constraints on Kleinian groups
under this extra assumption:

\medskip
1. Kleinian groups are residually finite and virtually torsion-free.\footnote{It is widely believed that there are
Gromov-hyperbolic groups which are not residually finite.}
(This, of course, holds for all finitely generated matrix groups.)

\medskip
2. Kleinian groups satisfy the {\em Haagerup property}, in particular, infinite Kleinian groups
do not satisfy property (T), see \cite{haagerup}.

\medskip
3. If a Kleinian group $\Ga$ is {\em K\"ahler}, then $\Ga$ is virtually isomorphic to the fundamental group
of a compact Riemann surface. This is a deep theorem of J.~Carlson and D.~Toledo \cite{Carlson-Toledo}, who proved that every homomorphism of
a    K\"ahler group to $\mob$ either factors through a virtually surface group, or its image fixes a point in $\B^{n+1}$.

\medskip
Recall that a topological group $G$ is said to satisfy the {\em Haagerup property} if it admits a (metrically)
proper continuous isometric action on a Hilbert space $H$. An action of a metrizable
topological group $G$ on $H$ is {\em metrically proper}  if for every bounded subset $B\subset H$,
the set
$$
\{g\in G: g(B)\cap B\ne \emptyset\}$$
is a bounded subset of $G$.  Since $\mob$ satisfies the Haagerup property for every $n$ (see e.g.
\cite{haagerup}), all Kleinian groups also do.

A group $\pi$ is called {\em K\"ahler} if it is isomorphic to the fundamental group of a compact K\"ahler manifold.
For instance, every uniform lattice in $\C\H^n$ is  K\"ahler; therefore it cannot be an abstract Kleinian group
unless $n=1$.

\begin{rem}\label{haag}
A  (finitely-generated) group satisfies the Haagerup property if and only if it admits an isometric (metrically) properly
discontinuous action on the  infinite dimensional hyperbolic space $\H^{\infty}$, see  \cite[7.A.III]{Gromov(1993)}.
The result of Carlson and Toledo shows that (for Gromov-hyperbolic groups) there are nontrivial obstructions to
replacing these infinite-dimensional actions with finite-dimensional ones.
\end{rem}

\begin{observation}
All currently known nontrivial restrictions on abstract Kleinian groups can be traced to 1, 2 or 3.
\end{observation}

\begin{problem}
Find other restrictions on abstract Kleinian groups.
\end{problem}

Potentially, some new restrictions would follow from the Rips-Thurston compactness theorem. The
difficulty comes from the following. Let $\Ga$ be a Gromov-hyperbolic group which admits no nontrivial
isometric actions on $\R$-trees. Then (see \cite{Kapovich00}) there exists $C<\infty$, such that
for every sequence $[\rho_j]\in {\mathcal D}_n(\Ga)$, we obtain a uniform bound
\begin{equation}\label{inequality2}
B_n([\rho_j])\le C,
\end{equation}
where $B_n:  {\mathcal D}_n(\Ga)\to \R$ is the minimax function defined in Problem \ref{constructive}. If $n$
were fixed, then the sequence $(\rho_j)$ would subconverge to a representation to $\mob$ (for some choice
of representations $\rho_j$ in the classes $[\rho_j]$). However, since we
are not fixing the dimension of the hyperbolic space on which our $\Ga$ is supposed to act,  the
inequality (\ref{inequality2}) does not seem to yield any useful information. By taking an ultralimit of $\rho_j$'s
we will get an action of $\Ga$ on an infinite-dimensional hyperbolic space.
This action, however, may have a fixed point,  since
$$
\lim_{n\to \infty} \mu_n=0,
$$
where $\mu_n$ is the Margulis constant for $\H^n$. See also Remark \ref{haag}.

\begin{example}
Let $M^3$ be a closed non-Haken hyperbolic 3-manifold, so that $\Pi:=\pi_1(M)$ contains a maximal
1-fuchsian subgroup $F$. For each automorphism $\phi: F\to F$ we define the HNN extension
$$
\Ga_\phi:= \Pi*_{F\cong_\phi F}= \< \Pi, t| t g t^{-1}= \phi(g), \forall g\in F\>.
$$
Then $\Ga_\phi$ is Gromov-hyperbolic for all pseudo-Anosov automorphisms $\phi$, see \cite{Bestvina-Feighn(1995)}.
It is a direct corollary of Theorem \ref{morg} that for every $n$, only finitely many of the groups $\Ga_\phi$ embed in
$\mob$ as Kleinian subgroups. Is it true that there exists a pseudo-Anosov automorphism $\phi$ such that
$\Ga_\phi$ is not an abstract Kleinian group?
\end{example}

Infinite finitely-generated Gromov-hyperbolic Coxeter groups are all linear, satisfy the Haagerup property and are not
K\"ahler (except for the virtually surface groups).

\begin{problem}
Is it true that every finitely-generated Gromov-hyperbolic Coxeter group $\Ga$ is an abstract Kleinian group?
\end{problem}

Note that there are Gromov-hyperbolic Coxeter groups $\Ga$ which do not admit discrete embeddings
$\rho:\Ga\to \mob$ (for any $n$), so that the Coxeter generators of $\Ga$ act as reflections in the faces of a
fundamental domain of $\rho(\Ga)$, see \cite{Felikson-Tumarkin}.

The answer to the next question is probably negative, but the examples would be tricky to construct:

\begin{question}
Is it true that a group weakly commensurable to a Kleinian group is also a Kleinian group? Even more ambitiously:
Is the property of being Kleinian a quasi-isometry invariant of a group?
\end{question}

Recall that two groups $\Ga$ and $\Ga'$ are called {\em weakly commensurable}  if there exists a chain of groups and homomorphisms
$$
\Ga=\Ga_0 \to \Ga_1 \leftarrow \Ga_2 \to \Ga_3 .... \leftarrow \Ga_{k-1} \to \Ga_k=\Ga',
$$
where each arrow $\Ga_i\to \Ga_{i\pm 1}$ is a homomorphism whose kernel and cokernel are finite.

\smallskip
There are few more known algebraic restrictions on geometrically finite Kleinian groups.
All such groups are relatively hyperbolic.

We recall that a group $\Ga$ is called {\em cohopfian}
if every injective endomorphism $\Ga\to \Ga$ is also surjective.

\begin{rem}
A group $\Ga$ is called {\em hopfian} if  every epimorphism
$\Ga\to \Ga$ is injective. Every residually finite group $\Ga$ is
hopfian, see \cite{Malcev(1940)}. In particular, every Kleinian
group is hopfian.
\end{rem}

For instance, free groups and free abelian groups
are not cohopfian. More generally, if $\Ga$ splits as a nontrivial free product,
$$
\Ga\cong \Ga_1 * \Ga_2,
$$
then $\Ga$ is not cohopfian: Indeed, for nontrivial elements
$\ga_1\in \Ga_1, \ga_2\in \Ga_2$, set $\al:=\ga_1\ga_2$, and
$$
\Ga_1':= \al \Ga_1 \al^{-1}.$$
Then
$$
\Ga\cong \Ga_1' * \Ga_2
$$
is a proper subgroup of $\Ga$. On the other hand, lattices in $\Isom(\H^n)$, $n\ge 3$,
and uniform lattices in $\Isom(\H^2)$ are cohopfian. Indeed, Mostow rigidity theorem implies that if
$M_1, M_2$ are hyperbolic $n$-manifolds of finite volume (and $n\ge 3$) or compact hyperbolic surfaces, and
$M_1\to M_2$ is a $d$-fold covering, then
$$
Vol(M_1)= d Vol(M_2).
$$
On the other hand, if $M$ is a hyperbolic manifold of finite
volume (or a compact hyperbolic surfaces), then every proper embedding
$$
\pi_1(M)\to \pi_1(M)
$$
induces a $d$-fold covering $M\to M$, with $d\in \{2, 3,...,
\infty\}$. Hence $\pi_1(M)$ is cohopfian.

If a Kleinian group $\Ga\subset \mob$ fails to be cohopfian, we can iterate a proper embedding $\phi: \Ga\to \Ga$,
thereby obtaining a sequence of discrete and faithful representations
$$
\rho_i =\underbrace{\phi\circ ... \circ \phi}_{i \hbox{~~times}}
$$
of $\Ga$ into $\mob$. By analyzing such sequences, T.~Delzant and L.~Potyagailo \cite{Delzant-Potyagailo}
obtained a characterization of geometrically
finite Kleinian groups which are {\em cohopfian}. We will need two definitions in order to describe their result.

\begin{defn}
If $\Ga$ is a Kleinian group and $\Del\subset \Ga$ is an elementary subgroup, let $\t\Del$ denote
the maximal elementary subgroup of $\Ga$ containing $\Del$.
\end{defn}

\begin{defn}
Suppose a group $\Ga$ splits as a graph of groups
\begin{equation}
\Ga\cong \pi_1({\mathcal G}, \Ga_e, \Ga_v),
\end{equation}
and suppose that edge groups $\Ga_e$ of this graph are  elementary.
We say that the edge group $\Ga_e$ is {\em essentially
non-maximal} if the subgroup $\tilde\Ga_e\subset \Ga$, is not
conjugate into any of the vertex subgroups of the graph of groups ${\mathcal G}$. The splitting  is
{\em essentially non-maximal} if there
exists at least one such an edge group.
Otherwise  we say that the splitting is essentially maximal.
\end{defn}

For instance, if every edge subgroup is a maximal elementary
subgroup of $\Ga$, then the splitting is essentially maximal.

\begin{thm}
(T.~Delzant and L.~Potyagailo \cite{Delzant-Potyagailo}.)  Let $\Ga$ be a non-elementary, geometrically finite,
one-ended  Kleinian group without $2$-torsion. Then $\Ga$ is cohopfian if and
only if the following two conditions are satisfied:

\begin{itemize}
\item[\sf 1)]  $\Ga$ has no essentially non-maximal splittings.

\item [\sf 2)]  $\Ga$ does not split as an amalgamated free product
$$
\Ga= \Ga_1*_{\Ga_3} \tilde{\Ga}_3,
$$
with $\tilde{\Ga}_3$  maximal elementary,  such that  the normal closure of the
subgroup  $\Ga_3$ in $\tilde \Ga_3$ is of infinite index in $\tilde{\Ga}_3$.
\end{itemize}
\end{thm}

One of the ingredients in the proof of this theorem was the fact that nonelementary geometrically
finite groups $\Ga$ do not contain subgroups $\Ga'$, which are conjugate to $\Ga$ in $\mob$,
see \cite{Wang-Zhou}.

\begin{question}
[L. Potyagailo] Let $\Ga\subset \mob$ be a finitely generated non-elementary Kleinian
group. Suppose that $\al\in \mob$ is such that
$$
\Ga'=\al \Ga \al^{-1}\subset \Ga.
$$
Does it follow that $\Ga'=\Ga$?
\end{question}

The affirmative answer to this question for $n=2$ was given in a paper of L.~Potyagailo and
K.-I.~Ohshika \cite{Ohshika-Potyagailo} (modulo Tameness Conjecture,  Theorem \ref{tamenessconjecture}).

\begin{question}
Is the isomorphism problem solvable within the class of all finitely-presented Kleinian groups?
Note that the work of F.~Dahmani and
D.~Groves \cite{Dahmani-Groves} implies solvability of the isomorphism problem in the category
of geometrically finite Kleinian groups.
\end{question}

It was proved by M.~Bonk and O.~Schramm \cite{Bonk-Schramm} that every Gromov-hyperbolic group $\Ga$ embeds quasi-isometrically
in the usual hyperbolic space $\H^n$ for some $n=n(\Ga)$. A natural question is if one can prove an {\em equivariant}
version of this result. Note that there are many Gromov-hyperbolic
groups which are not Kleinian, e.g. groups with property (T) and Gromov-hyperbolic K\"ahler groups.
Therefore one has to relax the {\em isometric} assumption. The natural category for this is the
{\em uniformly quasiconformal} actions. Such action is a monomorphism
$$
\rho:\Ga\hook QC(\S^n)
$$
whose image consists of $K$-quasiconformal homeomorphisms with $K$ depending only on $\rho$.

\begin{problem}
Let $\Ga$ be a Gromov-hyperbolic group. Does $\Ga$ admit a uniformly quasiconformal
discrete action on $\S^n$ for some $n$? For instance, is there such an action if $\Ga$ is a uniform lattice in
$PU(n,1)$ or satisfies the property (T)?
\end{problem}

T.~Farrell and J.~Lafont \cite{Farrell-Lafont} proved that the
topological counterpart of this problem has positive solution. A
corollary of their work is that every Gromov-hyperbolic group
$\Ga$ admits a {\em convergence} action $\rho$ on the closed
$n$-ball, so that the limit set of $\Ga'=\rho(\Ga)$ is
equivariantly homeomorphic to the ideal boundary of $\Ga$ and
$\Om(\Ga')/\Ga'$ is compact and connected. We refer the reader to
\cite{Gehring-Martin85} for the definition of a convergence
action.

\subsection{Topological constraints}

The basic problem here is to find topological restrictions on the
hyperbolic manifold $\H^{n+1}/\Ga$ and on the conformal\-ly-flat
manifold $\Om(\Ga)/\Ga$, which do not follow from the algebraic
restrictions on the group $\Ga$ and from the general algebraic
topology restrictions (e.g. vanishing of the characteristic
classes). There are only few nontrivial results in this direction.
For $n=3$ we have:

\begin{thm}
(M. Kapovich \cite{Kapovich(1992a)}.) There exists a function $c(\chi)$ with the following property.
Let $S$ be a closed hyperbolic surface.  Suppose that $M^4$ is a complete hyperbolic 4-manifold  which
is homeomorphic to the total space of an $\R^2$--bundle $\xi: E\to S$ with the Euler number $e(\xi)$.
Then
$$
|e(\xi)|\le c(\chi(S)).
$$
\end{thm}

More generally,

\begin{thm}
\label{inter}
(M. Kapovich \cite{Kapovich(1992a)}.) There exists a function $C(\chi_1, \chi_2)$ with the following property.
Suppose that $M^4$ is a complete oriented hyperbolic $4$-manifold.
Let $\sigma _j: \Sigma _j \to M^4$ ($j = 1, 2$) be
$\pi_1$--injective maps of closed oriented surfaces $\Sigma_j$. Then
$$
|\< \sigma _1 , \sigma _2 \>| \le
 C(\chi (\Sigma _1 ), \chi (\Sigma _2 )).
$$
\end{thm}

Here $\<  ,  \>$ is the intersection pairing on $H_2(M^4)$. The bounds appearing in these theorems are explicit but
astronomically high. The expected bounds are linear in $\chi(S)$ and $\chi(S_i)$, $i=1, 2$, cf. \cite{Gromov-Lawson-Thurston(1988)}.

\medskip
Other known restrictions are applications of the compactness theorem \ref{morg} and therefore explicit bounds
in the following theorems are unknown.

\begin{thm}
(M. Kapovich \cite{Kapovich(1993b)}.)
Given a closed hyperbolic $n$-manifold $B$ ($n\ge 3$) there exists a number $c(B)$ so that the following is true.
Suppose that  $M^{2n}$ a complete hyperbolic $2n$-manifold which
is homeomorphic to the total space of an $\R^n$--bundle $\xi\!: E\to B$ with the Euler number $e(\xi)$.
Then
$$
|e(\xi)|\le c(B).
$$
\end{thm}

I.~Belegradek greatly improved this result:

\begin{thm}
(I.~Belegradek \cite{Belegradek(1998a)}.)
Given a closed hyperbolic $n$-manifold $B$ ($n\ge 3$) there exists a number $C(B, k)$ so that
the number of inequivalent $\R^k$--bundles $\xi: E\to B$ whose total space admits a complete
hyperbolic metric, is at most $C(B, k)$.
\end{thm}

Given a group $\pi$, let  ${\mathcal M}_{\pi, n}$ denote the set of $n$-manifolds,
whose fundamental group is isomorphic to $\pi$ and which admit complete hyperbolic metrics.

\begin{thm}
(I.~Belegradek \cite{Belegradek(2001)}.)
Suppose that $\pi$ is a finitely-presented group with finite Betti numbers.
Assume that $\pi$ does not split as an amalgam over a virtually abelian subgroup.
The set ${\mathcal M}_{\pi, n}$ breaks into finitely many intersection preserving homotopy types.
\end{thm}

\bigskip
%\newpage
%\bibliography{/u/ma/kapovich/BIB/lit}
\bibliography{lit}

\bibliographystyle{siam}
%\addcontentsline{toc}{subsection}{References}

\end{document}